\pdfoutput=1
\RequirePackage{ifpdf}
\ifpdf 
\documentclass[pdftex]{sigma}
\else
\documentclass{sigma}
\fi

\usepackage{bm}
\usepackage{stmaryrd}
\usepackage[all]{xy}
\usepackage{tikz}
\usetikzlibrary{calc}
\usepackage{longtable}
\usepackage{multirow}

\usepackage{tabu}

\CompileMatrices

\newtheorem{thm}{Theorem}[section]

\newtheorem{coro}[thm]{Corollary}
\newtheorem{prop}[thm]{Proposition}

\theoremstyle{definition}
\newtheorem{construction}[thm]{Construction}

\newtheorem{exam}[thm]{Example}
\newtheorem{rema}[thm]{Remark}
\newtheorem{setting}[thm]{Setting}
\newtheorem{classification}[thm]{Classification}
\newtheorem{class-list}[thm]{Classification list}

\def\vector2#1#2{\left(\begin{array}{c} #1 \\ #2 \end{array}\right)}
\def\Cl{{\rm Cl}}

\def\C{{\rm C}}

\def\CC{{\mathbb C}}
\def\KK{{\mathbb K}}
\def\TT{{\mathbb T}}
\def\ZZ{{\mathbb Z}}

\def\QQ{{\mathbb Q}}
\def\PP{{\mathbb P}}

\def\Mat{{\rm Mat}}

\def\div{{\rm div}}
\def\quot{/\!\!/}

\DeclareMathOperator\HP{HP}

\def\im{{\rm im}}

\def\bangle#1{{\langle #1 \rangle}}

\def\Aut{{\rm Aut}}

\def\Spec{\operatorname{Spec}}

\def\cone{{\rm cone}}

\def\im{{\rm im}}

\numberwithin{equation}{section}

\begin{document}
\allowdisplaybreaks

\newcommand{\arXivNumber}{2108.03029}

\renewcommand{\PaperNumber}{088}

\FirstPageHeading

\ShortArticleName{On Gorenstein Fano Threefolds with an Action of a Two-Dimensional Torus}

\ArticleName{On Gorenstein Fano Threefolds\\ with an Action of a Two-Dimensional Torus}

\Author{Andreas B\"AUERLE and J\"urgen HAUSEN}

\AuthorNameForHeading{A.~B\"auerle and J.~Hausen}

\Address{Mathematisches Institut, Universit\"at T\"ubingen,\\
 Auf der Morgenstelle 10, 72076 T\"ubingen, Germany}
\Email{\href{mailto:baeuerle@math.uni-tuebingen.de}{baeuerle@math.uni-tuebingen.de}, \href{mailto:juergen.hausen@uni-tuebingen.de}{juergen.hausen@uni-tuebingen.de}}

\ArticleDates{Received April 20, 2022, in final form November 07, 2022; Published online November 16, 2022}

\Abstract{We classify the non-toric, $\mathbb Q$-factorial, Gorenstein, log terminal Fano threefolds of Picard number one that admit an effective action of a two-dimensional algebraic torus.}


\Classification{14J45; 14J35; 14L30}

\section{Introduction}

This article contributes to the classification
of singular Fano threefolds.
We work over an algebraically closed field~$\KK$
of characteristic zero.
By a Fano variety we mean a normal projective
variety~$X$ over $\KK$ admitting an ample
anticanonical divisor~$-\mathcal{K}_X$.
In the smooth case, the classifications by
Iskovskikh~\cite{Is1,Is2} and Mori--Mukai~\cite{MoMu}
provide us with a far developed picture in dimension
three.
In the singular case, the situation is less explored.
As a landmark, we have in dimension two the
classifications by Alexeev/Nikulin~\cite{AlNi}
and Nakayama~\cite{Nak} of the log terminal del Pezzo
surfaces $X$ of Gorenstein index $\imath_X \le 2$.
Here, log terminal means discrepancies greater
than $-1$ and $\imath_X$ is the smallest positive
integer with $\imath_X\mathcal{K}_X$ Cartier;
so, $\imath_X=1$ merely means that $X$ is
Gorenstein.
In dimension three, the classification problem for
singular Fano varieties is widely open.
The Mori--Fano threefolds, that means the terminal
$\QQ$-factorial Fano threefolds of Picard number one,
are intensely studied; see in particular Prokhorov's
classifications for higher index and degree
cases~\cite{Pr1,Pr2,Pr3,Pr4}.

Once we restrict to Fano varieties with many symmetries,
the singular case is more accessible.
A sample case are toric Fano varieties, where we mention
Kasprzyk's classification of the canonical toric Fano
threefolds~\cite{Ka}, comprising in particular
the toric Mori--Fano threefolds.
Going one step beyond the toric case, one
considers threefolds coming with an effective action
of a two-dimensional torus.
In this setting, the Mori--Fano threefolds have
been classified in~\cite{BeHaHuNi},
using the so-called anticanonical complex; see
also~\cite{HiWr} for generalizations.
Moreover, in~\cite{Il}, a~classification algorithm
for Gorenstein canonical Fano varieties with a torus
action of complexity one has been proposed using
the approach via polyhedral divisors~\cite{AlHa}.
However, as soon as we leave the surface case,
feasability becomes a serious question.

In the present article we classify the non-toric $\QQ$-factorial
Gorenstein log terminal Fano threefolds $X$ of Picard
number one that come with an effective action of a
two-dimensional torus.
We use the Cox ring based approach to rational
varieties with a torus action of complexity one developed
in~\cite{HaHe, HaSu}; see also
Section~\ref{sec:tvar} for a brief reminder.
The Cox ring of a normal projective variety $X$
with finitely generated divisor class group $\Cl(X)$
is defined as
\[
\mathcal{R}(X)
 =
\bigoplus_{\Cl(X)} \Gamma(X,\mathcal{O}_X(D)),
\]
where we refer to~\cite{ArDeHaLa} for the details.
For our Fano threefolds $X$ of Picard number one
acted on by a two-dimensional torus,
the divisor class group $\Cl(X)$ is of the form
$\ZZ \oplus \Gamma$ with a finite abelian torsion
part~$\Gamma$ and the Cox ring $\mathcal{R}(X)$ is a finitely
generated complete intersection ring with a~very
specific system of trinomial relations.
Moreover, the variety $X$ can be reconstructed
from the list of generator degrees in $\Cl(X)$
and the defining relations of the Cox ring
$\mathcal{R}(X)$ which allows us to encode $X$
via these Cox ring data in a compact manner.
Let us briefly describe the main result of the
article; the detailed classification lists, containing
in particular the Cox ring data,
are provided in Section~\ref{sec:class-lists}.

\begin{classification}
\label{classif:gorlt3folds}
We obtain $538$ families of non-toric, $\QQ$-factorial,
Gorenstein, log terminal Fano threefolds of Picard
number one acted on effectively by a two-dimensional
torus. Listed according to the possible divisor class
groups, we have:
\begin{table}[h!]\centering\setlength{\tabcolsep}{10.5pt}
\renewcommand{\arraystretch}{1.2}\small
\begin{tabular}{c|c|c}
\hline
Divisor class group
&Sporadic varieties
&True families
\\[1ex]\hline
$\hyperref[class-list(0)]{\ZZ}$ & \text{242} & \text{3 one-dimensional}
\\[1ex]
$\hyperref[class-list(0,2)]{\ZZ \times \ZZ/ 2 \ZZ}$ & \text{163} & \text{4 one-dimensional}
\\[1ex]
$\hyperref[class-list(0,2,2)]{\ZZ \times (\ZZ/ 2 \ZZ)^2}$ & \text{46} &
\text{5 one-dimensional,} \\ &&\text{1 two-dimensional}
\\[2ex]
$\hyperref[class-list(0,2,2,2)]{\ZZ \times (\ZZ/ 2 \ZZ)^3}$ & \text{6} & \text{1 one-dimensional}
\\[1ex]
$\hyperref[class-list(0,2,4)]{\ZZ \times \ZZ/ 2 \ZZ \times \ZZ/ 4 \ZZ}$ & \text{4} & \text{1 one-dimensional }
\\[1ex]
$\hyperref[class-list(0,2,6)]{\ZZ \times \ZZ/ 2 \ZZ \times \ZZ/ 6 \ZZ}$ & \text{1} & \text{0}
\\[1ex]
$\hyperref[class-list(0,3)]{\ZZ \times \ZZ/ 3 \ZZ}$ & \text{26} & \text{1 one-dimensional}
\\[1ex]
$\hyperref[class-list(0,3,3)]{\ZZ \times (\ZZ/ 3 \ZZ)^2}$ & \text{1} & \text{0}
\\[1ex]
$\hyperref[class-list(0,4)]{\ZZ \times \ZZ/ 4 \ZZ}$ & \text{18} & \text{1 one-dimensional}
\\[1ex]
$\hyperref[class-list(0,5)]{\ZZ \times \ZZ/ 5 \ZZ}$ & \text{4} & \text{0}
\\[1ex]
$\hyperref[class-list(0,6)]{\ZZ \times \ZZ/ 6 \ZZ}$ & \text{8} & \text{0}
\\[1ex]
$\hyperref[class-list(0,8)]{\ZZ \times \ZZ/ 8 \ZZ}$ & \text{2} & \text{0}
\\ \hline
\end{tabular}
\end{table}

\noindent
Moreover, every non-toric, $\QQ$-factorial, Gorenstein, log terminal
Fano threefold of Picard number one with an
effective action of a two-dimensional torus is isomorphic to
precisely one member of these $538$ families.
\end{classification}

Note that being Gorenstein and log terminal,
all varieties from
Classification~\ref{classif:gorlt3folds}
are canonical.
The overlap with the classification of
non-toric Mori--Fano threefolds coming with
an action of a~two-dimensisonal torus given
in~\cite{BeHaHuNi} consists
precisely of the smooth quadric in $\PP_4$.
The defining data of each of our 538
families are stored in the file~\cite{BaHa}.
Moreover, we store in this file geometric
invariants such as genus, codimension,
anticanonical self intersection, Hilbert series,
etc., which allows to extract varieties with
given properties.

Let us spend a few words on the methods used
in the classifications of non-toric Fano threefolds
$X$ of Picard number one coming with an action of
a two-dimensional torus.
As mentioned before, the main tool of~\cite{BeHaHuNi},
which settles the terminal case, is the
anticanonical complex $\mathcal{A}_X^c$
associated with $X$, a polyhedral complex
extending directly the features of the Fano
polytope from toric geometry: $X$ is
terminal if and only if $\mathcal{A}_X^c$
has only the origin as an interior lattice
point. This allows to bound the possible
Cox ring data via the volumes of
suitable lattice polytopes constructed
out of the complex.
In the log terminal Gorenstein case, one
could think of proceeding analogously by
using canonicity, which, however, appears
to end in overflowing computations, even
when building on the classification of
canonical threefold singularities with
action of a two-dimensional torus provided
in~\cite{BrHae}.
Instead we can benefit in a completely
different way and much more directly from
the Gorenstein property: it gives rise to
unit fraction identities involving
the Cox ring data that admit only
a finite number of integral solutions; see
Proposition~\ref{prop:G-matrix} for a sample.
Moreover, the computation of these integral
solutions turns out to be easily feasible,
which at the end makes the classification
possible.

\begin{rema}
The following figure shows how the 538 families from
the Classification~\ref{classif:gorlt3folds} are
distributed over the genus-codimension landscape of
Fano threefolds presented in~\cite[Figure~1]{BrKa}:\vspace{-5ex}
\begin{table}[h!]\centering
\small\setlength{\tabcolsep}{1.5pt} \setlength{\arraycolsep}{0pt}
	\begin{longtable}{c|c|ccccc|ccccc|ccccc|ccccc|ccccc|cc}\hline
	\multicolumn{29}{c}{Genus}
	\\ \hline
	\multirow{24}{*}{\rotatebox[origin=c]{90}{Codimension}}
\qquad	& &${ 2}$ &${ 3}$ &${ 4}$ &${ 5}$ &${ 6}$ & ${ 7}$ &${ 8}$ &${ 9}$ &${ 10}$ &${ 11}$ &${ 12}$ &${ 13}$ &${ 14}$ &${ 15}$ &${ 16}$ &${ 17}$ &${ 18}$ &${ 19}$ &${ 20}$ &${ 21}$ &${ 22}$ &${ 23}$ &${ 24}$ &${ 25}$ &${ 26}$ &${ 27}$ &${ 28}$
	\\
	\cline{2-29}
	&${ 1}$ & \hyperref[T(1,2)]{${ 21}$} &\hyperref[T(1,3)]{${ 5}$} & & & & & & & & & & & & & & & & & & & & & & & & &
	\\
	&${ 2}$ &\hyperref[T(2,2)]{${ 15}$} & \hyperref[T(2,3)]{${ 54}$} & \hyperref[T(2,4)]{${ 31}$} & & & & & & & & & & & & & & & & & & & & & & & &
	\\
	&${ 3}$ & & & & \hyperref[T(3,5)]{${ 22}$} & & & & & & & & & & & & & & & & & & & & & & &
	\\
	&${ 4}$ & & & \hyperref[T(4,4)]{${ 37}$} & & \hyperref[T(4,6)]{${ 17}$} & & & & & & & & & & & & & & & & & & & & & &
	\\
	&${ 5}$ & & & & & & \hyperref[T(5,7)]{${ 62}$} & & & & & & & & & & & & & & & & & & & & &
	\\
	&${ 6}$ & & & & \hyperref[T(6,5)]{${ 71}$} & & & & & & & & & & & & & & & & & & & & & & &
	\\
	&${ 7}$ & & & & & & & & \hyperref[T(7,9)]{${ 51}$} & & & & & & & & & & & & & & & & & & &
	\\
	\cline{2-29}
	&${ 8}$ & & & & & & & & & \hyperref[T(8,10)]{${ 41}$} & & & & & & & & & & & & & & & & & &
	\\
	&${ 9}$ & & & & & & & & & & \hyperref[T(9,11)]{${ 10}$} & & & & & & & & & & & & & & & & &
	\\
	&${ 10}$ & & & & & & & & & & & \hyperref[T(10,12)]{${ 1}$} & & & & & & & & & & & & & & & &
	\\
	&${ 11}$ & & & & & & & & & & & & \hyperref[T(11,13)]{${ 41}$} & & & & & & & & & & & & & & &
	\\
	&${ 12}$ & & & & & & & & & & & & & & & & & & & & & & & & & & &
	\\
	&${ 13}$ & & & & & & & & & & & & & & & & & & & & & & & & & & &
	\\
	&${ 14}$ & & & & & & & & & & & & & & & \hyperref[T(14,16)]{${ 9}$} & & & & & & & & & & & &
	\\
	\cline{2-29}
	&${ 15}$ & & & & & & & & & & & & & & & & \hyperref[T(15,17)]{${ 20}$} & & & & & & & & & & &
	\\
	&${ 16}$ & & & & & & & & & & & & & & & & & & & & & & & & & & &
	\\
	&${ 17}$ & & & & & & & & & & & & & & & & & & \hyperref[T(17,19)]{${ 12}$} & & & & & & & & &
	\\
	&${ 18}$ & & & & & & & & & & & & & & & & & & & & & & & & & & &
	\\
	&${ 19}$ & & & & & & & & & & & & & & & & & & & & \hyperref[T(19,21)]{${ 4}$} & & & & & & &
	\\
	&${ 20}$ & & & & & & & & & & & & & & & & & & & & & \hyperref[T(20,22)]{${ 3}$} & & & & & &
	\\
	&${ 21}$ & & & & & & & & & & & & & & & & & & & & & & & & & & &
	\\
	\cline{2-29}
	&${ 22}$ & & & & & & & & & & & & & & & & & & & & & & & & & & &
	\\
	&${ 23}$ & & & & & & & & & & & & & & & & & & & & & & & & \hyperref[T(23,25)]{${ 4}$} & & &
	\\
	&${ 24}$ & & & & & & & & & & & & & & & & & & & & & & & & & \hyperref[T(24,26)]{${ 4}$} & &
	\\
	&${ 25}$ & & & & & & & & & & & & & & & & & & & & & & & & & & &
	\\
	&${ 26}$ & & & & & & & & & & & & & & & & & & & & & & & & & & &\hyperref[T(26,28)]{${ 3}$}
	\\
	\hline
 \end{longtable}
\end{table}

\noindent
Here the genus of a Fano threefold~$X$ is $h^0(X,-\mathcal{K}_X)-2$
and the codimension is taken with respect to embedding into
a weighted projective space by means of a minimal system of
homogeneous generators of the anticanonical ring
\[
A_X = \bigoplus_{n \in \ZZ_{\ge 0}} \Gamma(X,-n\mathcal{K}_X).
\]
\end{rema}

The article is organized as follows.
Section~\ref{sec:tvar} serves to provide
the necessary background on the approach
to rational projective varieties $X$
with a torus action of complexity one
via the Cox ring
based on~\cite{HaHe,HaSu}.
In Picard number one, this approach
represents any family of $\QQ$-factorial
varieties~$X$ in terms of an integral
matrix~$P$.
Very first constraints arise from
log terminality:
Proposition~\ref{rem:piconeconfigs},
originally due to~\cite{BeHaHuNi},
shows that log terminality leaves us
with eight types of matrices~$P$
to consider.
In Section~\ref{sec:sample}, we exemplarily
discuss one of these eight cases,
showing basically all the necessary
arguments.
In particular, we see how to establish
in several refining steps appropriate
bounds on the entries of the defining
matrix~$P$ making a computational
treatment feasible.
The full elaboration of all cases will
be presented elsewhere.
Section~\ref{sec:class-lists} presents
the classification
tables and in Section~\ref{sec:hilbert-series},
we compute the Hilbert--Poincar\'{e} series of our
varieties (also accessible via clicking on
the items in the genus-codimension
landscape).

\section{Torus actions of complexity one}
\label{sec:tvar}

We recall the necessary background on rational
varieties with a torus action of complexity one
and fix our notation.
The reader is assumed to be familiar with the
very basics of toric geometry, in particular
the correspondence between fans and toric
varieties; see~\cite{CoLiSc,Dan,Ful}.
We restrict ourselves to spending just a few
words on Cox's quotient presentation~\cite{Cox}
of a toric variety arising from a fan.

\begin{construction}
\label{constr:coxtoric}
Let $Z$ be the toric variety defined by
a fan $\Sigma$ in a lattice $N$ such that
the primitive generators $v_1,\dots, v_r$
of the rays of $\Sigma$ span the
rational vector space $N_\QQ = N \otimes_\ZZ \QQ$.
We have a linear map
\[
P \colon \quad \ZZ^r \to N,\qquad
e_i \mapsto v_i.
\]
In case $N = \ZZ^n$, we also speak
of the \emph{generator matrix}
$P = [v_1, \dots, v_r]$ of $\Sigma$.
The divisor class group and the Cox ring
of $Z$ are
\[
\Cl(Z) = K := \ZZ^r / \im(P^*),
\qquad
\mathcal{R}(Z) = \KK[T_1,\dots,T_r],
\qquad \deg(T_i) = Q(e_i),
\]
where $P^*$ denotes the dual map of $P$ and
$Q \colon \ZZ^r \to K$ the projection.
Now, one defines a fan $\hat \Sigma$ in $\ZZ^r$
consisting of faces of the positive orthant
of $\QQ^r$ by
\[
\hat \Sigma
 :=
\{\delta_0 \preceq \QQ_{\ge 0}^r; \,
P(\delta_0) \subseteq \sigma \text{ for some } \sigma \in \Sigma\}.
\]
The toric variety $\hat Z$ associated with $\hat \Sigma$
is an open toric subset in $\bar Z := \KK^r$.
As $P$ is a map of the fans~$\hat \Sigma$ and
$\Sigma$, it defines a toric morphism
$p \colon \hat Z \to Z$.
The quasitorus
\[
H = \Spec \KK[K] = \ker(p) \subseteq \TT^r = (\KK^*)^r
\]
acts as a subgroup of the torus $\TT^r$ on $\hat Z$ and the
morphism $p \colon \hat Z \to Z$ turns out to
be a good quotient with respect to the $H$-action.
\end{construction}

The quotient presentation of toric varieties is
a central piece in the Cox ring based approach
to rational varieties with a torus action of
complexity one provided by~\cite{HaHe,HaSu};
see also~\cite[Section~3.4]{ArDeHaLa}.
The first step, however, is the following purely
algebraic construction of a certain class of graded
algebras; see~\cite[Construction~3.4.2.1]{ArDeHaLa}
and more generally~\cite[Constructions~3.5 and~6.3]{HaHiWr}.

\begin{construction}
\label{constr:RAP}
Fix $r \in \ZZ_{\ge 1}$, a sequence
$n_0, \dots, n_r \in \ZZ_{\ge 1}$, set
$n := n_0 + \dots + n_r$, and fix
integers $m \in \ZZ_{\ge 0}$ and $0 < s < n+m-r$.
The input data are matrices
\[
A =[a_0, \dots, a_r] \in\Mat(2,r+1;\KK),\qquad
P =\begin{bmatrix}
L & 0
\\
d & d'
\end{bmatrix}
 \in
\Mat(r+s,n+m; \ZZ),
\]
where $A$ has pairwise linearly independent
columns and $P$ is built from an
$(s \times n)$-block $d$, an $(s \times m)$-block
$d'$ and an $(r \times n)$-block $L$
of the shape
\[
L = \begin{bmatrix}
-l_0 & l_1 & \dots & 0
\\
\vdots & \vdots & \ddots & \vdots
\\
-l_0 & 0 &\dots & l_{r}
\end{bmatrix}\!,\qquad
l_i = (l_{i1}, \dots, l_{in_i}) \in \ZZ_{\ge 1}^{n_i}
\]
such that the columns $v_{ij}$, $v_k$ of
$P$ are pairwise different primitive vectors
generating $\QQ^{r+s}$ as a~cone.
Consider the polynomial algebra
\[
\KK[T_{ij},S_k] :=
\KK[T_{ij},S_k; \; 0 \le i \le r, \, 1 \le j \le n_i, 1 \le k \le m].
\]
Denote by $\mathfrak{I}$ the set of
all triples $I = (i_1,i_2,i_3)$ with
$0 \le i_1 < i_2 < i_3 \le r$
and define for any $I \in \mathfrak{I}$
a trinomial
\[
g_I := g_{i_1,i_2,i_3} := \det
\begin{bmatrix}
T_{i_1}^{l_{i_1}} & T_{i_2}^{l_{i_2}} & T_{i_3}^{l_{i_3}}
\\
a_{i_1} & a_{i_2} & a_{i_3}
\end{bmatrix}\!,\qquad
T_i^{l_i} := T_{i1}^{l_{i1}} \cdots T_{in_i}^{l_{in_i}}.
\]
Consider the factor group
$K := \ZZ^{n+m}/\rm{im}(P^*)$
and the projection $Q \colon \ZZ^{n+m} \to K$.
We define a~$K$-grading on
$\KK[T_{ij},S_k]$ by setting
\[
\deg(T_{ij}) :=\omega_{ij} :=Q(e_{ij}),\qquad
\deg(S_{k}) :=\omega_k :=Q(e_{k}).
\]
Then the trinomials $g_I$ just introduced
are $K$-homogeneous and they all share the
same $K$-degree.
In particular, we obtain a $K$-graded factor
algebra
\[
R(A,P) :=\KK[T_{ij},S_k] / \bangle{g_I; \ I \in \mathfrak{I}}.
\]
\end{construction}

\begin{exam}
\label{ex:running-ex-1}
We choose $r = 2$, moreover $n_0 = 2$, $n_1=n_2 = 1$ and
$m=1$ and, finally $s=2$. In this setting, consider
the defining matrices
\[A :=
\begin{bmatrix}
0 & 1 & 1
\\
1 & 1 & 0
\end{bmatrix}\!,
\qquad
P :=\begin{bmatrix}
-1 & -1 & $\phantom{$-$}$4 & 0 & 0
\\
-1 & -1 & $\phantom{$-$}$0 & 2 & 0
\\
$\phantom{$-$}$0 & $\phantom{$-$}$0 & -3 & 1 & 1
\\
$\phantom{$-$}$0 & -2 & $\phantom{$-$}$4 & 0 & 0
\end{bmatrix}\!.
\]
The algebra $R(A,P)$ arising from these matrices
comes due to $r=2$ with a single trinomial
relation and is explicitly given by
\[
R(A,P)
 =
\KK[T_{01},T_{02},T_{11},T_{21},S_1] / \big\langle T_{01}T_{02} + T_{11}^4 + T_{21}^2\big\rangle.
\]
We have
$K = \ZZ^5/{\rm im}(P^*) = \ZZ \oplus \ZZ/2\ZZ \oplus \ZZ / 2\ZZ$
and the degrees of the $T_{ij}$ and $S_1$ are
the columns of the \emph{degree matrix}
\[
Q =\begin{bmatrix}
2 & 2 & 1 & 2 & 1
\\
\bar 1 & \bar 1 & \bar 0 & \bar 0 & \bar 0
\\
\bar 1 & \bar 1 & \bar 1 & \bar 1 & \bar 0
\end{bmatrix}\!.
\]
\end{exam}

\begin{thm}[{see~\cite[Theorem~3.4.2.3]{ArDeHaLa}, also~\cite[Theorems~3.10 and 6.5]{HaHiWr}}]
The ring $R(A,P)$ produced by
Construction~{\rm \ref{constr:RAP}}
is a normal complete intersection ring
and its ideal of relations is generated
by the trinomials $g_i = g_{i,i+1,i+2}$,
where $i = 0, \dots, r-2$.
\end{thm}

\begin{rema}
We call a defining matrix $P$ \emph{irredundant} if we have
$l_{i1}n_i \ge 2$ for $i = 0,\dots, r$.
Each $R(A,P)$ is isomorphic as a graded algebra to
some $R(A',P')$ with $P'$ irredundant.
Note that for $r \ge 2$ and an irredundant $P$,
the ring $R(A,P)$ is not a polynomial ring.
\end{rema}

\begin{rema}
Consider a defining matrix $P$ as in Construction~\ref{constr:RAP}.
By an \emph{admissible operation} on the matrix $P$ we mean
one of the following:
\begin{enumerate}\itemsep=0pt
\item[$(i)$]
adding a multiple of one of the upper $r$ rows
to one of the lower $s$ rows,
\item[$(ii)$]
applying a unimodular matrix from the left to the $(d,d')$
block,
\item[$(iii)$]
swapping two columns $v_{ij_1}$ and $v_{ij_2}$
inside a leaf $v_{i1}, \dots, v_{in_i}$,
\item[$(iv)$]
swapping two leafs $v_{i1}, \dots, v_{in_i}$ and
$v_{j1}, \dots, v_{jn_j}$ and rearranging the
$L$-block by elementary row operations into its
required shape,
\item[$(v)$]
swapping two columns $v_{k_1}$ and $v_{k_2}$ of the
$d'$-block.
\end{enumerate}
If $P'$ arises from $P$ via admissible operations,
then with a suitable $A'$, the graded rings
$R(A,P)$ and $R(A',P')$ are isomorphic.
\end{rema}

\begin{rema}
The matrix $A$ of a ring $R(A,P)$ is responsible
for the coefficients of the defining trinomials
$g_i = g_{i,i+1,i+2}$. By rescaling variables
we can always reduce to defining relations
of the shape
\[
T_0^{l_0} + T_1^{l_1} + T_2^{l_2},\quad
\lambda_1T_1^{l_1} + T_2^{l_2} + T_3^{l_3},\quad \dots, \quad
\lambda_{r-2}T_{r-2}^{l_{r-2}} + T_{r-1}^{l_{r-1}} + T_r^{l_r}
\]
with pairwise distinct $1 \ne \lambda_i \in \CC^*$.
In particular, in case of a single defining relation,
there~is no need to care about the coefficients.
The matrix $A$ is motivated by the geometry behind
$R(A,P)$, see Remark~\ref{rem:divisorsonX}.
\end{rema}

We enter the second step, producing
rational normal varieties~$X$ with torus action
$\TT^s \times X \to X$ of complexity one.
Each of the resulting $X$ comes embedded in a toric
variety~$Z$, is defined in homogeneous coordinates
by the above trinomials $g_0, \dots, g_{r-1}$ and
the torus~$\TT^s$ acting on $X$ is a subtorus of the
acting torus $\TT^{r+s}$ of~$Z$.
The original references are again~\cite{HaHe, HaSu};
see also \cite[Construction~3.4.3.6]{ArDeHaLa}
as well as the more
general~\cite[Constructions~3.5 and~6.13]{HaHiWr}.

\begin{construction}
\label{constr:RAPdown}
In the situation of Construction~\ref{constr:RAP},
assume $r \ge 2$ and that $P$ is irredundant.
Consider the common zero set of the
defining relations of $R(A,P)$:
\[
\bar{X}
 :=
V(g_I; \ I \in \mathfrak{I})
 \subseteq
\bar{Z}
 :=
\KK^{n+m}.
\]
Let $\Sigma$ be any fan in $N = \ZZ^{r+s}$
having the columns of $P$ as the primitive
generators of its rays.
Then $\hat X := \bar X \cap \hat Z$ and
Construction~\ref{constr:coxtoric} yield
a commutative diagram
\[
\xymatrix@R=1.5em{
{\bar{X}}
\ar@{}[r]|\subseteq
\ar@{}[d]|{\rotatebox[origin=c]{90}{$\scriptstyle\subseteq$}}
&
{\bar{Z}}
\ar@{}[d]|{\rotatebox[origin=c]{90}{$\scriptstyle\subseteq$}}
\\
{\hat{X}}
\ar[r]
\ar[d]_{\quot H}^p
&
{\hat{Z}}
\ar[d]^{\quot H}_p
\\
X
\ar[r]
&
Z,}
\]
where $X := X(A,P,\Sigma) := p(\hat X)$ is a non-toric,
closed subvariety
of the toric variety~$Z$ arising from~$\Sigma$.
Dimension, divisor class group
and Cox ring of $X$ are
\[
\dim(X) = s+1,\qquad
\Cl(X) \cong K,\qquad
\mathcal{R}(X) \cong R(A,P).
\]
The subtorus $\TT^s \subseteq \TT^{r+s}$ of the
acting torus of $Z$ associated with the
sublattice $\ZZ^{s} \subseteq \ZZ^{r+s}$
leaves~$X$ invariant and the induced $T$-action
on $X$ is of complexity one.
\end{construction}

\begin{exam}
\label{ex:running-ex-2}
We continue Example~\ref{ex:running-ex-1}.
Let $\Sigma$ be the fan in $\ZZ^4$ having
$P$ as its generator matrix and
the maximal cones
\begin{alignat*}{4}
& \cone(v_{02},v_{11},v_{21},v_1),
\qquad &&
\cone(v_{01},v_{11},v_{21},v_1),
\qquad&&
\cone(v_{01},v_{02},v_{21},v_1),&
\\
& \cone(v_{01},v_{02},v_{11},v_1),
\qquad &&
\cone(v_{01},v_{02},v_{11},v_{21}).&&&
\end{alignat*}
The associated toric variety $Z$ is a
four-dimensional fake weighted projective
space with divisor class group
\[
\Cl(Z) = K = \ZZ \oplus \ZZ/2\ZZ \oplus \ZZ / 2\ZZ.
\]
Moreover, $H = \KK^* \times \{\pm 1\} \times \{\pm 1\}$
acts on $\bar Z = \KK^5$ via the weights given by
the columns of the degree matrix $Q$ and
Construction~\ref{constr:RAPdown} becomes
\[
\xymatrix@R=1.5em{
V\big(T_{01}T_{02} + T_{11}^4 + T_{21}^2\big)
\ar@{}[r]|(.75){=}
\ar@{}[d]|{\rotatebox[origin=c]{90}{$\scriptstyle\subseteq$}}
&
{\bar{X}}
\ar@{}[r]|\subseteq
\ar@{}[d]|{\rotatebox[origin=c]{90}{$\scriptstyle\subseteq$}}
&
{\bar{Z}}
\ar@{}[d]|{\rotatebox[origin=c]{90}{$\scriptstyle\subseteq$}}
\ar@{}[r]|{=}
&
{\KK^5}
\ar@{}[d]|{\rotatebox[origin=c]{90}{$\scriptstyle\subseteq$}}
\\
{\bar{X} \setminus \{0\}}
\ar@{}[r]|(.75){=}
&
{\hat{X}}
\ar[r]
\ar[d]_{\quot H}^p
&
{\hat{Z}}
\ar[d]^{\quot H}_p
\ar@{}[r]|{=}
&
{\KK^5 \setminus \{0\}}.
\\
&
X
\ar[r]
&
Z
&
}
\]
\end{exam}

\begin{thm}[{see~\cite[Theorem~4.4.1.6]{ArDeHaLa}
and~\cite[Theorems~3.10 and~6.18]{HaHiWr}}]
Every non-toric rational normal projective variety
with a torus action of complexity one is equivariantly
isomorphic to some $X(A,P,\Sigma)$ arising
from Construction~$\ref{constr:RAP}$.
\end{thm}

Any variety $X = X(A,P,\Sigma)$ inherits many
geometric properties from its ambient toric
variety~$Z$.
A first observation concerns the restriction
of the invariant divisors from $Z$ to $X$;
see~\mbox{\cite[Proposition~3.2.4.5]{ArDeHaLa}}.

\begin{rema}
Consider $X = X(A,P,\Sigma)$ as in
Construction~\ref{constr:RAPdown}.
The columns~$v_{ij}$ and $v_k$ of~$P$
define prime divisors
$D_{ij} = V_Z(T_{ij})$ and $D_k = V_Z(T_k)$
on $Z$.
The restrictions
$D^X_{ij} = V_X(T_{ij})$ and $D^X_k = V_X(S_k)$
are prime divisors on $X$ and in
$\Cl(Z) = K = \Cl(X)$, we have
\[
[D_{ij}] = \deg(T_{ij}) = \bigl[D_{ij}^X\bigr],\qquad
[D_k] = \deg(T_k) = \bigl[D_k^X\bigr].
\]
\end{rema}

We recover the divisors $D_{ij}^X$ as the components of
the critical values $c_i \in \PP_1$ of a certain quotient
map; see~\cite[Proposition~3.16]{HaHiWr} for a general treatment.

\begin{rema}
\label{rem:divisorsonX}
Consider $X = X(A,P,\Sigma)$ as in
Construction~\ref{constr:RAPdown}.
Consider the open sets of points
having finite isotropy groups with
respect to the $\TT^s$-action:
\[
Z_0 = \{z \in Z; \ \TT^s_z \text{ is finite} \},\qquad
X_0 = X \cap Z_0 = \{x \in X; \ \TT^s_x \text{ is finite} \}.
\]
Then $Z_0 \subseteq Z$ is invariant under the
torus $\TT^{r+s}$ and $X_0 \subseteq X$ is invariant
under $\TT^s$.
Moreover, we have a commutative diagram
\[
\xymatrix@R=1.5em{
{X}
\ar@{}[r]|\subseteq
\ar@{}[d]|{\rotatebox[origin=c]{90}{$\scriptstyle\subseteq$}}
&
{Z}
\ar@{}[d]|{\rotatebox[origin=c]{90}{$\scriptstyle\subseteq$}}
\\
{X_0}
\ar[r]
\ar[d]_{\quot {\TT^s}}^{\pi_X}
&
{Z_0}
\ar[d]^{\quot {\TT^s}}_{\pi_Z}
\\
{\PP_1}
\ar[r]
&
{\PP_r},}
\]
where $\pi_X$ and $\pi_Z$ are categorical quotients
with respect to the actions of $\TT^s$ on $X$ and $Z$
respectively and $\pi_Z$ is a toric morphism.
Moreover, we obtain
\[
\overline{\pi_X^{-1}(c_i)} = \bigcup_{j=1}^{n_i} D_{ij}^X \subseteq X,
\qquad
\overline{\pi_Z^{-1} (C_i)} = \bigcup_{j=1}^{n_i} D_{ij} \subseteq Z
\]
with the toric divisors
$C_0, \dots, \C_r \subseteq \PP_r$
and the points $c_i \in \PP_1$ having the
$i$-th column of $A$ as its homogeneous
coordinates. Finally,
\[
\vert \TT^s_{x_{ij}} \vert = l_{ij}
\]
holds for the order of the isotropy group $\TT^s_{x_{ij}}$ of the
action of the torus $\TT^s$ at any general point $x_{ij} \in D_{ij}^X$.
\end{rema}

The divisors from Remark~\ref{rem:divisorsonX}
also allow an explicit presentation of an anticanonical
divisor; see~\cite[Proposition~3.4.4.1]{ArDeHaLa}.

\begin{rema}
\label{rem:anticanclass}
Let $X = X(A,P,\Sigma)$ arise from
Construction~\ref{constr:RAPdown}.
Then the anticanonical divisor
class of $X$ is given as
\[
- \mathcal{K}_X =
\sum_{i,j} \deg(T_{ij}) + \sum_k \deg(S_k)-
(r-1) \sum_{i=1}^{n_0} l_{0j} \deg(T_{0j})
 \in K =\Cl(X).
\]
In particular, due to $\deg(T_{ij}) = [D^X_{ij}]$
and $\deg(T_k) = [D^X_k]$, we have the
following anticanonical divisor on $X$:
\[
D^X_0
 :=
\sum_{i,j} D^X_{ij} + \sum_k D^X_k
- (r-1) \sum_{j=1}^{n_0} l_{0 j} D^X_{0 j}.
\]
\end{rema}

\begin{exam}
For the variety $X$ from Example~\ref{ex:running-ex-2},
we can compute the anticanonical class~as
\begin{gather*}
- \mathcal{K}_X =
\deg(T_{11})+\deg(T_{21})+\deg(S_1) =
(4,\bar 0, \bar 0) \in
\ZZ \oplus \ZZ/2\ZZ \oplus \ZZ/2\ZZ =\Cl(X).
\end{gather*}
In particular, we see that the anticanonical class
is ample and, consequently, $X$ is a Fano variety.
\end{exam}

\begin{rema}
Let $X = X(A,P,\Sigma)$ arise from
Construction~\ref{constr:RAPdown}.
We call $\sigma \in \Sigma$ an \emph{$X$-cone}
if the corresponding toric
orbit $\TT^{r+s} \cdot z_\sigma \subseteq Z$
meets $X \subseteq Z$.
A cone $\sigma \in \Sigma$ is an $X$-cone
if and only if one of the following holds:
\begin{enumerate}\itemsep=0pt
\item[$(i)$]
$\sigma$ is a \emph{big cone}, that means
$v_{0j_0}, \dots, v_{rj_r} \in \sigma$
for some $j_0, \dots, j_r$,
\item[$(ii)$]
$\sigma$ is a \emph{leaf cone}, that means
$\sigma \subseteq \cone(v_{i1}, \dots, v_{in_i},v_1,\dots,v_m)$
for some $i$.
\end{enumerate}
Every $X$-cone $\sigma \in \Sigma$ defines an affine
open subvariety $X_\sigma =: X \cap Z_\sigma$
in $X$ by cutting down the corresponding affine toric
chart $Z_\sigma \subseteq Z$.
Note that $X$ is covered by the $X_\sigma$, where $\sigma$ runs
through the $X$-cones of $\Sigma$.
\end{rema}

\begin{exam}
\label{ex:running-ex-3}
Consider again the variety $X$
from Example~\ref{ex:running-ex-2}.
Then the fan $\Sigma$ has exactly four
maximal $X$-cones, namely
\begin{alignat*}{4}
&\cone(v_{02},v_{11},v_{21},v_1),
\qquad&&
\cone(v_{01},v_{11},v_{21},v_1),
\\
&\cone(v_{01},v_{02},v_{11},v_{21}),
\qquad&&
\cone(v_{01},v_{02},v_1).
\end{alignat*}
The first three are big cones, whereas the
fourth one is a leaf cone.
Thus, $X$ is covered by four open
affine subvarieties, given by the maximal
$X$-cones of $\Sigma$.
\end{exam}

Let us see how to detect Cartier divisors,
that means locally principal Weil divisors,
on a variety $X = X(A,P,\Sigma)$ in terms of
the defining data.

\begin{prop}
\label{prop:gorchar}
Let $X = X(A,P,\Sigma)$ arise from
Construction~$\ref{constr:RAPdown}$.
Consider on~$Z$ and~$X$ the Weil
divisors
\[
D =\sum a_{ij} D_{ij} + \sum a_k D_k,\qquad
D^X =\sum a_{ij} D^X_{ij} + \sum a_k D^X_k,
\]
in $K = \Cl(Z) = \Cl(X)$ the classes $\omega = [D] = [D^X]$,
$\omega_{ij} = [D_{ij}] = [D_{ij}^X]$,
$\omega_k = [D_k] = [D_k^X]$
and an $X$-cone $\sigma \in \Sigma$.
Then the following statements are equivalent:
\begin{enumerate}\itemsep=0pt
\item[$(i)$]
The divisor $D^X$ is Cartier on $X_\sigma$.
\item[$(ii)$]
The divisor $D$ is Cartier on $Z_\sigma$.
\item[$(iii)$]
We have $D = \div(\chi^u)$ on $Z_\sigma$ with
a character $\chi^u$ of $\TT^{r+s}$.
\item[$(iv)$]
There is $u \in \ZZ^{r+s}$ with
$\bangle{u,v_{ij}} = a_{ij}$ and
$\bangle{u,v_k} = a_k$ for all
$v_{ij}, v_k \in \sigma$.
\item[$(v)$]
We have
$\omega \in \bangle{\omega_{ij}, \omega_k; \ v_{ij}, v_k \not\in \sigma}$
in $K = \Cl(X)$.
\end{enumerate}
In particular, $D$ is a Cartier divisor on $X$
if and only if one of these conditions holds
for all maximal $X$-cones $\sigma \in \Sigma$.
\end{prop}

\begin{proof}
The equivalence of~$(i)$, $(ii)$ and~$(v)$ follows from
Proposition~\cite[Proposition~3.3.1.5]{ArDeHaLa}.
The rest is basic toric geometry.
\end{proof}

A normal variety $X$ is \emph{$\QQ$-factorial} if every
Weil divisor $D$ on $X$ admits a Cartier
multiple~$nD$ with $n \in \ZZ_{\ge 1}$.

\begin{coro}
A variety $X = X(A,P,\Sigma)$ as in
Construction~$\ref{constr:RAPdown}$
is $\QQ$-factorial if and only if
each $X$-cone $\sigma \in \Sigma$
is simplicial.
\end{coro}

Now, recall that a variety is \emph{Gorenstein} if its
canonical class is Cartier.
Combining Remark~\ref{rem:anticanclass} and
Proposition~\ref{prop:gorchar}, we obtain the
following characterization.

\begin{coro}
\label{cor:gor2div}
Consider $X = X(A,P,\Sigma)$ and let
$D^X = \sum a_{ij} D^X_{ij} + \sum a_k D^X_k$
be an anticanonical divisor on $X$.
Then $X$ is Gorenstein if and
only if for every maximal $X$-cone
$\sigma$, there is a linear form $u \in \ZZ^{r+s}$
with
\[
\bangle{u,v_{ij}} = a_{ij},\qquad
\bangle{u,v_k} = a_k\qquad
\text{for all}\quad v_{ij}, v_k \in \sigma.
\]
\end{coro}

\begin{exam}
\label{ex:running-ex-4}
Consider again the variety $X$
from Example~\ref{ex:running-ex-2}
and the four maximal $X$-cones
given in Example~\ref{ex:running-ex-3}.
Listed accordingly, we have linear forms
\[
(2,0,1,-1),
\quad
(0,0,1,1),
\quad
(-2,2,-3,0),
\quad
(-1,1,1,0)
\]
representing the anticanonical divisor $D^X_0$
on the corresponding affine open subvarieties of $X$.
In particular, $X$ is Gorenstein.
\end{exam}

If $X$ is a $\QQ$-factorial Fano variety of Picard
number one, then, the divisor class group $\Cl(X)$
allows a positive splitting into a free cyclic part and its
torsion part~$\Gamma$, that means that we have
an isomorphism
\[
\Cl(X) \cong \ZZ \oplus \Gamma
\]
such that for the anticanonical class
$\omega_X = (w_X,\eta_X)$, the $\ZZ$-part
$w_X$ is positive.
Note that in this setting the $\ZZ$-part
of any divisor class $\omega = (w,\eta)$
does not depend on the particular choice of the
splitting.

\begin{coro}
\label{prop:gor2div2}
Let $X = X(A,P,\Sigma)$ be $\QQ$-factorial, Gorenstein, Fano
and of Picard number one.
Then, for every maximal $X$-cone $\sigma$,
the $\ZZ$-parts $w_{ij}$, $w_k$ of the generator degrees
and $w_X$ of the anticanonical class satisfy
\[
\gcd(w_{ij},w_k; \ v_{ij} \not\in \sigma, v_k \not\in \sigma)
\mid
w_X.
\]
\end{coro}

\begin{proof}
As $X$ is Gorenstein, the canonical class $\omega_X$
represents a Cartier divisor.
Proposition~\ref{prop:gorchar} tells us that for
every maximal $X$-cone $\sigma$,
the $\omega_X$ lies in the subgroup of $\Cl(X) =K$
generated by the classes
$\omega_{ij}$, $\omega_k$, where $v_{ij} \not\in \sigma$,
$v_k \not\in \sigma$.
Thus, the $\ZZ$-part $w_X$ lies in the ideal of $\ZZ$
generated by the $\ZZ$-parts $w_{ij}$, $w_k$,
where $v_{ij} \not\in \sigma$, $v_k \not\in \sigma$.
The assertion follows.
\end{proof}

Finally, we discuss log terminality.
Recall that given any resolution of
singularities $\pi \colon X' \to X$ of a normal
variety, we have the ramification
formula
\[
\mathcal{K}_{X'} - \pi^*\mathcal{K}_{X}
 =
\sum_{i=1}^r a_i E_i
\]
with canonical divisors on $X'$
and $X$ and the exceptional divisors
$E_i, \dots, E_r$.
Then~$X$ is called \emph{log terminal} if
we have $a_i > -1$ for $i = 1, \dots, r$.

We characterize log terminality of a given
$\QQ$-factorial Fano variety $X = X(A,P,\Sigma)$.
A \emph{platonic tuple} is a tuple
$(l_0,\dots,l_r)$ of positive integers
such that after re-ordering
the $l_i$ decreasingly, we obtain a tuple
$(a,b,c,1, \dots, 1)$ with $(a,b,c)$ one of
\[
(x,y,1),
\quad
(y,2,2),
\quad
(5,3,2),
\quad
(4,3,2),
\quad
(3,3,2).
\]

\begin{prop}[{see~\cite[Theorem~3.13]{ArBrHaWr}}]
\label{rem:platonictuples}
A $\QQ$-factorial Fano variety $X = X(A,P,\Sigma)$
has at most log terminal
singularities if and only if
for any $X$-cone
$\sigma = \cone(v_{0j_0}, \dots, v_{rj_r})$
the exponents $l_{0j_0}, \dots, l_{rj_r}$
form a platonic tuple.
\end{prop}

\begin{exam}
For the variety $X = X(A,P,\Sigma)$
from Example~\ref{ex:running-ex-2}
we have to consider the $X$-cones
\[
\cone(v_{02},v_{11},v_{21}),
\qquad
\cone(v_{01},v_{11},v_{21}).
\]
Both of them yield the exponent
tuple $(1,4,2)$ which is platonic.
Consequently, $X$ is log terminal.
\end{exam}

Log terminality leads to the following first
constraints on the defining matrix $P$ of our
Fano varieties $X = X(A,P,\Sigma)$.

\begin{prop}[{see~\cite[Lemma~5.2]{BeHaHuNi}}]
\label{rem:piconeconfigs}
Let $X = X(A,P,\Sigma)$ a non-toric,
$\QQ$-factorial, log terminal Fano threefold
of Picard number one, where $P$ is irredundant.
Then, after suitable admissible operations,
$P$ fits into one of the following cases:
\begin{enumerate}\itemsep=0pt
\item[$(i)$]
 $m=0$, $r=2$ and $n=5$, where $n_0 = n_1 = 2$, $n_2 = 1$,
\item[$(ii)$]
 $m=0$, $r=3$ and $n=6$, where $n_0 = n_1 = 2$, $n_2 = n_3 = 1$,
\item[$(iii)$]
 $m=0$, $r=4$ and $n=7$, where $n_0 = n_1 = 2$, $n_2 = n_3 = n_4 = 1$,
\item[$(iv)$]
 $m=0$, $r=2$ and $n=5$, where $n_0 = 3$, $n_1 = n_2 = 1$,
\item[$(v)$]
 $m=0$, $r=3$ and $n=6$, where $n_0 = 3$, $n_1 = n_2 = n_3 = 1$,
\item[$(vi)$]
$m=1$, $r=2$ and $n=4$, where $n_0 = 2$, $n_1 = n_2 = 1$,
\item[$(vii)$]
$m=1$, $r=3$ and $n=5$, where $n_0 = 2$, $n_1 = n_2 = n_3 = 1$,
\item[$(viii)$]
$m=2$, $r=2$ and $n=3$, where $n_0 = n_1 = n_2 = 1$.
\end{enumerate}
\end{prop}

\begin{rema}
Every rational Gorenstein del Pezzo surface has
at most canonical singularities and thus is in
particular log terminal; see~\cite{HiWa}.
For Gorenstein Fano varieties of higher dimension
even the latter property need not hold. For instance,
\[
P =\begin{bmatrix}
-3 & -1 & 3 & 1 & $\phantom{$-$}$0
\\
-3 & -1 & 0 & 0 & $\phantom{$-$}$k
\\
-4 & -1 & 1 & 0 & $\phantom{$-$}$k
\\
 $\phantom{$-$}$1 & $\phantom{$-$}$0 & 0 & 0 & -1
\end{bmatrix}
\]
defines for each $k \ge 4$ a $\QQ$-factorial
Fano threefold $X = X(A,P,\Sigma)$ of Picard number
one, which is not log terminal by
Proposition~\ref{rem:platonictuples}.
More explicitly, $\Sigma$ consists of all
pointed cones generated by columns of $P$
and we have
\begin{gather*}
X =
V\big(T_{01}^3T_{02} + T_{11}^3T_{12} + T_{21}^k\big)
 \subseteq\PP_{1,\, k-3, \, 1, \, k-3, \, 1},
\\[1ex]
[ - \mathcal{K}_X ] =k - 3 \in\ZZ =\Cl(X).
\end{gather*}
This example series shows moreover that the
Gorenstein and Fano condition together are
even in the specific setting of threefolds
with an action of a two-dimensional torus
not enough to guarantee finiteness in fixed
dimension and Picard number.
\end{rema}

\section{Elaboration of a sample case}
\label{sec:sample}

Proposition~\ref{rem:piconeconfigs} divides
the proof of the classification theorem
into cases~(i) to~(viii).
Here we give a treatment of case~(vi)
as a sample, showing all types of
arguments of the full proof.
The setting is the following.

\begin{setting}
\label{setting:r=2-n0=2-n1=1-n2=1-m=1}
Let $X = X(A,P,\Sigma)$ be a $\QQ$-factorial
Fano threefold of Picard number one with
$r=2$, $n_0=2$, $n_1=n_2=1$ and $m=1$.
Then
\begin{equation*}
P =[v_{01},v_{02},v_{11},v_{21},v_1] =
\begin{bmatrix}
-l_{01} & -l_{02} & l_{11} & 0 & 0
\\
-l_{01} & -l_{02} & 0 & l_{21} & 0
\\
d_{011} & d_{021} & d_{111} & d_{211} & d_{11}
\\
d_{012} & d_{022} & d_{112} & d_{212} & d_{12}
\end{bmatrix}
\end{equation*}
with pairwise different primitive columns
$v_{01}$, $v_{02}$, $v_{11}$, $v_{21}$ and $v_1$
generating $\QQ^4$ as a cone.
The maximal $X$-cones of the fan $\Sigma$ of $Z$ are given by
\begin{gather*}
\sigma_{01} = \cone(v_{02},v_{11},v_{21},v_1),\qquad
\sigma_{02} = \cone(v_{01},v_{11},v_{21},v_1),
\\
\sigma_{1} = \cone(v_{01},v_{02},v_{11},v_{21}),\qquad
\tau_{0} = \cone(v_{01},v_{02},v_1).
\end{gather*}
We have $K = \ZZ \oplus \Gamma$ with the torsion part
$\Gamma$ and denote $\deg(T_{ij}) = (w_{ij},\eta_{ij})$ as well
as $\deg(T_{k}) = (w_{k},\eta_{k})$ accordingly.
In particular, we write
\begin{equation*}
Q^0
 =
[w_{01},w_{02},w_{11},w_{21},w_1]
\end{equation*}
for the $\ZZ$-part of the degree matrix~$Q$. Note that the
vector~$(w_{01},w_{02},w_{11},w_{21},w_1)$ is primitive
in $\ZZ^5$ and generates $\ker(P)$.
\end{setting}

Our first series of constraints arising
from the log terminality and the Gorenstein
property directly aims for entries of
the defining matrix~$P$.

\begin{prop}\label{prop:sample-1}\samepage
Consider $X = X(A,P,\Sigma)$ as in
Setting~$\ref{setting:r=2-n0=2-n1=1-n2=1-m=1}$.
Assume that $X$ is non-toric, log terminal
and Gorenstein.
\begin{enumerate}\itemsep=0pt
\item[$(i)$]
$D_0^X = (1-l_{01})D_{01}^X + (1-l_{02})D_{02}^X + D_{11}^X + D_{21}^X + D_1^X$
is an anticanonical divisor on $X$.
In~particular, the free part of the anticanonical divisor class
is given by
\[
w_X = (1-l_{01})w_{01} + (1-l_{02})w_{02} + w_{11} + w_{21} + w_{1}.
\]
\item[$(ii)$]
Admissible column and row operations turn
the defining matrix $P$ into the shape
\begin{gather*}
P
 =
\begin{bmatrix}
-l_{01} & -l_{02} & l_{11} & 0 & 0
\\
-l_{01} & -l_{02} & 0 & l_{21} & 0
\\
1-l_{01} & 1-l_{02} & d_{111} & d_{211} & 1
\\
d_{012} & d_{022} & d_{112} & d_{212} & 0
\end{bmatrix},
\\[1ex]
l_{01} \ge l_{02} \ge 1, \qquad
 l_{11} \ge l_{21} \ge 2,
\\[1ex]
 0 \le d_{211}, d_{212} < l_{21},\qquad
 0 \le d_{012} < l_{01},
\\[1ex]
 -\frac{w_X}{w_{02}} \le d_{022} < l_{02} + \frac{w_X}{w_{02}},
\end{gather*}
where $w_{02} \mid w_X$ and the tuple of exponents
$(l_{01},l_{11},l_{21})$
fits into precisely one of the following constellations:
\begin{alignat*}{5}
& (1,x,y), \qquad &&x \ge y > 1; \qquad && (2,z,3), \qquad &&3 \le z \le 5;&
\\
&(y,2,2), \qquad &&y \ge 2; \qquad && (3,z,2), \qquad &&3 \le z \le 5;&
\\
&(2,y,2), \qquad &&y \ge 3; \qquad && (z,3,2), \qquad && 4 \le z \le 5.&
\end{alignat*}
\end{enumerate}
\end{prop}

\begin{proof}
For the first assertion, note that we have
$r = 2$ and that the defining relation of
the Cox ring is given as
\[
g = T_{01}^{l_{01}}T_{02}^{l_{02}} + T_{11}^{l_{11}} + T_{21}^{l_{21}} .
\]
Thus, $\mu_X = \deg(g) = l_{01}\deg(T_{01}) + l_{02}\deg(T_{02})$
and Remark~\ref{rem:anticanclass} shows that
the anticanonical divisor $D^X_0$ is as claimed.

We prove~$(ii)$.
Swapping, if necessary, the first two
columns of $P$, we achieve
$l_{01} \ge l_{02}$.
Similarly, exchanging the data of
$v_{11}$ and $v_{21}$, we ensure
$l_{11} \ge l_{21}$.
Since~$X$ is non-toric, we must have
$l_{21} > 1$.
As we assume $X$ to be log terminal, we
can apply Remark~\ref{rem:platonictuples}
to $\cone(v_{01},v_{11},v_{21})$, showing that
$(l_{01},l_{11},l_{21})$ is as in the assertion.

We care about the entries of the $(d,d')$ block of~$P$.
Since $v_1 \in \ZZ^4$ is primitive, we can
apply a suitable unimodular $2 \times 2$ matrix from
the left to the $(d,d')$ block to ensure
\[
d_{11} = 1,\qquad
d_{12} = 0.
\]
We begin to make use of the assumption that
$X$ is Gorenstein.
First consider the $X$-cone $\tau_0 = \cone(v_{01},v_{02},v_1)$.
Then Corollary~\ref{cor:gor2div} provides a linear form
$u \in \ZZ^4$ such that
\[
\bangle{u,v_{01}} = 1-l_{01},\qquad
\bangle{u,v_{02}} = 1-l_{02},\qquad
\bangle{u,v_{1}} = 1.
\]
The last equation tells us in particular $u_3 = 1$.
Plugging this into the first two equations yields
\[
d_{011} = l_{01}(u_1+u_2) - u_4 d_{012} + 1-l_{01},\qquad
d_{021} = l_{02}(u_1+u_2) - u_4 d_{022} + 1-l_{02}.
\]
Thus, adding the $(u_1+u_2)$-fold of the first
and the $u_4$-fold of the fourth row of $P$
to the third one, we obtain
\[
d_{011} = 1-l_{01},\qquad
d_{021} = 1-l_{02}.
\]
Moreover, adding an appropriate multiple of the first row
of $P$ to the fourth one, we achieve
\[
0 \le d_{012} < l_{01}.
\]
Now consider the maximal $X$-cone
$\sigma_{02} = \cone(v_{01},v_{11},v_{21},v_1)$.
Let $u \in \ZZ^4$ be a linear form
representing $D_0^X$ on $X_{\sigma_{02}}$ according to
Corollary~\ref{cor:gor2div}($iii$).
Then
\[
0=Q^0 \cdot P^* \cdot u =
\sum \bangle{u,v_{ij}}w_{ij} + \bangle{u,v_1}w_1
 =w_X + (u_4 d_{022} - l_{02}(u_1 + u_2)) w_{02}.
\]
In particular, we see that $w_{02}$ divides $w_X$.
Moreover, we must have $u_3=1$. We obtain
\[
1-l_{01} =\bangle{u,v_{01}} =
-l_{01}u_1 - l_{01}u_2 + 1-l_{01} + u_4d_{012}.
\]
This merely means $l_{01}(u_1 + u_2)=u_4d_{012}$.
Plugging this into the previous equation yields
\[
- l_{01} \frac{w_X}{w_{02}} =
u_4 (d_{022}l_{01}-d_{012}l_{02}).
\]
Thus, $(d_{022}l_{01}-d_{012}l_{02})$ divides $l_{01}\frac{w_X}{w_{02}}$.
As a consequence, we can estimate $d_{022}$ as follows:
\[
l_{02}\frac{d_{012}}{l_{01}} - \frac{w_X}{w_{02}}
 \le
d_{022}
 \le
l_{02}\frac{d_{012}}{l_{01}} + \frac{w_X}{w_{02}}.
\]
Combining this with $0 \leq d_{012} < l_{01}$, we arrive
at the desired bounds for $d_{022}$.
Finally, we achieve
\begin{gather*}
0 \le d_{211}, d_{212} < l_{21}
\end{gather*}
by adding suitable multiples of the difference of the first
two rows of $P$ to third and the fourth one.
\end{proof}

The second series of constraints shows that all entries of
the $\ZZ$-part of the degree matrix
$Q^0 = [w_{01},w_{02},w_{11},w_{21},w_1]$ are bounded.

\begin{prop}
\label{prop:G-matrix}
Let $X = X(A,P,\Sigma)$ be as in
Setting~$\ref{setting:r=2-n0=2-n1=1-n2=1-m=1}$.
Assume that $X$ is non-toric, log terminal
and Gorenstein.
\begin{enumerate}\itemsep=0pt
\item[$(i)$]
Let $\alpha_{01}$, $\alpha_{02}$ and $\beta_1$ be any three positive
integers and consider the $5 \times 5$ matrix
\[
G :=
\begin{bmatrix}
1-l_{01}-\alpha_{01} & 1-l_{02} & 1 & 1 & 1
\\
1-l_{01} &1-l_{02}-\alpha_{02} & 1 & 1 & 1
\\
1-l_{01} &1-l_{02} & 1 & 1 & 1-\beta_1
\\
-l_{01} & -l_{02} & l_{11} & 0 & 0
\\
-l_{01} & -l_{02} & 0 & l_{21} & 0
\end{bmatrix}\!.
\]
Then $G$ is of rank at least four.
Moreover, $\det(G) = 0$ if and only
if $\alpha_{01}$, $\alpha_{02}$, $\beta_1$ and
$l_{01}$, $l_{02}$, $l_{11}$, $l_{21}$ satisfy the identity
\[
\frac{1}{\beta_1}+\frac{1}{\alpha_{01}}+\frac{1}{\alpha_{02}}
+
\biggl( \frac{l_{01}}{\alpha_{01}} + \frac{l_{02}}{\alpha_{02}} \biggr)
\biggl( \frac{1}{l_{11}} + \frac{1}{l_{21}} - 1 \biggr) =1.
\]
\item[$(ii)$]
There are unique $\alpha_{01}, \alpha_{02}, \beta_1 \in \ZZ_{\ge 1}$
with $\alpha_{01} w_{01} = \alpha_{02} w_{02} = \beta_1 w_1 = w_X$
and the corresponding matrix $G$ from~$(i)$ satisfies
\[
\ker(G) =\ker(P) =\ZZ \cdot (w_{01},w_{02},w_{11},w_{21},w_1).
\]
\item[$(iii)$]
According to the possible constellations of
$(l_{01},l_{11},l_{02})$ from~Proposition~$\ref{prop:sample-1}(ii)$
we have the following bounds for the
entries of the matrix $G$:
$$\text{\centering\setlength{\tabcolsep}{10.5pt}
\renewcommand{\arraystretch}{1.4}\rm\small\tabulinesep=1mm\begin{tabu}{c|c}
\hline
$(l_{01},l_{11},l_{21})$& bounds for
$(l_{01},l_{02},l_{11},l_{21},\alpha_{01},\alpha_{02},\beta_1)$
\\
\hline
$(1,x,y)$&\qquad\qquad $(1,1,30,8,35,4,36)$
\\
$(y,2,2)$&\qquad\qquad $(11,6,2,2,6,6,6)$
\\
$(2,y,2)$&\qquad\qquad $(2,2,18,2,28,10,30)$
\\
$(2,z,3)$&\qquad\qquad $(2,2,4,3,4,4,6)$
\\
$(3,z,2)$&\qquad\qquad $(3,3,5,2,12,5,24)$
\\
$(z,3,2)$&\qquad\qquad $(5,4,3,2,14,5,18)$
\\
\hline
\end{tabu}}$$
\end{enumerate}
\end{prop}

\begin{proof}
We verify~$(i)$. In order to see that $G$ is of rank at least
four, we just compute the minor
\[
G_{3,1} =\det
\begin{bmatrix}
1-l_{02} & 1 & 1 & 1
\\
1-l_{02}-\alpha_{02} & 1 & 1 & 1
\\
-l_{02} & l_{11} & 0 & 0
\\
-l_{02} & 0 & l_{21} & 0
\end{bmatrix}
 =\alpha_{02} l_{11}l_{21} \ne 0.
\]
Moreover, suitably rearranging the equation
$\det(G) = 0$, we arrive at the displayed identity on
$\alpha_{01}$, $\alpha_{02}$, $\beta_1$ and
$l_{01}$, $l_{02}$, $l_{11}$, $l_{21}$.

We show~$(ii)$.
Applying Corollary~\ref{prop:gor2div2} to the
three maximal $X$-cones $\sigma_{01}$, $\sigma_{02}$
and $\sigma_1$ we see that each of $w_{01}$, $w_{02}$ and $w_1$
is a multiple of $w_X$ and hence we obtain
positive integers $\alpha_{01}$, $\alpha_{02}$ and $\beta_1$
with
\begin{equation*}
\alpha_{01} w_{01}
 =
\alpha_{02} w_{02}
 =
\beta_1 w_1
 =
(1-l_{01}) w_{01} + (1-l_{02}) w_{02} + w_{11} + w_{21} + w_1.
\end{equation*}
Moreover, by homogeneity of the defining relation
$g$ we have
\begin{equation*}
l_{01} w_{01} + l_{02}w_{02} =l_{11} w_{11} =l_{21} w_{21}.
\end{equation*}
Now, $G$ from~$(i)$ is the coefficient matrix of
the corresponding system of linear equations.
In particular, for any choice of
$\alpha_{01}$, $\alpha_{02}$ and $\beta_1$
the integral matrix $G$ has kernel generated by
the primitive vector
$(w_{01},w_{02},w_{11},w_{21},w_1) \in \ZZ^5$.

We turn to~$(iii)$. We exemplarily treat the
configuration $(l_{01},l_{11},l_{21}) = (1,x,y)$,
where $x \ge y \ge 2$.
In this situation the condition $\det(G) = 0$
reads as
\begin{equation*}
\frac{1}{\beta_1}+
\bigg( \frac{1}{\alpha_{01}} + \frac{1}{\alpha_{02}} \bigg)
\bigg( \frac{1}{x} + \frac{1}{y} \bigg) =1.
\end{equation*}
We may assume $\alpha_{01} \ge \alpha_{02}$. We directly see $\beta_1 \ge 2$.
We arrive at 16 possible choices for
$\alpha_{02}$ and $y$ given by the conditions
\[
1 \le \alpha_{02} \le 4,\qquad
2 \le y \le 8,\qquad
\alpha_{02} \le 2\quad \text{or}\quad y = 2,
\]
We establish preliminary bounds.
First observe that the equation $\det(G) = 0$ can
be resolved for each of the variables.
For instance, we have
\[
x
 =
\frac{(\alpha_{01} + \alpha_{02})\beta_1 y}{\alpha_{01} \alpha_{02} \beta_1 y - \alpha_{01} \alpha_{02} y - \alpha_{01} \beta_1 - \alpha_{02} \beta_1} .
\]
Plugging the above 16 choices for $a = \alpha_{02}$
and $b = \beta_1$ into this expression, we obtain more
explicit presentations for the possible values of~$x$:
\begin{alignat*}{5}
&\frac{2 + 2a}{\big(1 - \frac{2}{b}\big) a - 1},\qquad&&
\frac{4 + 2a}{\big(3 - \frac{4}{b}\big) a - 2},\qquad&&
\frac{6 + 2a}{\big(5 - \frac{6}{b}\big) a - 3},\qquad&&
\frac{8 + 2a}{\big(7 - \frac{8}{b}\big) a - 4},&
\\
&\frac{3 + 3a}{\big(2 - \frac{3}{b}\big) a - 1},\qquad&&
\frac{6 + 3a}{\big(5 - \frac{6}{b}\big) a - 2},\qquad&&
\frac{4 + 4a}{\big(3 - \frac{4}{b}\big) a - 1},\qquad&&
\frac{8 + 4a}{\big(7 - \frac{8}{b}\big) a - 2},&
\\
&\frac{5 + 5a}{\big(4 - \frac{5}{b}\big) a - 1},\qquad&&
\frac{10 + 5a}{\big(9 - \frac{10}{b}\big) a - 2},\qquad&&
\frac{6 + 6a}{\big(5 - \frac{6}{b}\big) a - 1},\qquad&&
\frac{ 12 + 6a}{\big(11 - \frac{12}{b}\big) a - 2},&
\\
&\frac{14 + 7a}{\big(6 - \frac{7}{b}\big) a - 1},\qquad&&
\frac{14 + 7 a}{\big(13 - \frac{14}{b}\big) a - 2},\qquad&&
\frac{8 + 8 a}{\big(7 - \frac{8}{b}\big) a - 1},\qquad&&
\frac{16 + 8 a}{\big(15 - \frac{16}{b}\big) a - 2}.&
\end{alignat*}
Each of these expressions can be maximized,
having in mind $a,b,x \in \ZZ_{\ge 1}$.
We arrive at $x \le 48$, attained by the first
expression with $b=3$ and $a=4$.

In a similar way, we obtain preliminary bounds for
$\alpha_{01}$ and $\beta_1$.
Now, running $\det(G)=0$ with the preliminary
bounds for the involved variables, we arrive at
the bounds from the assertion.
\end{proof}

Besides constraints on the defining data,
we also need criteria to decide computationally
whether or not given defining data lead to
isomorphic varieties.
For this, we say that a defining matrix $P$ as in
Construction~\ref{constr:RAPdown}
has \emph{ordered exponents} if
we have
\begin{enumerate}\itemsep=0pt
 \item[$(i)$] $n_0 \ge \dots \ge n_r$,
 \item[$(ii)$] $l_{i1} \ge \dots \ge l_{in_i}$ for each $i = 0, \dots, r$ and
 \item[$(iii)$] if $n_i = n_{i+1}$ then $l_{i1} \ge l_{i+1, 1}$.
\end{enumerate}
If $P$ has ordered exponents,
then we call the data
$(n_0, \dots, n_r,m)$
the \emph{format} of~$P$.
Note that via admissible operations,
we can always assume that $P$ has ordered
exponents.

\begin{prop}
\label{prop:non-iso-crit}
Let $(A,P,\Sigma)$ and $(A',P',\Sigma')$
be as in Construction~$\ref{constr:RAPdown}$
such that the associated varieties
$X$ and $X'$ are isomorphic to each other.
\begin{enumerate}\itemsep=0pt
\item[$(i)$]
There is an isomorphism $\varphi \colon X \to X'$
which is equivariant with respect to the torus
actions.
\item[$(ii)$]
If $P$ and $P'$ have ordered exponents, then
they share the same format and for each $i$
there is an $i'$ with $n_{i'} = n_i$ and
$(l_{i1}, \dots, l_{in_i}) = (l_{i'1}', \dots, l_{i'n_i}')$
such that
\[
\bangle{\deg(T_{ij}); \, j = 1, \dots, n_i}
 \cong
\bangle{\deg(T_{i'j}); \, j = 1, \dots, n_i}
\]
holds for the subgroups in $\Cl(X)$
and $\Cl(X')$, respectively, generated by the
corresponding degrees.
\end{enumerate}
\end{prop}

\begin{proof}
For the first assertion, observe that for
any isomorphism $\varphi \colon X \to X'$
of varieties, we can install a torus action
on $X'$ making $\varphi$ equivariant.
Now, any torus action of complexity one on
the non-toric $X'$ corresponds to a maximal
torus in the affine algebraic group $\Aut(X')$;
see for instance~\cite[Theorem~2.1]{ArHaHeLi}.
Thus, the assertion follows from the fact that
any two maximal tori in an affine algebraic
group are conjugate.
The second assertion follows from the first one
and the fact that any equivariant isomorphism
respects the data described in
Remark~\ref{rem:divisorsonX}.
\end{proof}

This allows us in particular to associate
with any $X \cong X(A,P,\Sigma)$ a \emph{format}
$(n_0, \dots, n_r,m)$ by taking $P$ with
ordered exponents.
Bringing together our constraints on the entries of $P$,
the bounds on the entries of the matrix $G$ and the
distinction criterion just shown allows us to perform the
classification computationally.

\begin{coro}
There is a list of $155$ explicitly given
matrices $P$ of format $(2,1,1,1)$, each of them defining
a non-toric $\QQ$-factorial,
Gorenstein, log terminal Fano threefold $X(A,P,\Sigma)$
of Picard number one.
$$ \text{\centering\setlength{\tabcolsep}{4.5pt}\rm\small
\renewcommand{\arraystretch}{1.2}\tabulinesep=1mm\begin{tabu}{c|ccccccccc|c}
\multicolumn{11}{c}{Number of members $P$ of the list according
to divisor class group and exponent configuration.}
\\
\hline
&$\ZZ$&$\ZZ+\ZZ_2$&$\ZZ+\ZZ_3$&$\ZZ+\ZZ_4$&$\ZZ+\ZZ_5$&$\ZZ+\ZZ_6$&$\ZZ+\ZZ_2^2$&$\ZZ+\ZZ_2+\ZZ_4$
&$\ZZ+\ZZ_2^3$&sum
\\
\hline
$(1,x,y)$&13&19&6&4&1&3&4&2&&52
\\
$(y,2,2)$&&12&&3&&&9&&1&25
\\
$(2,y,2)$&13&24&&&&&8&&1&46
\\
$(2,z,3)$&2&&3&&&1&&&&6
\\
$(3,z,2)$&7&6&&&&&1&&&14
\\
$(z,3,2)$&5&6&&&&&1&&&12
\\
\hline
sum&40&67&9&7&1&4&23&2&2&155
\\ \hline
\end{tabu}}$$

\noindent
Distinct matrices from the list
yield non-isomorphic varieties
and every non-toric, $\QQ$-factorial, log terminal
and Gorenstein Fano threefold of Picard number one
of format $(2,1,1,1)$ is isomorphic to an $X = X(A,P,\Sigma)$
with $P$ from the list.
\end{coro}

\begin{proof}
Proposition~\ref{prop:G-matrix} allows us
to write down explicitly all possible
matrices $G$ and hence to determine all
possible
$Q^0 = [w_{01},w_{02},w_{11},w_{21},w_1]$
by computer.
Now, recall that $P$ annihilates the
transpose of $Q^0$.
This enables us to determine in the
matrix $P$, adjusted according to
Proposition~\ref{prop:sample-1},
all the remaining variables.
So, we are left with a finite list of
explicitly given possible defining
matrices~$P$.
Checking for the necessary properties
by means of~\cite{HaKe} and reducing via
Proposition~\ref{prop:non-iso-crit}
to data defining pairwise non-isomorphic
varieties,
we obtain the list presented in
the assertion.
\end{proof}

The remaining cases from Proposition~\ref{rem:piconeconfigs}
are treated analogously to our sample case, using basically
the same arguments in the details. We give a summarizing
overview on the strategy. The explicit elaboration will be
presented elsewhere.

\begin{rema}
Each of the cases~$(i)$ to $(v)$ and $(vii)$, $(viii)$ from
Proposition~\ref{rem:piconeconfigs} can be
treated according to the following pattern.
\begin{itemize}
\item
Establish the analogue of
Setting~\ref{setting:r=2-n0=2-n1=1-n2=1-m=1}
and Proposition~\ref{prop:sample-1}.
Observe that each row of the resulting
matrix $P$ admits at most one entry,
which is not bounded by other entries of~$P$.
\item
With each matrix~$P$ obtained so far
associate a matrix $G$ as in
Proposition~\ref{prop:G-matrix} by following
the lines of the proof of Proposition~\ref{prop:G-matrix}$(ii)$.
Bound the entries of~$G$, using that $G$ is not
of full rank as in the proof of Proposition~\ref{prop:G-matrix}$(iii)$.
\item
Go through the possible values of $G$,
compute $Q^0 = \ker(G) = \ker(P)$ explicitly
and use this identity to determine the remaining
unbounded values of $P$, one per row, as mentioned.
\item
From the resulting list of explicit matrices $P$ remove
those not defining a Gorenstein Fano variety
and remove redundancies using
Proposition~\ref{prop:non-iso-crit}.
\end{itemize}
\end{rema}

In the present classification, we comfortably
get along by establishing bounds for the integral
solutions of the unit fraction identities such as
$\det(G) = 0$ in Proposition~\ref{prop:G-matrix}
and naively going through the possibilities.
For more extensive projects, for instance, higher
Gorenstein indices or higher dimension, it
turns out that establishing algorithms for
the computation of specific unit fraction decompositions
is considerably more efficient; see~\cite{Ba,HaHaHaSp}.

\newpage

\section{Classification lists}
\label{sec:class-lists}

Here we provide the detailed presentation
of our classification result.
Let us briefly recall the background.
Each non-toric, $\QQ$-factorial, Gorenstein,
log terminal Fano threefold $X$ of Picard number
one coming with an effective action of a
two-dimensional torus is uniquely determined
by its Cox ring.
In particular, $X$ can be encoded by the degree
matrix $Q$, that means the list of degrees
of the Cox ring generators in $\Cl(X)$
and the defining trinomial relations
$g_0,\dots,g_{r-1}$.
For instance, our example variety~$X$
from~Examples~\ref{ex:running-ex-1},
\ref{ex:running-ex-2},
\ref{ex:running-ex-3}
and~\ref{ex:running-ex-4} is encoded by
\[
Q =
\begin{bmatrix}
2 & 2 & 1 & 2 & 1
\\
\bar 1 & \bar 1 & \bar 0 & \bar 0 & \bar 0
\\
\bar 1 & \bar 1 & \bar 1 & \bar 1 & \bar 0
\end{bmatrix}\!,
\qquad
g_0 = T_1T_2 + T_3^4 + T_4^2 ,
\]
where the columns of $Q$ live in
$\ZZ \oplus \ZZ/2\ZZ \oplus \ZZ / 2\ZZ$.
Indeed, the defining matrix $P$ is determined
up to admissible operations by $Q$, the
format $(2,1,1,1)$ and the list of exponents
of $g_0$.
Alternatively, $X$ is the hypersurface defined
by $g_0$ in the fake weighted projective
space $Z = \hat Z / H$, where
$\hat Z = \KK^5 \setminus \{0\}$
and the quasitorus $H$ and its action on $\KK^5$
are given by
\[
H = \KK^* \times \{\pm 1\} \times \{\pm 1\},\qquad
(t,\zeta,\eta) \cdot z =
\big(t^2\zeta\eta z_1, t^2\zeta\eta z_2, t\eta z_3,t^2\eta z_4, t z_5\big).
\]

We turn to the classification lists.
Every non-toric, $\QQ$-factorial, Gorenstein,
log terminal Fano threefold $X$ of Picard number
one coming with an effective action of a
two-dimensional torus is isomorphic to
precisely one of the listed varieties.
Conversely, each of the listed data defines a~non-toric, $\QQ$-factorial, Gorenstein,
log terminal Fano threefold $X$ of Picard number
one coming with an effective action of a
two-dimensional torus.

Each of the lists represents a possible
divisor class group and format.
Each variety in such a~list is specified
by its matrix $Q$ of generator degrees and
its defining trinomial relations; observe
that we number the variables of the relation
conventionally and not in accordance with
the double-indexed enumeration of the columns
of the associated defining matrix~$P$.
Besides the specifying data, we list the
anticanonical self intersection.
A file containing also the defining
matrices~$P$ and further invariants is
available at~\cite{BaHa}.


\begin{class-list}\label{class-list(0)}
Non-toric, $\QQ$-factorial, Gorenstein, log terminal Fano
threefolds of Picard number one with an effective two-torus action:
Specifying data for divisor class group $\ZZ$ and format $(2,2,1,0)$.
\begin{table}[h!]\centering\setlength{\tabcolsep}{3.5pt}
\renewcommand{\arraystretch}{1.3}\small\tabulinesep=1mm
\begin{tabu}{c|c|c|c||c|c|c|c}
\hline
ID&relations&gd-matrix&$-\mathcal{K}^3$&ID&relations&gd-matrix&$-\mathcal{K}^3$
\\
\hline
$10$&$T_{1}^{5}T_{2}+T_{3}^{3}T_{4}^{3}+T_{5}^{2}$
&$\arraycolsep=1mm\begin{bmatrix}1&1&1&1&3\end{bmatrix}$&${2}$
&${32}$&${T_{1}^{4}T_{2}+T_{3}^{2}T_{4}^{2}+T_{5}^{3}}$
&$\arraycolsep=1mm\begin{bmatrix}{1} & {2} & {2} & {1} & {2}\end{bmatrix}$&${6}$
\\
\hline
${11}$&${T_{1}^{4}T_{2}^{3}+T_{3}^{3}T_{4}^{2}+T_{5}^{2}}$
&$\arraycolsep=1mm\begin{bmatrix}{1} & {2} & {2} & {2} & {5}\end{bmatrix}$&${2}$
&${33}$&${T_{1}^{5}T_{2}+T_{3}^{2}T_{4}+T_{5}^{3}}$
&$\arraycolsep=1mm\begin{bmatrix}{1} & {1} & {2} & {2} & {2}\end{bmatrix}$&${6}$
\\
\hline
${14}$&${T_{1}^{5}T_{2}+T_{3}^{3}T_{4}^{2}+T_{5}^{2}}$
&$\arraycolsep=1mm\begin{bmatrix}{2} & {4} & {4} & {1} & {7}\end{bmatrix}$&${4}$
&${34}$&${T_{1}^{4}T_{2}^{2}+T_{3}^{2}T_{4}+T_{5}^{3}}$
&$\arraycolsep=1mm\begin{bmatrix}{1} & {1} & {2} & {2} & {2}\end{bmatrix}$&${6}$
\\
\hline
${15}$&${T_{1}^{3}T_{2}^{2}+T_{3}^{3}T_{4}^{2}+T_{5}^{2}}$
&$\arraycolsep=1mm\begin{bmatrix}{2} & {4} & {4} & {1} & {7}\end{bmatrix}$&${4}$
&${35}$&${T_{1}^{3}T_{2}+T_{3}^{2}T_{4}^{2}+T_{5}^{3}}$
&$\arraycolsep=1mm\begin{bmatrix}{2} & {6} & {3} & {3} & {4}\end{bmatrix}$&${6}$
\\
\hline
${16}$&${T_{1}^{3}T_{2}^{2}+T_{3}^{3}T_{4}^{2}+T_{5}^{2}}$
&$\arraycolsep=1mm\begin{bmatrix}{2} & {1} & {2} & {1} & {4}\end{bmatrix}$&${4}$
&${36}$&${T_{1}^{3}T_{2}^{2}+T_{3}^{2}T_{4}+T_{5}^{3}}$
&$\arraycolsep=1mm\begin{bmatrix}{1} & {6} & {6} & {3} & {5}\end{bmatrix}$&${6}$
\\
\hline
${17}$&${T_{1}^{5}T_{2}^{3}+T_{3}^{3}T_{4}+T_{5}^{2}}$
&$\arraycolsep=1mm\begin{bmatrix}{1} & {1} & {2} & {2} & {4}\end{bmatrix}$&${4}$
&${37}$&${T_{1}^{3}T_{2}^{2}+T_{3}^{2}T_{4}+T_{5}^{3}}$
&$\arraycolsep=1mm\begin{bmatrix}{2} & {3} & {3} & {6} & {4}\end{bmatrix}$&${6}$
\\
\hline
${18}$&${T_{1}^{3}T_{2}+T_{3}^{3}T_{4}^{2}+T_{5}^{2}}$
&$\arraycolsep=1mm\begin{bmatrix}{4} & {2} & {4} & {1} & {7}\end{bmatrix}$&${4}$
&${38}$&${T_{1}^{16}T_{2}+T_{3}^{2}T_{4}+T_{5}^{2}}$
&$\arraycolsep=1mm\begin{bmatrix}{1} & {2} & {6} & {6} & {9}\end{bmatrix}$&${6}$
\\
\hline
${31}$&${T_{1}^{3}T_{2}+T_{3}^{3}T_{4}+T_{5}^{2}}$
&$\arraycolsep=1mm\begin{bmatrix}{3} & {1} & {3} & {1} & {5}\end{bmatrix}$&${6}$
&${39}$&${T_{1}^{12}T_{2}^{3}+T_{3}^{2}T_{4}+T_{5}^{2}}$
&$\arraycolsep=1mm\begin{bmatrix}{1} & {2} & {6} & {6} & {9}\end{bmatrix}$&${6}$
\\
\hline
\end{tabu}
\end{table}

\begin{table}[t!]\centering\setlength{\tabcolsep}{3.5pt}
\renewcommand{\arraystretch}{1.25}\small\tabulinesep=1mm
\begin{tabu}{c|c|c|c||c|c|c|c}
\hline
ID&relations&gd-matrix&$-\mathcal{K}^3$&ID&relations&gd-matrix&$-\mathcal{K}^3$
\\
\hline
${40}$&${T_{1}^{8}T_{2}^{5}+T_{3}^{2}T_{4}+T_{5}^{2}}$
&$\arraycolsep=1mm\begin{bmatrix}{1} & {2} & {6} & {6} & {9}\end{bmatrix}$&${6}$
&${88}$&${T_{1}^{2}T_{2}+T_{3}^{2}T_{4}+T_{5}^{2}}$
&$\arraycolsep=1mm\begin{bmatrix}{3} & {2} & {3} & {2} & {4}\end{bmatrix}$&${12}$
\\
\hline
${41}$&${T_{1}^{7}T_{2}^{4}+T_{3}^{2}T_{4}+T_{5}^{2}}$
&$\arraycolsep=1mm\begin{bmatrix}{2} & {1} & {6} & {6} & {9}\end{bmatrix}$&${6}$
&${89}$&${T_{1}^{21}T_{2}+T_{3}T_{4}+T_{5}^{3}}$
&$\arraycolsep=1mm\begin{bmatrix}{1} & {3} & {12} & {12} & {8}\end{bmatrix}$&${12}$
\\
\hline
${42}$&${T_{1}^{4}T_{2}+T_{3}^{2}T_{4}+T_{5}^{2}}$
&$\arraycolsep=1mm\begin{bmatrix}{3} & {2} & {6} & {2} & {7}\end{bmatrix}$&${6}$
&${90}$&${T_{1}^{18}T_{2}^{2}+T_{3}T_{4}+T_{5}^{3}}$
&$\arraycolsep=1mm\begin{bmatrix}{1} & {3} & {12} & {12} & {8}\end{bmatrix}$&${12}$
\\
\hline
${52}$&${T_{1}^{4}T_{2}+T_{3}^{3}T_{4}^{3}+T_{5}^{2}}$
&$\arraycolsep=1mm\begin{bmatrix}{1} & {2} & {1} & {1} & {3}\end{bmatrix}$&${8}$
&${91}$&${T_{1}^{12}T_{2}^{4}+T_{3}T_{4}+T_{5}^{3}}$
&$\arraycolsep=1mm\begin{bmatrix}{1} & {3} & {12} & {12} & {8}\end{bmatrix}$&${12}$
\\
\hline
${53}$&${T_{1}^{4}T_{2}^{3}+T_{3}^{3}T_{4}+T_{5}^{2}}$
&$\arraycolsep=1mm\begin{bmatrix}{1} & {2} & {2} & {4} & {5}\end{bmatrix}$&${8}$
&${92}$&${T_{1}^{10}T_{2}+T_{3}T_{4}+T_{5}^{4}}$
&$\arraycolsep=1mm\begin{bmatrix}{1} & {2} & {6} & {6} & {3}\end{bmatrix}$&${12}$
\\
\hline
${54}$&${T_{1}^{11}T_{2}+T_{3}^{2}T_{4}+T_{5}^{2}}$
&$\arraycolsep=1mm\begin{bmatrix}{1} & {1} & {4} & {4} & {6}\end{bmatrix}$&${8}$
&${93}$&${T_{1}^{9}T_{2}^{5}+T_{3}T_{4}+T_{5}^{3}}$
&$\arraycolsep=1mm\begin{bmatrix}{1} & {3} & {12} & {12} & {8}\end{bmatrix}$&${12}$
\\
\hline
${55}$&${T_{1}^{9}T_{2}^{3}+T_{3}^{2}T_{4}+T_{5}^{2}}$
&$\arraycolsep=1mm\begin{bmatrix}{1} & {1} & {4} & {4} & {6}\end{bmatrix}$&${8}$
&${94}$&${T_{1}^{9}T_{2}+T_{3}T_{4}+T_{5}^{6}}$
&$\arraycolsep=1mm\begin{bmatrix}{1} & {3} & {6} & {6} & {2}\end{bmatrix}$&${12}$
\\
\hline
${56}$&${T_{1}^{8}T_{2}+T_{3}^{2}T_{4}+T_{5}^{2}}$
&$\arraycolsep=1mm\begin{bmatrix}{1} & {2} & {4} & {2} & {5}\end{bmatrix}$&${8}$
&${95}$&${T_{1}^{7}T_{2}^{3}+T_{3}T_{4}+T_{5}^{3}}$
&$\arraycolsep=1mm\begin{bmatrix}{3} & {1} & {12} & {12} & {8}\end{bmatrix}$&${12}$
\\
\hline
${57}$&${T_{1}^{7}T_{2}^{5}+T_{3}^{2}T_{4}+T_{5}^{2}}$
&$\arraycolsep=1mm\begin{bmatrix}{1} & {1} & {4} & {4} & {6}\end{bmatrix}$&${8}$
&${96}$&${T_{1}^{6}T_{2}^{3}+T_{3}T_{4}+T_{5}^{4}}$
&$\arraycolsep=1mm\begin{bmatrix}{1} & {2} & {6} & {6} & {3}\end{bmatrix}$&${12}$
\\
\hline
${58}$&${T_{1}^{4}T_{2}^{3}+T_{3}^{2}T_{4}+T_{5}^{2}}$
&$\arraycolsep=1mm\begin{bmatrix}{1} & {2} & {4} & {2} & {5}\end{bmatrix}$&${8}$
&${97}$&${T_{1}^{5}T_{2}^{2}+T_{3}T_{4}+T_{5}^{4}}$
&$\arraycolsep=1mm\begin{bmatrix}{2} & {1} & {6} & {6} & {3}\end{bmatrix}$&${12}$
\\
\hline
${65}$&${T_{1}^{3}T_{2}+T_{3}^{2}T_{4}+T_{5}^{5}}$
&$\arraycolsep=1mm\begin{bmatrix}{1} & {2} & {2} & {1} & {1}\end{bmatrix}$&${10}$
&${98}$&${T_{1}^{3}T_{2}^{2}+T_{3}T_{4}+T_{5}^{12}}$
&$\arraycolsep=1mm\begin{bmatrix}{2} & {3} & {6} & {6} & {1}\end{bmatrix}$&${12}$
\\
\hline
${66}$&${T_{1}^{2}T_{2}+T_{3}^{2}T_{4}+T_{5}^{5}}$
&$\arraycolsep=1mm\begin{bmatrix}{2} & {1} & {2} & {1} & {1}\end{bmatrix}$&${10}$
&${99}$&${T_{1}^{3}T_{2}+T_{3}T_{4}+T_{5}^{3}}$
&$\arraycolsep=1mm\begin{bmatrix}{4} & {3} & {3} & {12} & {5}\end{bmatrix}$&${12}$
\\
\hline
${67}$&${T_{1}^{15}T_{2}+T_{3}T_{4}+T_{5}^{5}}$
&$\arraycolsep=1mm\begin{bmatrix}{1} & {5} & {10} & {10} & {4}\end{bmatrix}$&${10}$
&${100}$&${T_{1}^{2}T_{2}+T_{3}T_{4}+T_{5}^{2}}$
&$\arraycolsep=1mm\begin{bmatrix}{3} & {4} & {4} & {6} & {5}\end{bmatrix}$&${12}$
\\
\hline
${68}$&${T_{1}^{10}T_{2}^{2}+T_{3}T_{4}+T_{5}^{5}}$
&$\arraycolsep=1mm\begin{bmatrix}{1} & {5} & {10} & {10} & {4}\end{bmatrix}$&${10}$
&${106}$&${T_{1}^{4}T_{2}+T_{3}^{3}T_{4}^{2}+T_{5}^{2}}$
&$\arraycolsep=1mm\begin{bmatrix}{1} & {4} & {2} & {1} & {4}\end{bmatrix}$&${16}$
\\
\hline
${69}$&${T_{1}^{5}T_{2}^{3}+T_{3}T_{4}+T_{5}^{5}}$
&$\arraycolsep=1mm\begin{bmatrix}{1} & {5} & {10} & {10} & {4}\end{bmatrix}$&${10}$
&${107}$&${T_{1}^{3}T_{2}+T_{3}^{3}T_{4}^{2}+T_{5}^{2}}$
&$\arraycolsep=1mm\begin{bmatrix}{2} & {8} & {4} & {1} & {7}\end{bmatrix}$&${16}$
\\
\hline
${70}$&${T_{1}^{5}T_{2}+T_{3}T_{4}+T_{5}^{5}}$
&$\arraycolsep=1mm\begin{bmatrix}{2} & {5} & {5} & {10} & {3}\end{bmatrix}$&${10}$
&${108}$&${T_{1}^{3}T_{2}+T_{3}^{3}T_{4}+T_{5}^{2}}$
&$\arraycolsep=1mm\begin{bmatrix}{1} & {1} & {1} & {1} & {2}\end{bmatrix}$&${16}$
\\
\hline
${77}$&${T_{1}^{3}T_{2}+T_{3}^{3}T_{4}^{3}+T_{5}^{2}}$
&$\arraycolsep=1mm\begin{bmatrix}{2} & {6} & {3} & {1} & {6}\end{bmatrix}$&${12}$
&${109}$&${T_{1}^{3}T_{2}+T_{3}^{2}T_{4}+T_{5}^{4}}$
&$\arraycolsep=1mm\begin{bmatrix}{1} & {1} & {1} & {2} & {1}\end{bmatrix}$&${16}$
\\
\hline
${78}$&${T_{1}^{3}T_{2}+T_{3}^{3}T_{4}^{2}+T_{5}^{2}}$
&$\arraycolsep=1mm\begin{bmatrix}{2} & {12} & {4} & {3} & {9}\end{bmatrix}$&${12}$
&${110}$&${T_{1}^{7}T_{2}+T_{3}^{2}T_{4}+T_{5}^{2}}$
&$\arraycolsep=1mm\begin{bmatrix}{1} & {1} & {2} & {4} & {4}\end{bmatrix}$&${16}$
\\
\hline
${79}$&${T_{1}^{5}T_{2}+T_{3}^{2}T_{4}+T_{5}^{3}}$
&$\arraycolsep=1mm\begin{bmatrix}{1} & {4} & {4} & {1} & {3}\end{bmatrix}$&${12}$
&${111}$&${T_{1}^{5}T_{2}^{3}+T_{3}^{2}T_{4}+T_{5}^{2}}$
&$\arraycolsep=1mm\begin{bmatrix}{1} & {1} & {2} & {4} & {4}\end{bmatrix}$&${16}$
\\
\hline
${80}$&${T_{1}^{10}T_{2}+T_{3}^{2}T_{4}+T_{5}^{2}}$
&$\arraycolsep=1mm\begin{bmatrix}{1} & {2} & {3} & {6} & {6}\end{bmatrix}$&${12}$
&${112}$&${T_{1}^{4}T_{2}+T_{3}^{2}T_{4}+T_{5}^{2}}$
&$\arraycolsep=1mm\begin{bmatrix}{1} & {2} & {2} & {2} & {3}\end{bmatrix}$&${16}$
\\
\hline
${81}$&${T_{1}^{7}T_{2}+T_{3}^{2}T_{4}+T_{5}^{2}}$
&$\arraycolsep=1mm\begin{bmatrix}{2} & {4} & {3} & {12} & {9}\end{bmatrix}$&${12}$
&${113}$&${T_{1}^{8}T_{2}+T_{3}T_{4}+T_{5}^{4}}$
&$\arraycolsep=1mm\begin{bmatrix}{1} & {4} & {4} & {8} & {3}\end{bmatrix}$&${16}$
\\
\hline
${82}$&${T_{1}^{6}T_{2}^{3}+T_{3}^{2}T_{4}+T_{5}^{2}}$
&$\arraycolsep=1mm\begin{bmatrix}{1} & {2} & {3} & {6} & {6}\end{bmatrix}$&${12}$
&${114}$&${T_{1}^{7}T_{2}+T_{3}T_{4}+T_{5}^{4}}$
&$\arraycolsep=1mm\begin{bmatrix}{1} & {1} & {4} & {4} & {2}\end{bmatrix}$&${16}$
\\
\hline
${83}$&${T_{1}^{5}T_{2}^{2}+T_{3}^{2}T_{4}+T_{5}^{2}}$
&$\arraycolsep=1mm\begin{bmatrix}{2} & {4} & {3} & {12} & {9}\end{bmatrix}$&${12}$
&${115}$&${T_{1}^{6}T_{2}+T_{3}T_{4}+T_{5}^{8}}$
&$\arraycolsep=1mm\begin{bmatrix}{1} & {2} & {4} & {4} & {1}\end{bmatrix}$&${16}$
\\
\hline
${84}$&${T_{1}^{5}T_{2}^{2}+T_{3}^{2}T_{4}+T_{5}^{2}}$
&$\arraycolsep=1mm\begin{bmatrix}{2} & {1} & {3} & {6} & {6}\end{bmatrix}$&${12}$
&${116}$&${T_{1}^{5}T_{2}^{3}+T_{3}T_{4}+T_{5}^{4}}$
&$\arraycolsep=1mm\begin{bmatrix}{1} & {1} & {4} & {4} & {2}\end{bmatrix}$&${16}$
\\
\hline
${85}$&${T_{1}^{4}T_{2}+T_{3}^{2}T_{4}+T_{5}^{2}}$
&$\arraycolsep=1mm\begin{bmatrix}{4} & {2} & {3} & {12} & {9}\end{bmatrix}$&${12}$
&${117}$&${T_{1}^{3}T_{2}^{2}+T_{3}T_{4}+T_{5}^{8}}$
&$\arraycolsep=1mm\begin{bmatrix}{2} & {1} & {4} & {4} & {1}\end{bmatrix}$&${16}$
\\
\hline
${86}$&${T_{1}^{3}T_{2}^{3}+T_{3}^{2}T_{4}+T_{5}^{2}}$
&$\arraycolsep=1mm\begin{bmatrix}{2} & {4} & {3} & {12} & {9}\end{bmatrix}$&${12}$
&${126}$&${T_{1}^{4}T_{2}+T_{3}^{3}T_{4}^{2}+T_{5}^{2}}$
&$\arraycolsep=1mm\begin{bmatrix}{1} & {6} & {2} & {2} & {5}\end{bmatrix}$&${18}$
\\
\hline
${87}$&${T_{1}^{2}T_{2}+T_{3}^{2}T_{4}+T_{5}^{3}}$
&$\arraycolsep=1mm\begin{bmatrix}{4} & {1} & {4} & {1} & {3}\end{bmatrix}$&${12}$
&${127}$&${T_{1}^{5}T_{2}+T_{3}^{3}T_{4}+T_{5}^{2}}$
&$\arraycolsep=1mm\begin{bmatrix}{1} & {1} & {1} & {3} & {3}\end{bmatrix}$&${18}$
\\
\hline
\end{tabu}\vspace{-4mm}
\end{table}

\begin{table}[t!]\centering\setlength{\tabcolsep}{3.5pt}
\renewcommand{\arraystretch}{1.25}\small\tabulinesep=1mm
\begin{tabu}{c|c|c|c||c|c|c|c}
\hline
ID&relations&gd-matrix&$-\mathcal{K}^3$&ID&relations&gd-matrix&$-\mathcal{K}^3$
\\
\hline
${128}$&${T_{1}^{3}T_{2}+T_{3}^{3}T_{4}^{3}+T_{5}^{2}}$
&$\arraycolsep=1mm\begin{bmatrix}{1} & {3} & {1} & {1} & {3}\end{bmatrix}$&${18}$
&${164}$&${T_{1}^{2}T_{2}+T_{3}^{2}T_{4}+T_{5}^{6}}$
&$\arraycolsep=1mm\begin{bmatrix}{1} & {4} & {2} & {2} & {1}\end{bmatrix}$&${24}$
\\
\hline
${129}$&${T_{1}^{8}T_{2}+T_{3}^{2}T_{4}+T_{5}^{2}}$
&$\arraycolsep=1mm\begin{bmatrix}{1} & {2} & {2} & {6} & {5}\end{bmatrix}$&${18}$
&${165}$&${T_{1}^{2}T_{2}+T_{3}^{2}T_{4}+T_{5}^{4}}$
&$\arraycolsep=1mm\begin{bmatrix}{3} & {2} & {1} & {6} & {2}\end{bmatrix}$&${24}$
\\
\hline
${130}$&${T_{1}^{4}T_{2}^{3}+T_{3}^{2}T_{4}+T_{5}^{2}}$
&$\arraycolsep=1mm\begin{bmatrix}{1} & {2} & {2} & {6} & {5}\end{bmatrix}$&${18}$
&${166}$&${T_{1}^{2}T_{2}+T_{3}^{2}T_{4}+T_{5}^{3}}$
&$\arraycolsep=1mm\begin{bmatrix}{1} & {1} & {1} & {1} & {1}\end{bmatrix}$&${24}$
\\
\hline
${131}$&${T_{1}^{15}T_{2}+T_{3}T_{4}+T_{5}^{3}}$
&$\arraycolsep=1mm\begin{bmatrix}{1} & {9} & {6} & {18} & {8}\end{bmatrix}$&${18}$
&${167}$&${T_{1}^{22}T_{2}+T_{3}T_{4}+T_{5}^{2}}$
&$\arraycolsep=1mm\begin{bmatrix}{1} & {8} & {6} & {24} & {15}\end{bmatrix}$&${24}$
\\
\hline
${132}$&${T_{1}^{11}T_{2}+T_{3}T_{4}+T_{5}^{3}}$
&$\arraycolsep=1mm\begin{bmatrix}{1} & {1} & {6} & {6} & {4}\end{bmatrix}$&${18}$
&${168}$&${T_{1}^{16}T_{2}+T_{3}T_{4}+T_{5}^{2}}$
&$\arraycolsep=1mm\begin{bmatrix}{1} & {2} & {6} & {12} & {9}\end{bmatrix}$&${24}$
\\
\hline
${133}$&${T_{1}^{10}T_{2}^{2}+T_{3}T_{4}+T_{5}^{3}}$
&$\arraycolsep=1mm\begin{bmatrix}{1} & {1} & {6} & {6} & {4}\end{bmatrix}$&${18}$
&${169}$&${T_{1}^{13}T_{2}+T_{3}T_{4}+T_{5}^{2}}$
&$\arraycolsep=1mm\begin{bmatrix}{1} & {3} & {4} & {12} & {8}\end{bmatrix}$&${24}$
\\
\hline
${134}$&${T_{1}^{8}T_{2}^{4}+T_{3}T_{4}+T_{5}^{3}}$
&$\arraycolsep=1mm\begin{bmatrix}{1} & {1} & {6} & {6} & {4}\end{bmatrix}$&${18}$
&${170}$&${T_{1}^{12}T_{2}^{3}+T_{3}T_{4}+T_{5}^{2}}$
&$\arraycolsep=1mm\begin{bmatrix}{1} & {2} & {6} & {12} & {9}\end{bmatrix}$&${24}$
\\
\hline
${135}$&${T_{1}^{7}T_{2}^{5}+T_{3}T_{4}+T_{5}^{3}}$
&$\arraycolsep=1mm\begin{bmatrix}{1} & {1} & {6} & {6} & {4}\end{bmatrix}$&${18}$
&${171}$&${T_{1}^{9}T_{2}+T_{3}T_{4}+T_{5}^{3}}$
&$\arraycolsep=1mm\begin{bmatrix}{1} & {6} & {3} & {12} & {5}\end{bmatrix}$&${24}$
\\
\hline
${136}$&${T_{1}^{7}T_{2}+T_{3}T_{4}+T_{5}^{3}}$
&$\arraycolsep=1mm\begin{bmatrix}{1} & {2} & {3} & {6} & {3}\end{bmatrix}$&${18}$
&${172}$&${T_{1}^{8}T_{2}^{5}+T_{3}T_{4}+T_{5}^{2}}$
&$\arraycolsep=1mm\begin{bmatrix}{1} & {2} & {6} & {12} & {9}\end{bmatrix}$&${24}$
\\
\hline
${137}$&${T_{1}^{6}T_{2}^{2}+T_{3}T_{4}+T_{5}^{3}}$
&$\arraycolsep=1mm\begin{bmatrix}{1} & {9} & {6} & {18} & {8}\end{bmatrix}$&${18}$
&${173}$&${T_{1}^{7}T_{2}^{4}+T_{3}T_{4}+T_{5}^{2}}$
&$\arraycolsep=1mm\begin{bmatrix}{2} & {1} & {6} & {12} & {9}\end{bmatrix}$&${24}$
\\
\hline
${138}$&${T_{1}^{6}T_{2}+T_{3}T_{4}+T_{5}^{3}}$
&$\arraycolsep=1mm\begin{bmatrix}{2} & {9} & {3} & {18} & {7}\end{bmatrix}$&${18}$
&${174}$&${T_{1}^{7}T_{2}^{3}+T_{3}T_{4}+T_{5}^{2}}$
&$\arraycolsep=1mm\begin{bmatrix}{1} & {3} & {4} & {12} & {8}\end{bmatrix}$&${24}$
\\
\hline
${139}$&${T_{1}^{5}T_{2}^{2}+T_{3}T_{4}+T_{5}^{3}}$
&$\arraycolsep=1mm\begin{bmatrix}{1} & {2} & {3} & {6} & {3}\end{bmatrix}$&${18}$
&${175}$&${T_{1}^{6}T_{2}^{3}+T_{3}T_{4}+T_{5}^{2}}$
&$\arraycolsep=1mm\begin{bmatrix}{1} & {8} & {6} & {24} & {15}\end{bmatrix}$&${24}$
\\
\hline
${140}$&${T_{1}^{5}T_{2}+T_{3}T_{4}+T_{5}^{6}}$
&$\arraycolsep=1mm\begin{bmatrix}{1} & {1} & {3} & {3} & {1}\end{bmatrix}$&${18}$
&${176}$&${T_{1}^{6}T_{2}+T_{3}T_{4}+T_{5}^{2}}$
&$\arraycolsep=1mm\begin{bmatrix}{3} & {8} & {2} & {24} & {13}\end{bmatrix}$&${24}$
\\
\hline
${141}$&${T_{1}^{4}T_{2}+T_{3}T_{4}+T_{5}^{3}}$
&$\arraycolsep=1mm\begin{bmatrix}{2} & {1} & {3} & {6} & {3}\end{bmatrix}$&${18}$
&${177}$&${T_{1}^{5}T_{2}+T_{3}T_{4}+T_{5}^{4}}$
&$\arraycolsep=1mm\begin{bmatrix}{1} & {3} & {2} & {6} & {2}\end{bmatrix}$&${24}$
\\
\hline
${142}$&${T_{1}^{3}T_{2}+T_{3}T_{4}+T_{5}^{9}}$
&$\arraycolsep=1mm\begin{bmatrix}{2} & {3} & {3} & {6} & {1}\end{bmatrix}$&${18}$
&${178}$&${T_{1}^{5}T_{2}+T_{3}T_{4}+T_{5}^{3}}$
&$\arraycolsep=1mm\begin{bmatrix}{1} & {1} & {2} & {4} & {2}\end{bmatrix}$&${24}$
\\
\hline
${143}$&${T_{1}^{2}T_{2}+T_{3}T_{4}+T_{5}^{3}}$
&$\arraycolsep=1mm\begin{bmatrix}{2} & {2} & {3} & {3} & {2}\end{bmatrix}$&${18}$
&${179}$&${T_{1}^{5}T_{2}+T_{3}T_{4}+T_{5}^{2}}$
&$\arraycolsep=1mm\begin{bmatrix}{3} & {1} & {4} & {12} & {8}\end{bmatrix}$&${24}$
\\
\hline
${145}$&${T_{1}^{3}T_{2}+T_{3}^{3}T_{4}+T_{5}^{2}}$
&$\arraycolsep=1mm\begin{bmatrix}{2} & {10} & {5} & {1} & {8}\end{bmatrix}$&${20}$
&${180}$&${T_{1}^{4}T_{2}^{2}+T_{3}T_{4}+T_{5}^{3}}$
&$\arraycolsep=1mm\begin{bmatrix}{1} & {1} & {2} & {4} & {2}\end{bmatrix}$&${24}$
\\
\hline
${146}$&${T_{1}^{26}T_{2}+T_{3}T_{4}+T_{5}^{2}}$
&$\arraycolsep=1mm\begin{bmatrix}{1} & {4} & {10} & {20} & {15}\end{bmatrix}$&${20}$
&${181}$&${T_{1}^{4}T_{2}+T_{3}T_{4}+T_{5}^{6}}$
&$\arraycolsep=1mm\begin{bmatrix}{1} & {2} & {2} & {4} & {1}\end{bmatrix}$&${24}$
\\
\hline
${147}$&${T_{1}^{18}T_{2}^{3}+T_{3}T_{4}+T_{5}^{2}}$
&$\arraycolsep=1mm\begin{bmatrix}{1} & {4} & {10} & {20} & {15}\end{bmatrix}$&${20}$
&${182}$&${T_{1}^{4}T_{2}+T_{3}T_{4}+T_{5}^{2}}$
&$\arraycolsep=1mm\begin{bmatrix}{1} & {6} & {4} & {6} & {5}\end{bmatrix}$&${24}$
\\
\hline
${148}$&${T_{1}^{10}T_{2}^{5}+T_{3}T_{4}+T_{5}^{2}}$
&$\arraycolsep=1mm\begin{bmatrix}{1} & {4} & {10} & {20} & {15}\end{bmatrix}$&${20}$
&${183}$&${T_{1}^{4}T_{2}+T_{3}T_{4}+T_{5}^{2}}$
&$\arraycolsep=1mm\begin{bmatrix}{1} & {2} & {3} & {3} & {3}\end{bmatrix}$&${24}$
\\
\hline
${149}$&${T_{1}^{7}T_{2}^{2}+T_{3}T_{4}+T_{5}^{2}}$
&$\arraycolsep=1mm\begin{bmatrix}{4} & {1} & {10} & {20} & {15}\end{bmatrix}$&${20}$
&${184}$&${T_{1}^{4}T_{2}+T_{3}T_{4}+T_{5}^{2}}$
&$\arraycolsep=1mm\begin{bmatrix}{3} & {2} & {2} & {12} & {7}\end{bmatrix}$&${24}$
\\
\hline
${150}$&${T_{1}^{3}T_{2}^{2}+T_{3}T_{4}+T_{5}^{2}}$
&$\arraycolsep=1mm\begin{bmatrix}{4} & {5} & {2} & {20} & {11}\end{bmatrix}$&${20}$
&${185}$&${T_{1}^{3}T_{2}^{2}+T_{3}T_{4}+T_{5}^{3}}$
&$\arraycolsep=1mm\begin{bmatrix}{1} & {6} & {3} & {12} & {5}\end{bmatrix}$&${24}$
\\
\hline
${151}$&${T_{1}^{2}T_{2}+T_{3}T_{4}+T_{5}^{4}}$
&$\arraycolsep=1mm\begin{bmatrix}{5} & {2} & {2} & {10} & {3}\end{bmatrix}$&${20}$
&${186}$&${T_{1}^{3}T_{2}+T_{3}T_{4}+T_{5}^{3}}$
&$\arraycolsep=1mm\begin{bmatrix}{1} & {3} & {3} & {3} & {2}\end{bmatrix}$&${24}$
\\
\hline
${160}$&${T_{1}^{4}T_{2}+T_{3}^{3}T_{4}+T_{5}^{2}}$
&$\arraycolsep=1mm\begin{bmatrix}{1} & {6} & {3} & {1} & {5}\end{bmatrix}$&${24}$
&${187}$&${T_{1}^{2}T_{2}+T_{3}T_{4}+T_{5}^{8}}$
&$\arraycolsep=1mm\begin{bmatrix}{3} & {2} & {2} & {6} & {1}\end{bmatrix}$&${24}$
\\
\hline
${161}$&${T_{1}^{4}T_{2}+T_{3}^{2}T_{4}+T_{5}^{3}}$
&$\arraycolsep=1mm\begin{bmatrix}{1} & {2} & {1} & {4} & {2}\end{bmatrix}$&${24}$
&${190}$&${T_{1}^{3}T_{2}+T_{3}^{3}T_{4}+T_{5}^{2}}$
&$\arraycolsep=1mm\begin{bmatrix}{1} & {15} & {5} & {3} & {9}\end{bmatrix}$&${30}$
\\
\hline
${162}$&${T_{1}^{3}T_{2}+T_{3}^{2}T_{4}+T_{5}^{3}}$
&$\arraycolsep=1mm\begin{bmatrix}{1} & {12} & {6} & {3} & {5}\end{bmatrix}$&${24}$
&${191}$&${T_{1}^{2}T_{2}+T_{3}^{2}T_{4}+T_{5}^{3}}$
&$\arraycolsep=1mm\begin{bmatrix}{1} & {10} & {5} & {2} & {4}\end{bmatrix}$&${30}$
\\
\hline
${163}$&${T_{1}^{2}T_{2}^{2}+T_{3}^{2}T_{4}+T_{5}^{3}}$
&$\arraycolsep=1mm\begin{bmatrix}{1} & {2} & {1} & {4} & {2}\end{bmatrix}$&${24}$
&${192}$&${T_{1}^{13}T_{2}+T_{3}T_{4}+T_{5}^{2}}$
&$\arraycolsep=1mm\begin{bmatrix}{1} & {5} & {3} & {15} & {9}\end{bmatrix}$&${30}$
\\
\hline
\end{tabu}\vspace{-4mm}
\end{table}

\begin{table}[t!]\centering\setlength{\tabcolsep}{2.5pt}
\renewcommand{\arraystretch}{1.25}\small\tabulinesep=1mm
\begin{tabu}{c|c|c|c||c|c|c|c}
\hline
ID&relations&gd-matrix&$-\mathcal{K}^3$&ID&relations&gd-matrix&$-\mathcal{K}^3$
\\ \hline
${193}$&${T_{1}^{7}T_{2}+T_{3}T_{4}+T_{5}^{3}}$
&$\arraycolsep=1mm\begin{bmatrix}{1} & {5} & {2} & {10} & {4}\end{bmatrix}$&${30}$
&${221}$&${T_{1}^{3}T_{2}+T_{3}^{2}T_{4}+T_{5}^{2}}$
&$\arraycolsep=1mm\begin{bmatrix}{4} & {2} & {1} & {12} & {7}\end{bmatrix}$&${36}$
\\ \hline
${194}$&${T_{1}^{3}T_{2}^{3}+T_{3}T_{4}+T_{5}^{2}}$
&$\arraycolsep=1mm\begin{bmatrix}{1} & {5} & {3} & {15} & {9}\end{bmatrix}$&${30}$
&${222}$&${T_{1}^{10}T_{2}+T_{3}T_{4}+T_{5}^{2}}$
&$\arraycolsep=1mm\begin{bmatrix}{1} & {4} & {2} & {12} & {7}\end{bmatrix}$&${36}$
\\ \hline
${195}$&${T_{1}^{2}T_{2}^{2}+T_{3}T_{4}+T_{5}^{3}}$
&$\arraycolsep=1mm\begin{bmatrix}{1} & {5} & {2} & {10} & {4}\end{bmatrix}$&${30}$
&${223}$&${T_{1}^{7}T_{2}+T_{3}T_{4}+T_{5}^{2}}$
&$\arraycolsep=1mm\begin{bmatrix}{1} & {1} & {2} & {6} & {4}\end{bmatrix}$&${36}$
\\ \hline
${196}$&${T_{1}T_{2}+T_{3}T_{4}+T_{5}^{5}}$
&$\arraycolsep=1mm\begin{bmatrix}{2} & {3} & {2} & {3} & {1}\end{bmatrix}$&${30}$
&${224}$&${T_{1}^{5}T_{2}^{3}+T_{3}T_{4}+T_{5}^{2}}$
&$\arraycolsep=1mm\begin{bmatrix}{1} & {1} & {2} & {6} & {4}\end{bmatrix}$&${36}$
\\ \hline
${203}$&${T_{1}^{5}T_{2}+T_{3}^{2}T_{4}+T_{5}^{2}}$
&$\arraycolsep=1mm\begin{bmatrix}{1} & {1} & {1} & {4} & {3}\end{bmatrix}$&${32}$
&${225}$&${T_{1}^{3}T_{2}^{2}+T_{3}T_{4}+T_{5}^{2}}$
&$\arraycolsep=1mm\begin{bmatrix}{4} & {1} & {2} & {12} & {7}\end{bmatrix}$&${36}$
\\ \hline
${204}$&${T_{1}^{4}T_{2}+T_{3}^{2}T_{4}+T_{5}^{2}}$
&$\arraycolsep=1mm\begin{bmatrix}{2} & {2} & {1} & {8} & {5}\end{bmatrix}$&${32}$
&${226}$&${T_{1}^{3}T_{2}+T_{3}T_{4}+T_{5}^{4}}$
&$\arraycolsep=1mm\begin{bmatrix}{1} & {1} & {1} & {3} & {1}\end{bmatrix}$&${36}$
\\ \hline
${205}$&${T_{1}^{3}T_{2}^{3}+T_{3}^{2}T_{4}+T_{5}^{2}}$
&$\arraycolsep=1mm\begin{bmatrix}{1} & {1} & {1} & {4} & {3}\end{bmatrix}$&${32}$
&${227}$&${T_{1}T_{2}+T_{3}T_{4}+T_{5}^{2}}$
&$\arraycolsep=1mm\begin{bmatrix}{3} & {1} & {2} & {2} & {2}\end{bmatrix}$&${36}$
\\ \hline
${206}$&${T_{1}^{3}T_{2}^{2}+T_{3}^{2}T_{4}+T_{5}^{2}}$
&$\arraycolsep=1mm\begin{bmatrix}{2} & {2} & {1} & {8} & {5}\end{bmatrix}$&${32}$
&${230}$&${T_{1}^{3}T_{2}+T_{3}T_{4}+T_{5}^{5}}$
&$\arraycolsep=1mm\begin{bmatrix}{1} & {2} & {1} & {4} & {1}\end{bmatrix}$&${40}$
\\ \hline
${207}$&${T_{1}^{2}T_{2}+T_{3}^{2}T_{4}+T_{5}^{2}}$
&$\arraycolsep=1mm\begin{bmatrix}{1} & {2} & {1} & {2} & {2}\end{bmatrix}$&${32}$
&${231}$&${T_{1}^{2}T_{2}+T_{3}T_{4}+T_{5}^{5}}$
&$\arraycolsep=1mm\begin{bmatrix}{2} & {1} & {1} & {4} & {1}\end{bmatrix}$&${40}$
\\ \hline
${208}$&${T_{1}^{11}T_{2}+T_{3}T_{4}+T_{5}^{2}}$
&$\arraycolsep=1mm\begin{bmatrix}{1} & {1} & {4} & {8} & {6}\end{bmatrix}$&${32}$
&${232}$&${T_{1}^{5}T_{2}+T_{3}T_{4}+T_{5}^{2}}$
&$\arraycolsep=1mm\begin{bmatrix}{3} & {7} & {1} & {21} & {11}\end{bmatrix}$&${42}$
\\ \hline
${209}$&${T_{1}^{9}T_{2}^{3}+T_{3}T_{4}+T_{5}^{2}}$
&$\arraycolsep=1mm\begin{bmatrix}{1} & {1} & {4} & {8} & {6}\end{bmatrix}$&${32}$
&${233}$&${T_{1}^{4}T_{2}+T_{3}T_{4}+T_{5}^{3}}$
&$\arraycolsep=1mm\begin{bmatrix}{2} & {7} & {1} & {14} & {5}\end{bmatrix}$&${42}$
\\ \hline
${210}$&${T_{1}^{8}T_{2}+T_{3}T_{4}+T_{5}^{2}}$
&$\arraycolsep=1mm\begin{bmatrix}{1} & {2} & {2} & {8} & {5}\end{bmatrix}$&${32}$
&${234}$&${T_{1}^{2}T_{2}+T_{3}T_{4}+T_{5}^{7}}$
&$\arraycolsep=1mm\begin{bmatrix}{2} & {3} & {1} & {6} & {1}\end{bmatrix}$&${42}$
\\ \hline
${211}$&${T_{1}^{7}T_{2}^{5}+T_{3}T_{4}+T_{5}^{2}}$
&$\arraycolsep=1mm\begin{bmatrix}{1} & {1} & {4} & {8} & {6}\end{bmatrix}$&${32}$
&${235}$&${T_{1}^{5}T_{2}+T_{3}T_{4}+T_{5}^{3}}$
&$\arraycolsep=1mm\begin{bmatrix}{1} & {4} & {1} & {8} & {3}\end{bmatrix}$&${48}$
\\ \hline
${212}$&${T_{1}^{4}T_{2}^{3}+T_{3}T_{4}+T_{5}^{2}}$
&$\arraycolsep=1mm\begin{bmatrix}{1} & {2} & {2} & {8} & {5}\end{bmatrix}$&${32}$
&${236}$&${T_{1}^{2}T_{2}+T_{3}T_{4}+T_{5}^{3}}$
&$\arraycolsep=1mm\begin{bmatrix}{4} & {1} & {1} & {8} & {3}\end{bmatrix}$&${48}$
\\
\hline
${213}$&${T_{1}^{3}T_{2}+T_{3}T_{4}+T_{5}^{2}}$
&$\arraycolsep=1mm\begin{bmatrix}{1} & {1} & {2} & {2} & {2}\end{bmatrix}$&${32}$
&${237}$&${T_{1}T_{2}+T_{3}T_{4}+T_{5}^{3}}$
&$\arraycolsep=1mm\begin{bmatrix}{1} & {2} & {1} & {2} & {1}\end{bmatrix}$&${48}$
\\
\hline
${214}$&${T_{1}^{2}T_{2}+T_{3}T_{4}+T_{5}^{4}}$
&$\arraycolsep=1mm\begin{bmatrix}{1} & {2} & {2} & {2} & {1}\end{bmatrix}$&${32}$
&${238}$&${T_{1}T_{2}+T_{3}T_{4}+T_{5}^{2}}$
&$\arraycolsep=1mm\begin{bmatrix}{1} & {3} & {1} & {3} & {2}\end{bmatrix}$&${48}$
\\
\hline
${215}$&${T_{1}^{2}T_{2}+T_{3}T_{4}+T_{5}^{2}}$
&$\arraycolsep=1mm\begin{bmatrix}{1} & {4} & {2} & {4} & {3}\end{bmatrix}$&${32}$
&${241}$&${T_{1}^{5}T_{2}+T_{3}T_{4}+T_{5}^{2}}$
&$\arraycolsep=1mm\begin{bmatrix}{1} & {1} & {1} & {5} & {3}\end{bmatrix}$&${50}$
\\
\hline
${217}$&${T_{1}^{6}T_{2}+T_{3}^{2}T_{4}+T_{5}^{2}}$
&$\arraycolsep=1mm\begin{bmatrix}{1} & {2} & {1} & {6} & {4}\end{bmatrix}$&${36}$
&${242}$&${T_{1}^{3}T_{2}^{3}+T_{3}T_{4}+T_{5}^{2}}$
&$\arraycolsep=1mm\begin{bmatrix}{1} & {1} & {1} & {5} & {3}\end{bmatrix}$&${50}$
\\
\hline
${218}$&${T_{1}^{5}T_{2}+T_{3}^{2}T_{4}+T_{5}^{2}}$
&$\arraycolsep=1mm\begin{bmatrix}{2} & {4} & {1} & {12} & {7}\end{bmatrix}$&${36}$
&${243}$&${T_{1}^{7}T_{2}+T_{3}T_{4}+T_{5}^{2}}$
&$\arraycolsep=1mm\begin{bmatrix}{1} & {3} & {1} & {9} & {5}\end{bmatrix}$&${54}$
\\
\hline
${219}$&${T_{1}^{3}T_{2}^{2}+T_{3}^{2}T_{4}+T_{5}^{2}}$
&$\arraycolsep=1mm\begin{bmatrix}{2} & {4} & {1} & {12} & {7}\end{bmatrix}$&${36}$
&${244}$&${T_{1}^{3}T_{2}+T_{3}T_{4}+T_{5}^{2}}$
&$\arraycolsep=1mm\begin{bmatrix}{3} & {1} & {1} & {9} & {5}\end{bmatrix}$&${54}$
\\
\hline
${220}$&${T_{1}^{3}T_{2}^{2}+T_{3}^{2}T_{4}+T_{5}^{2}}$
&$\arraycolsep=1mm\begin{bmatrix}{2} & {1} & {1} & {6} & {4}\end{bmatrix}$&${36}$
&${245}$&${T_{1}T_{2}+T_{3}T_{4}+T_{5}^{2}}$
&$\arraycolsep=1mm\begin{bmatrix}{1} & {1} & {1} & {1} & {1}\end{bmatrix}$&${54}$
\\ \hline
\end{tabu}\vspace{-2mm}
\end{table}
\end{class-list}

\begin{class-list}
Non-toric, $\QQ$-factorial, Gorenstein, log terminal Fano
threefolds of Picard number one with an effective two-torus action:
Specifying data for divisor class group $\ZZ$ and format $(2,2,1,1,0)$.
\begin{table}[h!]\centering\setlength{\tabcolsep}{3.5pt}
\renewcommand{\arraystretch}{1.25}\small\tabulinesep=1mm
\begin{tabu}{c|c|c|c}
\hline
ID&relations&gd-matrix&$-\mathcal{K}^3$
\\ \hline
${30}$&${T_{1}^{3}T_{2}+T_{3}T_{4}+T_{5}^{3}},\ {\lambda T_{3}T_{4}+T_{5}^{3}+T_{6}^{2}}$
&$\arraycolsep=1mm\begin{bmatrix}{1} & {3} & {3} & {3} & {2} & {3}\end{bmatrix}$&${6}$
\\
\hline
${76}$&${T_{1}^{2}T_{2}+T_{3}T_{4}+T_{5}^{3}},\ {\lambda T_{3}T_{4}+T_{5}^{3}+T_{6}^{2}}$
&$\arraycolsep=1mm\begin{bmatrix}{1} & {4} & {2} & {4} & {2} & {3}\end{bmatrix}$&${12}$
\\ \hline
${189}$&${T_{1}T_{2}+T_{3}T_{4}+T_{5}^{3}},\ {\lambda T_{3}T_{4}+T_{5}^{3}+T_{6}^{2}}$
&$\arraycolsep=1mm\begin{bmatrix}{1} & {5} & {1} & {5} & {2} & {3}\end{bmatrix}$&${30}$
\\ \hline
\end{tabu}
\end{table}
\end{class-list}

\begin{class-list}
Non-toric, $\QQ$-factorial, Gorenstein, log terminal Fano
threefolds of Picard number one with an effective two-torus action:
Specifying data for divisor class group $\ZZ$ and format $(3,1,1,0)$.
\begin{table}[h!]\centering\setlength{\tabcolsep}{3.5pt}
\renewcommand{\arraystretch}{1.25}\small\tabulinesep=1mm
\begin{tabu}{c|c|c|c||c|c|c|c}\hline
ID&relations&gd-matrix&$-\mathcal{K}^3$&ID&relations&gd-matrix&$-\mathcal{K}^3$
\\ \hline
${3}$&${T_{1}^{5}T_{2}^{3}T_{3}^{2}+T_{4}^{3}+T_{5}^{2}}$
&$\arraycolsep=1mm\begin{bmatrix} {1} & {1} & {2} & {4} & {6}\end{bmatrix}$&${2}$
&${51}$&${T_{1}^{2}T_{2}^{2}T_{3}+T_{4}^{5}+T_{5}^{2}}$
&$\arraycolsep=1mm\begin{bmatrix}{2} & {1} & {4} & {2} & {5}\end{bmatrix}$&${8}$
\\
\hline
${4}$&${T_{1}^{4}T_{2}^{3}T_{3}^{2}+T_{4}^{3}+T_{5}^{2}}$
&$\arraycolsep=1mm\begin{bmatrix}{1} & {2} & {1} & {4} & {6}\end{bmatrix}$&${2}$
&${63}$&${T_{1}T_{2}T_{3}+T_{4}^{7}+T_{5}^{3}}$
&$\arraycolsep=1mm\begin{bmatrix}{1} & {10} & {10} & {3} & {7}\end{bmatrix}$&${10}$
\\ \hline
${5}$&${T_{1}^{5}T_{2}^{5}T_{3}+T_{4}^{3}+T_{5}^{2}}$
&$\arraycolsep=1mm\begin{bmatrix}{1} & {1} & {2} & {4} & {6}\end{bmatrix}$&${2}$
&${64}$&${T_{1}T_{2}T_{3}+T_{4}^{3}+T_{5}^{2}}$
&$\arraycolsep=1mm\begin{bmatrix}{2} & {5} & {5} & {4} & {6}\end{bmatrix}$&${10}$
\\ \hline
${6}$&${T_{1}^{5}T_{2}^{3}T_{3}+T_{4}^{3}+T_{5}^{2}}$
&$\arraycolsep=1mm\begin{bmatrix}{1} & {2} & {1} & {4} & {6}\end{bmatrix}$&${2}$
&${75}$&${T_{1}^{3}T_{2}^{3}T_{3}+T_{4}^{3}+T_{5}^{2}}$
&$\arraycolsep=1mm\begin{bmatrix}{3} & {1} & {12} & {8} & {12}\end{bmatrix}$&${12}$
\\ \hline
${7}$&${T_{1}^{5}T_{2}T_{3}+T_{4}^{3}+T_{5}^{2}}$
&$\arraycolsep=1mm\begin{bmatrix}{2} & {1} & {1} & {4} & {6}\end{bmatrix}$&${2}$
&${104}$&${T_{1}^{2}T_{2}T_{3}+T_{4}^{7}+T_{5}^{2}}$
&$\arraycolsep=1mm\begin{bmatrix}{1} & {4} & {8} & {2} & {7}\end{bmatrix}$&${16}$
\\ \hline
${8}$&${T_{1}^{4}T_{2}^{3}T_{3}+T_{4}^{3}+T_{5}^{2}}$
&$\arraycolsep=1mm\begin{bmatrix}{2} & {1} & {1} & {4} & {6}\end{bmatrix}$&${2}$
&${105}$&${T_{1}^{2}T_{2}T_{3}+T_{4}^{3}+T_{5}^{2}}$
&$\arraycolsep=1mm\begin{bmatrix}{1} & {2} & {2} & {2} & {3}\end{bmatrix}$&${16}$
\\ \hline
${9}$&${T_{1}^{3}T_{2}^{2}T_{3}+T_{4}^{5}+T_{5}^{2}}$
&$\arraycolsep=1mm\begin{bmatrix}{2} & {1} & {2} & {2} & {5}\end{bmatrix}$&${2}$
&${121}$&${T_{1}^{5}T_{2}T_{3}+T_{4}^{3}+T_{5}^{2}}$
&$\arraycolsep=1mm\begin{bmatrix}{1} & {6} & {1} & {4} & {6}\end{bmatrix}$&${18}$
\\ \hline
${13}$&${T_{1}^{2}T_{2}^{2}T_{3}+T_{4}^{7}+T_{5}^{2}}$
&$\arraycolsep=1mm\begin{bmatrix}{4} & {1} & {4} & {2} & {7}\end{bmatrix}$&${4}$
&${122}$&${T_{1}^{4}T_{2}^{2}T_{3}+T_{4}^{3}+T_{5}^{2}}$
&$\arraycolsep=1mm\begin{bmatrix}{1} & {1} & {6} & {4} & {6}\end{bmatrix}$&${18}$
\\ \hline
${25}$&${T_{1}^{5}T_{2}^{2}T_{3}+T_{4}^{3}+T_{5}^{2}}$
&$\arraycolsep=1mm\begin{bmatrix}{2} & {1} & {6} & {6} & {9}\end{bmatrix}$&${6}$
&${123}$&${T_{1}^{3}T_{2}^{3}T_{3}+T_{4}^{3}+T_{5}^{2}}$
&$\arraycolsep=1mm\begin{bmatrix}{1} & {1} & {6} & {4} & {6}\end{bmatrix}$&${18}$
\\ \hline
${26}$&${T_{1}^{4}T_{2}^{4}T_{3}+T_{4}^{3}+T_{5}^{2}}$
&$\arraycolsep=1mm\begin{bmatrix}{2} & {1} & {6} & {6} & {9}\end{bmatrix}$&${6}$
&${124}$&${T_{1}^{2}T_{2}T_{3}+T_{4}^{5}+T_{5}^{2}}$
&$\arraycolsep=1mm\begin{bmatrix}{1} & {2} & {6} & {2} & {5}\end{bmatrix}$&${18}$
\\ \hline
${27}$&${T_{1}^{4}T_{2}^{2}T_{3}+T_{4}^{3}+T_{5}^{2}}$
&$\arraycolsep=1mm\begin{bmatrix}{1} & {6} & {2} & {6} & {9}\end{bmatrix}$&${6}$
&${125}$&${T_{1}T_{2}T_{3}+T_{4}^{7}+T_{5}^{2}}$
&$\arraycolsep=1mm\begin{bmatrix}{1} & {9} & {18} & {4} & {14}\end{bmatrix}$&${18}$
\\ \hline
${28}$&${T_{1}^{3}T_{2}T_{3}+T_{4}^{3}+T_{5}^{2}}$
&$\arraycolsep=1mm\begin{bmatrix}{2} & {3} & {3} & {4} & {6}\end{bmatrix}$&${6}$
&${144}$&${T_{1}^{2}T_{2}^{2}T_{3}+T_{4}^{3}+T_{5}^{2}}$
&$\arraycolsep=1mm\begin{bmatrix}{4} & {1} & {20} & {10} & {15}\end{bmatrix}$&${20}$
\\ \hline
${29}$&${T_{1}^{2}T_{2}T_{3}+T_{4}^{7}+T_{5}^{2}}$
&$\arraycolsep=1mm\begin{bmatrix}{3} & {2} & {6} & {2} & {7}\end{bmatrix}$&${6}$
&${158}$&${T_{1}^{4}T_{2}T_{3}+T_{4}^{3}+T_{5}^{2}}$
&$\arraycolsep=1mm\begin{bmatrix}{1} & {12} & {2} & {6} & {9}\end{bmatrix}$&${24}$
\\ \hline
${46}$&${T_{1}^{5}T_{2}^{3}T_{3}+T_{4}^{3}+T_{5}^{2}}$
&$\arraycolsep=1mm\begin{bmatrix}{1} & {1} & {4} & {4} & {6}\end{bmatrix}$&${8}$
&${159}$&${T_{1}^{2}T_{2}^{2}T_{3}+T_{4}^{3}+T_{5}^{2}}$
&$\arraycolsep=1mm\begin{bmatrix}{2} & {1} & {12} & {6} & {9}\end{bmatrix}$&${24}$
\\ \hline
${47}$&${T_{1}^{4}T_{2}^{4}T_{3}+T_{4}^{3}+T_{5}^{2}}$
&$\arraycolsep=1mm\begin{bmatrix}{1} & {1} & {4} & {4} & {6}\end{bmatrix}$&${8}$
&${188}$&${T_{1}T_{2}T_{3}+T_{4}^{3}+T_{5}^{2}}$
&$\arraycolsep=1mm\begin{bmatrix}{1} & {5} & {30} & {12} & {18}\end{bmatrix}$&${30}$
\\
\hline
${48}$&${T_{1}^{4}T_{2}T_{3}+T_{4}^{3}+T_{5}^{2}}$
&$\arraycolsep=1mm\begin{bmatrix}{1} & {1} & {1} & {2} & {3}\end{bmatrix}$&${8}$
&${201}$&${T_{1}^{3}T_{2}T_{3}+T_{4}^{3}+T_{5}^{2}}$
&$\arraycolsep=1mm\begin{bmatrix}{1} & {8} & {1} & {4} & {6}\end{bmatrix}$&${32}$
\\
\hline
${49}$&${T_{1}^{3}T_{2}^{2}T_{3}+T_{4}^{3}+T_{5}^{2}}$
&$\arraycolsep=1mm\begin{bmatrix}{1} & {4} & {1} & {4} & {6}\end{bmatrix}$&${8}$
&${202}$&${T_{1}^{2}T_{2}^{2}T_{3}+T_{4}^{3}+T_{5}^{2}}$
&$\arraycolsep=1mm\begin{bmatrix}{1} & {1} & {8} & {4} & {6}\end{bmatrix}$&${32}$
\\
\hline
${50}$&${T_{1}^{3}T_{2}^{2}T_{3}+T_{4}^{3}+T_{5}^{2}}$
&$\arraycolsep=1mm\begin{bmatrix}{1} & {1} & {1} & {2} & {3}\end{bmatrix}$&${8}$
&${240}$&${T_{1}T_{2}T_{3}+T_{4}^{3}+T_{5}^{2}}$
&$\arraycolsep=1mm\begin{bmatrix}{1} & {1} & {10} & {4} & {6}\end{bmatrix}$&${50}$
\\
\hline
\end{tabu}\vspace{-4mm}
\end{table}
\end{class-list}

\begin{class-list}
Non-toric, $\QQ$-factorial, Gorenstein, log terminal Fano
threefolds of Picard number one with an effective two-torus action:
Specifying data for divisor class group $\ZZ$ and format $(2,1,1,1)$.
\begin{table}[h!]\centering\setlength{\tabcolsep}{3.5pt}
\renewcommand{\arraystretch}{1.3}\small\tabulinesep=1mm
\begin{tabu}{c|c|c|c||c|c|c|c}
\hline
ID&relations&gd-matrix&$-\mathcal{K}^3$&ID&relations&gd-matrix&$-\mathcal{K}^3$
\\
\hline
${1}$&${T_{1}^{5}T_{2}^{2}+T_{3}^{3}+T_{4}^{2}}$
&$\arraycolsep=1mm\begin{bmatrix}{2} & {1} & {4} & {6} & {1}\end{bmatrix}$&${2}$
&${20}$&${T_{1}^{3}T_{2}^{2}+T_{3}^{3}+T_{4}^{2}}$
&$\arraycolsep=1mm\begin{bmatrix}{2} & {6} & {6} & {9} & {1}\end{bmatrix}$&${6}$
\\
\hline
${2}$&${T_{1}^{3}T_{2}^{2}+T_{3}^{5}+T_{4}^{2}}$
&$\arraycolsep=1mm\begin{bmatrix}{2} & {2} & {2} & {5} & {1}\end{bmatrix}$&${2}$
&${21}$&${T_{1}^{3}T_{2}+T_{3}^{3}+T_{4}^{2}}$
&$\arraycolsep=1mm\begin{bmatrix}{3} & {3} & {4} & {6} & {2}\end{bmatrix}$&${6}$
\\
\hline
${12}$&${T_{1}^{2}T_{2}+T_{3}^{4}+T_{4}^{3}}$
&$\arraycolsep=1mm\begin{bmatrix}{4} & {4} & {3} & {4} & {1}\end{bmatrix}$&${4}$
&${22}$&${T_{1}^{2}T_{2}+T_{3}^{4}+T_{4}^{3}}$
&$\arraycolsep=1mm\begin{bmatrix}{3} & {6} & {3} & {4} & {2}\end{bmatrix}$&${6}$
\\
\hline
${19}$&${T_{1}^{3}T_{2}^{2}+T_{3}^{3}+T_{4}^{2}}$
&$\arraycolsep=1mm\begin{bmatrix}{2} & {3} & {4} & {6} & {3}\end{bmatrix}$&${6}$
&${23}$&${T_{1}^{2}T_{2}+T_{3}^{9}+T_{4}^{2}}$
&$\arraycolsep=1mm\begin{bmatrix}{6} & {6} & {2} & {9} & {1}\end{bmatrix}$&${6}$
\\
\hline
\end{tabu}\vspace{-2mm}
\end{table}

\begin{table}[t!]\centering\setlength{\tabcolsep}{3.5pt}
\renewcommand{\arraystretch}{1.3}\small\tabulinesep=1mm
\begin{tabu}{c|c|c|c||c|c|c|c}
\hline
ID&relations&gd-matrix&$-\mathcal{K}^3$&ID&relations&gd-matrix&$-\mathcal{K}^3$
\\
\hline
${24}$&${T_{1}^{2}T_{2}+T_{3}^{7}+T_{4}^{2}}$
&$\arraycolsep=1mm\begin{bmatrix}{6} & {2} & {2} & {7} & {3}\end{bmatrix}$&${6}$
&${119}$&${T_{1}^{2}T_{2}+T_{3}^{3}+T_{4}^{2}}$
&$\arraycolsep=1mm\begin{bmatrix}{2} & {2} & {2} & {3} & {3}\end{bmatrix}$&${18}$
\\
\hline
${43}$&${T_{1}^{5}T_{2}+T_{3}^{3}+T_{4}^{2}}$
&$\arraycolsep=1mm\begin{bmatrix}{1} & {1} & {2} & {3} & {1}\end{bmatrix}$&${8}$
&${120}$&${T_{1}T_{2}+T_{3}^{5}+T_{4}^{4}}$
&$\arraycolsep=1mm\begin{bmatrix}{2} & {18} & {4} & {5} & {9}\end{bmatrix}$&${18}$
\\
\hline
${44}$&${T_{1}^{3}T_{2}+T_{3}^{5}+T_{4}^{2}}$
&$\arraycolsep=1mm\begin{bmatrix}{2} & {4} & {2} & {5} & {1}\end{bmatrix}$&${8}$
&${152}$&${T_{1}T_{2}+T_{3}^{8}+T_{4}^{3}}$
&$\arraycolsep=1mm\begin{bmatrix}{2} & {22} & {3} & {8} & {11}\end{bmatrix}$&${22}$
\\
\hline
${45}$&${T_{1}^{2}T_{2}+T_{3}^{5}+T_{4}^{2}}$
&$\arraycolsep=1mm\begin{bmatrix}{4} & {2} & {2} & {5} & {1}\end{bmatrix}$&${8}$
&${153}$&${T_{1}^{3}T_{2}+T_{3}^{3}+T_{4}^{2}}$
&$\arraycolsep=1mm\begin{bmatrix}{1} & {3} & {2} & {3} & {3}\end{bmatrix}$&${24}$
\\
\hline
${59}$&${T_{1}^{2}T_{2}+T_{3}^{5}+T_{4}^{2}}$
&$\arraycolsep=1mm\begin{bmatrix}{5} & {10} & {4} & {10} & {1}\end{bmatrix}$&${10}$
&${154}$&${T_{1}^{3}T_{2}+T_{3}^{3}+T_{4}^{2}}$
&$\arraycolsep=1mm\begin{bmatrix}{2} & {12} & {6} & {9} & {1}\end{bmatrix}$&${24}$
\\
\hline
${60}$&${T_{1}^{2}T_{2}+T_{3}^{3}+T_{4}^{2}}$
&$\arraycolsep=1mm\begin{bmatrix}{5} & {2} & {4} & {6} & {5}\end{bmatrix}$&${10}$
&${155}$&${T_{1}T_{2}+T_{3}^{9}+T_{4}^{2}}$
&$\arraycolsep=1mm\begin{bmatrix}{6} & {12} & {2} & {9} & {1}\end{bmatrix}$&${24}$
\\
\hline
${61}$&${T_{1}T_{2}+T_{3}^{5}+T_{4}^{4}}$
&$\arraycolsep=1mm\begin{bmatrix}{10} & {10} & {4} & {5} & {1}\end{bmatrix}$&${10}$
&${156}$&${T_{1}T_{2}+T_{3}^{7}+T_{4}^{2}}$
&$\arraycolsep=1mm\begin{bmatrix}{2} & {12} & {2} & {7} & {3}\end{bmatrix}$&${24}$
\\
\hline
${62}$&${T_{1}T_{2}+T_{3}^{5}+T_{4}^{3}}$
&$\arraycolsep=1mm\begin{bmatrix}{5} & {10} & {3} & {5} & {2}\end{bmatrix}$&${10}$
&${157}$&${T_{1}T_{2}+T_{3}^{3}+T_{4}^{2}}$
&$\arraycolsep=1mm\begin{bmatrix}{3} & {3} & {2} & {3} & {1}\end{bmatrix}$&${24}$
\\
\hline
${71}$&${T_{1}^{4}T_{2}+T_{3}^{3}+T_{4}^{2}}$
&$\arraycolsep=1mm\begin{bmatrix}{3} & {12} & {8} & {12} & {1}\end{bmatrix}$&${12}$
&${197}$&${T_{1}^{4}T_{2}+T_{3}^{3}+T_{4}^{2}}$
&$\arraycolsep=1mm\begin{bmatrix}{1} & {8} & {4} & {6} & {1}\end{bmatrix}$&${32}$
\\
\hline
${72}$&${T_{1}^{2}T_{2}+T_{3}^{9}+T_{4}^{2}}$
&$\arraycolsep=1mm\begin{bmatrix}{3} & {12} & {2} & {9} & {4}\end{bmatrix}$&${12}$
&${198}$&${T_{1}^{2}T_{2}+T_{3}^{5}+T_{4}^{2}}$
&$\arraycolsep=1mm\begin{bmatrix}{1} & {8} & {2} & {5} & {2}\end{bmatrix}$&${32}$
\\
\hline
${73}$&${T_{1}T_{2}+T_{3}^{8}+T_{4}^{3}}$
&$\arraycolsep=1mm\begin{bmatrix}{12} & {12} & {3} & {8} & {1}\end{bmatrix}$&${12}$
&${199}$&${T_{1}^{2}T_{2}+T_{3}^{3}+T_{4}^{2}}$
&$\arraycolsep=1mm\begin{bmatrix}{1} & {4} & {2} & {3} & {4}\end{bmatrix}$&${32}$
\\
\hline
${74}$&${T_{1}T_{2}+T_{3}^{5}+T_{4}^{3}}$
&$\arraycolsep=1mm\begin{bmatrix}{3} & {12} & {3} & {5} & {4}\end{bmatrix}$&${12}$
&${200}$&${T_{1}T_{2}+T_{3}^{5}+T_{4}^{2}}$
&$\arraycolsep=1mm\begin{bmatrix}{2} & {8} & {2} & {5} & {1}\end{bmatrix}$&${32}$
\\
\hline
${101}$&${T_{1}^{4}T_{2}+T_{3}^{3}+T_{4}^{2}}$
&$\arraycolsep=1mm\begin{bmatrix}{1} & {2} & {2} & {3} & {2}\end{bmatrix}$&${16}$
&${216}$&${T_{1}^{2}T_{2}+T_{3}^{7}+T_{4}^{2}}$
&$\arraycolsep=1mm\begin{bmatrix}{1} & {12} & {2} & {7} & {4}\end{bmatrix}$&${36}$
\\
\hline
${102}$&${T_{1}^{2}T_{2}+T_{3}^{3}+T_{4}^{2}}$
&$\arraycolsep=1mm\begin{bmatrix}{2} & {2} & {2} & {3} & {1}\end{bmatrix}$&${16}$
&${228}$&${T_{1}T_{2}+T_{3}^{7}+T_{4}^{3}}$
&$\arraycolsep=1mm\begin{bmatrix}{1} & {20} & {3} & {7} & {10}\end{bmatrix}$&${40}$
\\
\hline
${103}$&${T_{1}T_{2}+T_{3}^{4}+T_{4}^{3}}$
&$\arraycolsep=1mm\begin{bmatrix}{4} & {8} & {3} & {4} & {1}\end{bmatrix}$&${16}$
&${229}$&${T_{1}T_{2}+T_{3}^{3}+T_{4}^{2}}$
&$\arraycolsep=1mm\begin{bmatrix}{1} & {5} & {2} & {3} & {5}\end{bmatrix}$&${40}$
\\
\hline
${118}$&${T_{1}^{2}T_{2}+T_{3}^{5}+T_{4}^{2}}$
&$\arraycolsep=1mm\begin{bmatrix}{2} & {6} & {2} & {5} & {1}\end{bmatrix}$&${18}$
&${239}$&${T_{1}^{2}T_{2}+T_{3}^{3}+T_{4}^{2}}$
&$\arraycolsep=1mm\begin{bmatrix}{1} & {10} & {4} & {6} & {1}\end{bmatrix}$&${50}$
\\
\hline
\end{tabu}
\end{table}
\end{class-list}

\begin{class-list}\label{class-list(0,2)}
Non-toric, $\QQ$-factorial, Gorenstein, log terminal Fano
threefolds of Picard number one with an effective two-torus action:
Specifying data for divisor class group $\ZZ \oplus \ZZ/2\ZZ$ and
format $(2,2,1,0)$.
\begin{table}[h!]\centering\setlength{\tabcolsep}{3.5pt}
\renewcommand{\arraystretch}{1.3}\small\tabulinesep=1mm
\begin{tabu}{c|c|c|c||c|c|c|c}
\hline
ID&relations&gd-matrix&$-\mathcal{K}^3$&ID&relations&gd-matrix&$-\mathcal{K}^3$
\\
\hline
${252}$&${T_{1}^{4}T_{2}^{2}+T_{3}^{3}T_{4}^{2}+T_{5}^{2}}$
&$\arraycolsep=1mm\begin{bmatrix}{1} & {2} & {2} & {1} & {4} \\ {\bar{0}} & {\bar{0}} & {\bar{0}} & {\bar{1}} & {\bar{1}}\end{bmatrix}$&${2}$
&${273}$&${T_{1}^{4}T_{2}+T_{3}^{2}T_{4}^{2}+T_{5}^{2}}$
&$\arraycolsep=1mm\begin{bmatrix}{2} & {2} & {1} & {4} & {5} \\ {\bar{1}} & {\bar{0}} & {\bar{0}} & {\bar{0}} & {\bar{1}}\end{bmatrix}$&${4}$
\\
\hline
${253}$&${T_{1}^{3}T_{2}^{2}+T_{3}^{3}T_{4}^{2}+T_{5}^{2}}$
&$\arraycolsep=1mm\begin{bmatrix}{2} & {1} & {2} & {1} & {4} \\ {\bar{0}} & {\bar{1}} & {\bar{0}} & {\bar{0}} & {\bar{1}}\end{bmatrix}$&${2}$
&${274}$&${T_{1}^{3}T_{2}^{2}+T_{3}^{2}T_{4}^{2}+T_{5}^{2}}$
&$\arraycolsep=1mm\begin{bmatrix}{2} & {2} & {1} & {4} & {5} \\ {\bar{0}} & {\bar{1}} & {\bar{1}} & {\bar{0}} & {\bar{0}}\end{bmatrix}$&${4}$
\\
\hline
${254}$&${T_{1}^{4}T_{2}^{4}+T_{3}^{3}T_{4}+T_{5}^{2}}$
&$\arraycolsep=1mm\begin{bmatrix}{1} & {1} & {2} & {2} & {4} \\ {\bar{1}} & {\bar{0}} & {\bar{0}} & {\bar{0}} & {\bar{1}}\end{bmatrix}$&${2}$
&${275}$&${T_{1}^{8}T_{2}^{4}+T_{3}^{2}T_{4}+T_{5}^{2}}$
&$\arraycolsep=1mm\begin{bmatrix}{1} & {1} & {4} & {4} & {6} \\ {\bar{1}} & {\bar{0}} & {\bar{0}} & {\bar{0}} & {\bar{1}}\end{bmatrix}$&${4}$
\\
\hline
${271}$&${T_{1}^{4}T_{2}+T_{3}^{3}T_{4}^{3}+T_{5}^{2}}$
&$\arraycolsep=1mm\begin{bmatrix}{1} & {2} & {1} & {1} & {3} \\ {\bar{0}} & {\bar{0}} & {\bar{1}} & {\bar{1}} & {\bar{0}}\end{bmatrix}$&${4}$
&${276}$&${T_{1}^{4}T_{2}^{2}+T_{3}^{2}T_{4}+T_{5}^{2}}$
&$\arraycolsep=1mm\begin{bmatrix}{2} & {1} & {4} & {2} & {5} \\ {\bar{1}} & {\bar{0}} & {\bar{0}} & {\bar{0}} & {\bar{1}}\end{bmatrix}$&${4}$
\\
\hline
${272}$&${T_{1}^{4}T_{2}^{2}+T_{3}^{3}T_{4}+T_{5}^{2}}$
&$\arraycolsep=1mm\begin{bmatrix}{2} & {1} & {2} & {4} & {5} \\ {\bar{1}} & {\bar{0}} & {\bar{0}} & {\bar{0}} & {\bar{1}}\end{bmatrix}$&${4}$
&${291}$&${T_{1}^{10}T_{2}+T_{3}^{2}T_{4}+T_{5}^{2}}$
&$\arraycolsep=1mm\begin{bmatrix}{1} & {2} & {3} & {6} & {6} \\ {\bar{0}} & {\bar{0}} & {\bar{1}} & {\bar{0}} & {\bar{1}}\end{bmatrix}$&${6}$
\\
\hline
\end{tabu}
\end{table}

\begin{table}[h!]\centering\setlength{\tabcolsep}{3.5pt}
\renewcommand{\arraystretch}{1.3}\small\tabulinesep=1mm
\begin{tabu}{c|c|c|c||c|c|c|c}
\hline
ID&relations&gd-matrix&$-\mathcal{K}^3$&ID&relations&gd-matrix&$-\mathcal{K}^3$
\\
\hline
${292}$&${T_{1}^{8}T_{2}^{2}+T_{3}^{2}T_{4}+T_{5}^{2}}$
&$\arraycolsep=1mm\begin{bmatrix}{1} & {2} & {3} & {6} & {6} \\ {\bar{0}} & {\bar{0}} & {\bar{1}} & {\bar{0}} & {\bar{1}}\end{bmatrix}$&${6}$
&${332}$&${T_{1}^{2}T_{2}+T_{3}^{2}T_{4}^{2}+T_{5}^{2}}$
&$\arraycolsep=1mm\begin{bmatrix}{2} & {2} & {2} & {1} & {3} \\ {\bar{1}} & {\bar{0}} & {\bar{1}} & {\bar{0}} & {\bar{0}}\end{bmatrix}$&${8}$
\\
\hline
${293}$&${T_{1}^{6}T_{2}^{3}+T_{3}^{2}T_{4}+T_{5}^{2}}$
&$\arraycolsep=1mm\begin{bmatrix}{1} & {2} & {3} & {6} & {6} \\ {\bar{0}} & {\bar{0}} & {\bar{1}} & {\bar{0}} & {\bar{1}}\end{bmatrix}$&${6}$
&${333}$&${T_{1}^{6}T_{2}^{2}+T_{3}T_{4}+T_{5}^{4}}$
&$\arraycolsep=1mm\begin{bmatrix}{1} & {1} & {4} & {4} & {2} \\ {\bar{0}} & {\bar{1}} & {\bar{0}} & {\bar{0}} & {\bar{1}}\end{bmatrix}$&${8}$
\\
\hline
${294}$&${T_{1}^{5}T_{2}^{2}+T_{3}^{2}T_{4}+T_{5}^{2}}$
&$\arraycolsep=1mm\begin{bmatrix}{2} & {1} & {3} & {6} & {6} \\ {\bar{0}} & {\bar{1}} & {\bar{0}} & {\bar{0}} & {\bar{1}}\end{bmatrix}$&${6}$
&${334}$&${T_{1}^{4}T_{2}^{2}+T_{3}T_{4}+T_{5}^{8}}$
&$\arraycolsep=1mm\begin{bmatrix}{1} & {2} & {4} & {4} & {1} \\ {\bar{1}} & {\bar{1}} & {\bar{0}} & {\bar{0}} & {\bar{0}}\end{bmatrix}$&${8}$
\\
\hline
${295}$&${T_{1}^{4}T_{2}^{4}+T_{3}^{2}T_{4}+T_{5}^{2}}$
&$\arraycolsep=1mm\begin{bmatrix}{1} & {2} & {3} & {6} & {6} \\ {\bar{0}} & {\bar{0}} & {\bar{1}} & {\bar{0}} & {\bar{1}}\end{bmatrix}$&${6}$
&${335}$&${T_{1}^{4}T_{2}^{2}+T_{3}T_{4}+T_{5}^{4}}$
&$\arraycolsep=1mm\begin{bmatrix}{1} & {4} & {4} & {8} & {3} \\ {\bar{0}} & {\bar{1}} & {\bar{0}} & {\bar{0}} & {\bar{1}}\end{bmatrix}$&${8}$
\\
\hline
${296}$&${T_{1}^{2}T_{2}+T_{3}^{2}T_{4}+T_{5}^{2}}$
&$\arraycolsep=1mm\begin{bmatrix}{3} & {2} & {3} & {2} & {4} \\ {\bar{0}} & {\bar{0}} & {\bar{1}} & {\bar{0}} & {\bar{1}}\end{bmatrix}$&${6}$
&${350}$&${T_{1}^{4}T_{2}+T_{3}^{3}T_{4}+T_{5}^{2}}$
&$\arraycolsep=1mm\begin{bmatrix}{1} & {6} & {3} & {1} & {5} \\ {\bar{0}} & {\bar{0}} & {\bar{1}} & {\bar{1}} & {\bar{0}}\end{bmatrix}$&${12}$
\\
\hline
${320}$&${T_{1}^{4}T_{2}+T_{3}^{3}T_{4}^{2}+T_{5}^{2}}$
&$\arraycolsep=1mm\begin{bmatrix}{1} & {4} & {2} & {1} & {4} \\ {\bar{0}} & {\bar{0}} & {\bar{0}} & {\bar{1}} & {\bar{1}}\end{bmatrix}$&${8}$
&${351}$&${T_{1}^{2}T_{2}^{2}+T_{3}^{2}T_{4}+T_{5}^{3}}$
&$\arraycolsep=1mm\begin{bmatrix}{1} & {2} & {1} & {4} & {2} \\ {\bar{1}} & {\bar{1}} & {\bar{0}} & {\bar{0}} & {\bar{0}}\end{bmatrix}$&${12}$
\\
\hline
${321}$&${T_{1}^{3}T_{2}+T_{3}^{3}T_{4}+T_{5}^{2}}$
&$\arraycolsep=1mm\begin{bmatrix}{1} & {1} & {1} & {1} & {2} \\ {\bar{0}} & {\bar{0}} & {\bar{1}} & {\bar{1}} & {\bar{0}}\end{bmatrix}$&${8}$
&${352}$&${T_{1}^{2}T_{2}+T_{3}^{2}T_{4}+T_{5}^{6}}$
&$\arraycolsep=1mm\begin{bmatrix}{2} & {2} & {1} & {4} & {1} \\ {\bar{1}} & {\bar{0}} & {\bar{1}} & {\bar{0}} & {\bar{0}}\end{bmatrix}$&${12}$
\\
\hline
${322}$&${T_{1}^{3}T_{2}+T_{3}^{2}T_{4}+T_{5}^{4}}$
&$\arraycolsep=1mm\begin{bmatrix}{1} & {1} & {1} & {2} & {1} \\ {\bar{1}} & {\bar{1}} & {\bar{0}} & {\bar{0}} & {\bar{0}}\end{bmatrix}$&${8}$
&${353}$&${T_{1}^{2}T_{2}+T_{3}^{2}T_{4}+T_{5}^{4}}$
&$\arraycolsep=1mm\begin{bmatrix}{1} & {6} & {3} & {2} & {2} \\ {\bar{0}} & {\bar{0}} & {\bar{1}} & {\bar{0}} & {\bar{1}}\end{bmatrix}$&${12}$
\\
\hline
${323}$&${T_{1}^{7}T_{2}+T_{3}^{2}T_{4}+T_{5}^{2}}$
&$\arraycolsep=1mm\begin{bmatrix}{1} & {1} & {2} & {4} & {4} \\ {\bar{0}} & {\bar{0}} & {\bar{1}} & {\bar{0}} & {\bar{1}}\end{bmatrix}$&${8}$
&${354}$&${T_{1}^{2}T_{2}+T_{3}^{2}T_{4}+T_{5}^{3}}$
&$\arraycolsep=1mm\begin{bmatrix}{1} & {1} & {1} & {1} & {1} \\ {\bar{1}} & {\bar{0}} & {\bar{1}} & {\bar{0}} & {\bar{0}}\end{bmatrix}$&${12}$
\\
\hline
${324}$&${T_{1}^{6}T_{2}^{2}+T_{3}^{2}T_{4}+T_{5}^{2}}$
&$\arraycolsep=1mm\begin{bmatrix}{1} & {1} & {2} & {4} & {4} \\ {\bar{0}} & {\bar{0}} & {\bar{1}} & {\bar{0}} & {\bar{1}}\end{bmatrix}$&${8}$
&${355}$&${T_{1}^{10}T_{2}^{4}+T_{3}T_{4}+T_{5}^{2}}$
&$\arraycolsep=1mm\begin{bmatrix}{1} & {2} & {6} & {12} & {9} \\ {\bar{0}} & {\bar{1}} & {\bar{0}} & {\bar{0}} & {\bar{1}}\end{bmatrix}$&${12}$
\\
\hline
${325}$&${T_{1}^{6}T_{2}^{2}+T_{3}^{2}T_{4}+T_{5}^{2}}$
&$\arraycolsep=1mm\begin{bmatrix}{1} & {1} & {2} & {4} & {4} \\ {\bar{0}} & {\bar{1}} & {\bar{1}} & {\bar{0}} & {\bar{0}}\end{bmatrix}$&${8}$
&${356}$&${T_{1}^{9}T_{2}+T_{3}T_{4}+T_{5}^{3}}$
&$\arraycolsep=1mm\begin{bmatrix}{1} & {6} & {3} & {12} & {5} \\ {\bar{0}} & {\bar{1}} & {\bar{1}} & {\bar{0}} & {\bar{1}}\end{bmatrix}$&${12}$
\\
\hline
${326}$&${T_{1}^{5}T_{2}^{3}+T_{3}^{2}T_{4}+T_{5}^{2}}$
&$\arraycolsep=1mm\begin{bmatrix}{1} & {1} & {2} & {4} & {4} \\ {\bar{0}} & {\bar{0}} & {\bar{1}} & {\bar{0}} & {\bar{1}}\end{bmatrix}$&${8}$
&${357}$&${T_{1}^{8}T_{2}^{2}+T_{3}T_{4}+T_{5}^{2}}$
&$\arraycolsep=1mm\begin{bmatrix}{2} & {1} & {6} & {12} & {9} \\ {\bar{1}} & {\bar{1}} & {\bar{0}} & {\bar{0}} & {\bar{0}}\end{bmatrix}$&${12}$
\\
\hline
${327}$&${T_{1}^{4}T_{2}^{4}+T_{3}^{2}T_{4}+T_{5}^{2}}$
&$\arraycolsep=1mm\begin{bmatrix}{1} & {1} & {2} & {4} & {4} \\ {\bar{0}} & {\bar{0}} & {\bar{1}} & {\bar{0}} & {\bar{1}}\end{bmatrix}$&${8}$
&${358}$&${T_{1}^{5}T_{2}+T_{3}T_{4}+T_{5}^{3}}$
&$\arraycolsep=1mm\begin{bmatrix}{1} & {1} & {2} & {4} & {2} \\ {\bar{0}} & {\bar{1}} & {\bar{1}} & {\bar{0}} & {\bar{1}}\end{bmatrix}$&${12}$
\\
\hline
${328}$&${T_{1}^{4}T_{2}^{4}+T_{3}^{2}T_{4}+T_{5}^{2}}$
&$\arraycolsep=1mm\begin{bmatrix}{1} & {1} & {2} & {4} & {4} \\ {\bar{1}} & {\bar{0}} & {\bar{0}} & {\bar{0}} & {\bar{1}}\end{bmatrix}$&${8}$
&${359}$&${T_{1}^{4}T_{2}^{4}+T_{3}T_{4}+T_{5}^{2}}$
&$\arraycolsep=1mm\begin{bmatrix}{1} & {3} & {4} & {12} & {8} \\ {\bar{0}} & {\bar{1}} & {\bar{0}} & {\bar{0}} & {\bar{1}}\end{bmatrix}$&${12}$
\\
\hline
${329}$&${T_{1}^{4}T_{2}+T_{3}^{2}T_{4}+T_{5}^{2}}$
&$\arraycolsep=1mm\begin{bmatrix}{1} & {2} & {2} & {2} & {3} \\ {\bar{0}} & {\bar{0}} & {\bar{1}} & {\bar{0}} & {\bar{1}}\end{bmatrix}$&${8}$
&${360}$&${T_{1}^{4}T_{2}^{2}+T_{3}T_{4}+T_{5}^{2}}$
&$\arraycolsep=1mm\begin{bmatrix}{2} & {3} & {2} & {12} & {7} \\ {\bar{1}} & {\bar{0}} & {\bar{0}} & {\bar{0}} & {\bar{1}}\end{bmatrix}$&${12}$
\\
\hline
${330}$&${T_{1}^{2}T_{2}^{2}+T_{3}^{2}T_{4}+T_{5}^{4}}$
&$\arraycolsep=1mm\begin{bmatrix}{1} & {1} & {1} & {2} & {1} \\ {\bar{0}} & {\bar{1}} & {\bar{1}} & {\bar{0}} & {\bar{0}}\end{bmatrix}$&${8}$
&${361}$&${T_{1}^{4}T_{2}+T_{3}T_{4}+T_{5}^{2}}$
&$\arraycolsep=1mm\begin{bmatrix}{1} & {2} & {3} & {3} & {3} \\ {\bar{0}} & {\bar{0}} & {\bar{1}} & {\bar{1}} & {\bar{0}}\end{bmatrix}$&${12}$
\\
\hline
${331}$&${T_{1}^{2}T_{2}+T_{3}^{2}T_{4}^{2}+T_{5}^{2}}$
&$\arraycolsep=1mm\begin{bmatrix}{2} & {2} & {2} & {1} & {3} \\ {\bar{1}} & {\bar{0}} & {\bar{0}} & {\bar{0}} & {\bar{1}}\end{bmatrix}$&${8}$
&${362}$&${T_{1}^{3}T_{2}+T_{3}T_{4}+T_{5}^{3}}$
&$\arraycolsep=1mm\begin{bmatrix}{1} & {3} & {3} & {3} & {2} \\ {\bar{0}} & {\bar{0}} & {\bar{1}} & {\bar{1}} & {\bar{0}}\end{bmatrix}$&${12}$
\\
\hline
\end{tabu}\vspace{6mm}
\end{table}

\begin{table}[h!]\centering\setlength{\tabcolsep}{3.5pt}
\renewcommand{\arraystretch}{1.3}\small\tabulinesep=1mm
\begin{tabu}{c|c|c|c||c|c|c|c}
\hline
ID&relations&gd-matrix&$-\mathcal{K}^3$&ID&relations&gd-matrix&$-\mathcal{K}^3$
\\
\hline
${363}$&${T_{1}^{3}T_{2}+T_{3}T_{4}+T_{5}^{3}}$
&$\arraycolsep=1mm\begin{bmatrix}{1} & {3} & {3} & {3} & {2} \\ {\bar{0}} & {\bar{1}} & {\bar{0}} & {\bar{1}} & {\bar{1}}\end{bmatrix}$&${12}$
&${391}$&${T_{1}^{2}T_{2}+T_{3}T_{4}+T_{5}^{4}}$
&$\arraycolsep=1mm\begin{bmatrix}{1} & {2} & {2} & {2} & {1} \\ {\bar{0}} & {\bar{0}} & {\bar{1}} & {\bar{1}} & {\bar{0}}\end{bmatrix}$&${16}$
\\
\hline
${364}$&${T_{1}^{2}T_{2}^{2}+T_{3}T_{4}+T_{5}^{4}}$
&$\arraycolsep=1mm\begin{bmatrix}{1} & {3} & {2} & {6} & {2} \\ {\bar{0}} & {\bar{1}} & {\bar{0}} & {\bar{0}} & {\bar{1}}\end{bmatrix}$&${12}$
&${392}$&${T_{1}^{2}T_{2}+T_{3}T_{4}+T_{5}^{2}}$
&$\arraycolsep=1mm\begin{bmatrix}{1} & {4} & {2} & {4} & {3} \\ {\bar{0}} & {\bar{0}} & {\bar{1}} & {\bar{1}} & {\bar{0}}\end{bmatrix}$&${16}$
\\
\hline
${381}$&${T_{1}^{5}T_{2}+T_{3}^{2}T_{4}+T_{5}^{2}}$
&$\arraycolsep=1mm\begin{bmatrix}{1} & {1} & {1} & {4} & {3} \\ {\bar{0}} & {\bar{0}} & {\bar{1}} & {\bar{0}} & {\bar{1}}\end{bmatrix}$&${16}$
&${400}$&${T_{1}^{6}T_{2}+T_{3}^{2}T_{4}+T_{5}^{2}}$
&$\arraycolsep=1mm\begin{bmatrix}{1} & {2} & {1} & {6} & {4} \\ {\bar{0}} & {\bar{0}} & {\bar{1}} & {\bar{0}} & {\bar{1}}\end{bmatrix}$&${18}$
\\
\hline
${382}$&${T_{1}^{4}T_{2}^{2}+T_{3}^{2}T_{4}+T_{5}^{2}}$
&$\arraycolsep=1mm\begin{bmatrix}{1} & {1} & {1} & {4} & {3} \\ {\bar{0}} & {\bar{0}} & {\bar{1}} & {\bar{0}} & {\bar{1}}\end{bmatrix}$&${16}$
&${401}$&${T_{1}^{4}T_{2}^{2}+T_{3}^{2}T_{4}+T_{5}^{2}}$
&$\arraycolsep=1mm\begin{bmatrix}{1} & {2} & {1} & {6} & {4} \\ {\bar{0}} & {\bar{0}} & {\bar{1}} & {\bar{0}} & {\bar{1}}\end{bmatrix}$&${18}$
\\
\hline
${383}$&${T_{1}^{4}T_{2}^{2}+T_{3}^{2}T_{4}+T_{5}^{2}}$
&$\arraycolsep=1mm\begin{bmatrix}{1} & {1} & {1} & {4} & {3} \\ {\bar{0}} & {\bar{1}} & {\bar{1}} & {\bar{0}} & {\bar{0}}\end{bmatrix}$&${16}$
&${402}$&${T_{1}^{3}T_{2}^{2}+T_{3}^{2}T_{4}+T_{5}^{2}}$
&$\arraycolsep=1mm\begin{bmatrix}{2} & {1} & {1} & {6} & {4} \\ {\bar{0}} & {\bar{1}} & {\bar{0}} & {\bar{0}} & {\bar{1}}\end{bmatrix}$&${18}$
\\
\hline
${384}$&${T_{1}^{4}T_{2}+T_{3}^{2}T_{4}+T_{5}^{2}}$
&$\arraycolsep=1mm\begin{bmatrix}{2} & {2} & {1} & {8} & {5} \\ {\bar{1}} & {\bar{0}} & {\bar{0}} & {\bar{0}} & {\bar{1}}\end{bmatrix}$&${16}$
&${403}$&${T_{1}^{4}T_{2}^{4}+T_{3}T_{4}+T_{5}^{2}}$
&$\arraycolsep=1mm\begin{bmatrix}{1} & {1} & {2} & {6} & {4} \\ {\bar{1}} & {\bar{0}} & {\bar{0}} & {\bar{0}} & {\bar{1}}\end{bmatrix}$&${18}$
\\
\hline
${385}$&${T_{1}^{3}T_{2}^{3}+T_{3}^{2}T_{4}+T_{5}^{2}}$
&$\arraycolsep=1mm\begin{bmatrix}{1} & {1} & {1} & {4} & {3} \\ {\bar{0}} & {\bar{0}} & {\bar{1}} & {\bar{0}} & {\bar{1}}\end{bmatrix}$&${16}$
&${405}$&${T_{1}^{3}T_{2}+T_{3}T_{4}+T_{5}^{5}}$
&$\arraycolsep=1mm\begin{bmatrix}{1} & {2} & {1} & {4} & {1} \\ {\bar{1}} & {\bar{1}} & {\bar{0}} & {\bar{0}} & {\bar{0}}\end{bmatrix}$&${20}$
\\
\hline
${386}$&${T_{1}^{3}T_{2}^{2}+T_{3}^{2}T_{4}+T_{5}^{2}}$
&$\arraycolsep=1mm\begin{bmatrix}{2} & {2} & {1} & {8} & {5} \\ {\bar{0}} & {\bar{1}} & {\bar{0}} & {\bar{0}} & {\bar{1}}\end{bmatrix}$&${16}$
&${408}$&${T_{1}^{5}T_{2}+T_{3}T_{4}+T_{5}^{3}}$
&$\arraycolsep=1mm\begin{bmatrix}{1} & {4} & {1} & {8} & {3} \\ {\bar{0}} & {\bar{1}} & {\bar{1}} & {\bar{0}} & {\bar{1}}\end{bmatrix}$&${24}$
\\
\hline
${387}$&${T_{1}^{2}T_{2}+T_{3}^{2}T_{4}+T_{5}^{2}}$
&$\arraycolsep=1mm\begin{bmatrix}{1} & {2} & {1} & {2} & {2} \\ {\bar{0}} & {\bar{0}} & {\bar{1}} & {\bar{0}} & {\bar{1}}\end{bmatrix}$&${16}$
&${409}$&${T_{1}T_{2}+T_{3}T_{4}+T_{5}^{3}}$
&$\arraycolsep=1mm\begin{bmatrix}{1} & {2} & {1} & {2} & {1} \\ {\bar{1}} & {\bar{1}} & {\bar{1}} & {\bar{1}} & {\bar{0}}\end{bmatrix}$&${24}$
\\
\hline
${388}$&${T_{1}^{8}T_{2}^{4}+T_{3}T_{4}+T_{5}^{2}}$
&$\arraycolsep=1mm\begin{bmatrix}{1} & {1} & {4} & {8} & {6} \\ {\bar{1}} & {\bar{0}} & {\bar{0}} & {\bar{0}} & {\bar{1}}\end{bmatrix}$&${16}$
&${410}$&${T_{1}T_{2}+T_{3}T_{4}+T_{5}^{3}}$
&$\arraycolsep=1mm\begin{bmatrix}{2} & {1} & {1} & {2} & {1} \\ {\bar{1}} & {\bar{1}} & {\bar{0}} & {\bar{0}} & {\bar{0}}\end{bmatrix}$&${24}$
\\
\hline
${389}$&${T_{1}^{4}T_{2}^{2}+T_{3}T_{4}+T_{5}^{2}}$
&$\arraycolsep=1mm\begin{bmatrix}{2} & {1} & {2} & {8} & {5} \\ {\bar{1}} & {\bar{0}} & {\bar{0}} & {\bar{0}} & {\bar{1}}\end{bmatrix}$&${16}$
&${411}$&${T_{1}T_{2}+T_{3}T_{4}+T_{5}^{2}}$
&$\arraycolsep=1mm\begin{bmatrix}{1} & {3} & {1} & {3} & {2} \\ {\bar{0}} & {\bar{0}} & {\bar{1}} & {\bar{1}} & {\bar{0}}\end{bmatrix}$&${24}$
\\
\hline
${390}$&${T_{1}^{3}T_{2}+T_{3}T_{4}+T_{5}^{2}}$
&$\arraycolsep=1mm\begin{bmatrix}{1} & {1} & {2} & {2} & {2} \\ {\bar{0}} & {\bar{0}} & {\bar{1}} & {\bar{1}} & {\bar{0}}\end{bmatrix}$&${16}$
\\
\cline{1-4}
\end{tabu}
\end{table}
\end{class-list}

\begin{class-list}
Non-toric, $\QQ$-factorial, Gorenstein, log terminal Fano
threefolds of Picard number one with an effective two-torus action:
Specifying data for divisor class group $\ZZ\oplus \ZZ/2\ZZ$ and
format $(2,2,1,1,0)$.
\begin{table}[h!]\centering\setlength{\tabcolsep}{6.5pt}
\renewcommand{\arraystretch}{1.3}\small\tabulinesep=1mm
\begin{tabu}{c|c|c|c}
\hline
ID&relations&gd-matrix&$-\mathcal{K}^3$
\\
\hline
${269}$&
$\setlength{\arraycolsep}{1pt}\begin{array}{r}{T_{1}^{2}T_{2}+T_{3}T_{4}+T_{5}^{4}},\\ {\lambda T_{3}T_{4}+T_{5}^{4}+T_{6}^{2}}\end{array}$
&$\arraycolsep=1mm\begin{bmatrix}{1} & {2} & {2} & {2} & {1} & {2} \\ {\bar{0}} & {\bar{0}} & {\bar{0}} & {\bar{0}} & {\bar{1}} & {\bar{1}}\end{bmatrix}$&${4}$
\\
\hline
${270}$&
$\setlength{\arraycolsep}{1pt}\begin{array}{r}{T_{1}^{2}T_{2}+T_{3}T_{4}+T_{5}^{2}},\\ {\lambda T_{3}T_{4}+T_{5}^{2}+T_{6}^{2}}\end{array}$
&$\arraycolsep=1mm\begin{bmatrix}{2} & {2} & {2} & {4} & {3} & {3} \\ {\bar{1}} & {\bar{0}} & {\bar{0}} & {\bar{0}} & {\bar{0}} & {\bar{1}}\end{bmatrix}$&${4}$
\\
\hline
${349}$&
$\setlength{\arraycolsep}{1pt}\begin{array}{r}{T_{1}T_{2}+T_{3}T_{4}+T_{5}^{4}},\\ {\lambda T_{3}T_{4}+T_{5}^{4}+T_{6}^{2}}\end{array}$
&$\arraycolsep=1mm\begin{bmatrix}{1} & {3} & {1} & {3} & {1} & {2} \\ {\bar{1}} & {\bar{1}} & {\bar{1}} & {\bar{1}} & {\bar{0}} & {\bar{1}}\end{bmatrix}$&${12}$
\\
\hline
\end{tabu}
\end{table}
\end{class-list}

\begin{class-list}
Non-toric, $\QQ$-factorial, Gorenstein, log terminal Fano
threefolds of Picard number one with an effective two-torus action:
Specifying data for divisor class group $\ZZ\oplus \ZZ/2\ZZ$ and
format $(3,1,1,0)$.
\begin{table}[h!]\centering\setlength{\tabcolsep}{3.5pt}
\renewcommand{\arraystretch}{1.3}\small\tabulinesep=1mm
\begin{tabu}{c|c|c|c||c|c|c|c}
\hline
ID&relations&gd-matrix&$-\mathcal{K}^3$&ID&relations&gd-matrix&$-\mathcal{K}^3$
\\
\hline
${248}$&${T_{1}^{3}T_{2}^{3}T_{3}+T_{4}^{4}+T_{5}^{2}}$
&$\arraycolsep=1mm\begin{bmatrix}{1} & {1} & {2} & {2} & {4} \\ {\bar{0}} & {\bar{0}} & {\bar{0}} & {\bar{1}} & {\bar{1}}\end{bmatrix}$&${2}$
&${312}$&${T_{1}^{2}T_{2}^{2}T_{3}+T_{4}^{4}+T_{5}^{2}}$
&$\arraycolsep=1mm\begin{bmatrix}{1} & {1} & {4} & {2} & {4} \\ {\bar{0}} & {\bar{0}} & {\bar{0}} & {\bar{1}} & {\bar{1}}\end{bmatrix}$&${8}$
\\
\hline
${249}$&${T_{1}^{3}T_{2}^{2}T_{3}+T_{4}^{4}+T_{5}^{2}}$
&$\arraycolsep=1mm\begin{bmatrix}{1} & {2} & {1} & {2} & {4} \\ {\bar{0}} & {\bar{0}} & {\bar{0}} & {\bar{1}} & {\bar{1}}\end{bmatrix}$&${2}$
&${313}$&${T_{1}^{2}T_{2}^{2}T_{3}+T_{4}^{4}+T_{5}^{2}}$
&$\arraycolsep=1mm\begin{bmatrix}{1} & {1} & {4} & {2} & {4} \\ {\bar{1}} & {\bar{0}} & {\bar{0}} & {\bar{0}} & {\bar{1}}\end{bmatrix}$&${8}$
\\
\hline
${250}$&${T_{1}^{3}T_{2}T_{3}+T_{4}^{4}+T_{5}^{2}}$
&$\arraycolsep=1mm\begin{bmatrix}{2} & {1} & {1} & {2} & {4} \\ {\bar{0}} & {\bar{0}} & {\bar{0}} & {\bar{1}} & {\bar{1}}\end{bmatrix}$&${2}$
&${314}$&${T_{1}^{2}T_{2}T_{3}+T_{4}^{8}+T_{5}^{2}}$
&$\arraycolsep=1mm\begin{bmatrix}{1} & {2} & {4} & {1} & {4} \\ {\bar{0}} & {\bar{0}} & {\bar{0}} & {\bar{1}} & {\bar{1}}\end{bmatrix}$&${8}$
\\
\hline
${251}$&${T_{1}^{2}T_{2}^{2}T_{3}+T_{4}^{8}+T_{5}^{2}}$
&$\arraycolsep=1mm\begin{bmatrix}{2} & {1} & {2} & {1} & {4} \\ {\bar{0}} & {\bar{0}} & {\bar{0}} & {\bar{1}} & {\bar{1}}\end{bmatrix}$&${2}$
&${315}$&${T_{1}^{2}T_{2}T_{3}+T_{4}^{4}+T_{5}^{2}}$
&$\arraycolsep=1mm\begin{bmatrix}{1} & {1} & {1} & {1} & {2} \\ {\bar{0}} & {\bar{0}} & {\bar{0}} & {\bar{1}} & {\bar{1}}\end{bmatrix}$&${8}$
\\
\hline
${263}$&${T_{1}^{4}T_{2}^{4}T_{3}+T_{4}^{3}+T_{5}^{2}}$
&$\arraycolsep=1mm\begin{bmatrix}{1} & {1} & {4} & {4} & {6} \\ {\bar{1}} & {\bar{0}} & {\bar{0}} & {\bar{0}} & {\bar{1}}\end{bmatrix}$&${4}$
&${316}$&${T_{1}^{2}T_{2}T_{3}+T_{4}^{4}+T_{5}^{2}}$
&$\arraycolsep=1mm\begin{bmatrix}{1} & {1} & {1} & {1} & {2} \\ {\bar{0}} & {\bar{1}} & {\bar{1}} & {\bar{1}} & {\bar{1}}\end{bmatrix}$&${8}$
\\
\hline
${264}$&${T_{1}^{4}T_{2}T_{3}+T_{4}^{3}+T_{5}^{2}}$
&$\arraycolsep=1mm\begin{bmatrix}{1} & {1} & {1} & {2} & {3} \\ {\bar{0}} & {\bar{1}} & {\bar{1}} & {\bar{0}} & {\bar{0}}\end{bmatrix}$&${4}$
&${317}$&${T_{1}^{2}T_{2}T_{3}+T_{4}^{3}+T_{5}^{2}}$
&$\arraycolsep=1mm\begin{bmatrix}{1} & {2} & {2} & {2} & {3} \\ {\bar{0}} & {\bar{1}} & {\bar{1}} & {\bar{0}} & {\bar{0}}\end{bmatrix}$&${8}$
\\
\hline
${265}$&${T_{1}^{3}T_{2}^{2}T_{3}+T_{4}^{3}+T_{5}^{2}}$
&$\arraycolsep=1mm\begin{bmatrix}{1} & {1} & {1} & {2} & {3} \\ {\bar{1}} & {\bar{0}} & {\bar{1}} & {\bar{0}} & {\bar{0}}\end{bmatrix}$&${4}$
&${318}$&${T_{1}^{4}T_{2}^{2}T_{3}+T_{4}^{2}+T_{5}^{2}}$
&$\arraycolsep=1mm\begin{bmatrix}{1} & {1} & {2} & {4} & {4} \\ {\bar{0}} & {\bar{1}} & {\bar{0}} & {\bar{1}} & {\bar{0}}\end{bmatrix}$&${8}$
\\
\hline
${266}$&${T_{1}^{2}T_{2}^{2}T_{3}+T_{4}^{6}+T_{5}^{2}}$
&$\arraycolsep=1mm\begin{bmatrix}{1} & {1} & {2} & {1} & {3} \\ {\bar{1}} & {\bar{0}} & {\bar{0}} & {\bar{0}} & {\bar{1}}\end{bmatrix}$&${4}$
&${319}$&${T_{1}^{3}T_{2}^{2}T_{3}+T_{4}^{2}+T_{5}^{2}}$
&$\arraycolsep=1mm\begin{bmatrix}{1} & {2} & {1} & {4} & {4} \\ {\bar{0}} & {\bar{1}} & {\bar{0}} & {\bar{1}} & {\bar{0}}\end{bmatrix}$&${8}$
\\
\hline
${267}$&${T_{1}^{2}T_{2}^{2}T_{3}+T_{4}^{5}+T_{5}^{2}}$
&$\arraycolsep=1mm\begin{bmatrix}{2} & {1} & {4} & {2} & {5} \\ {\bar{1}} & {\bar{0}} & {\bar{0}} & {\bar{0}} & {\bar{1}}\end{bmatrix}$&${4}$
&${337}$&${T_{1}T_{2}T_{3}+T_{4}^{8}+T_{5}^{2}}$
&$\arraycolsep=1mm\begin{bmatrix}{1} & {5} & {10} & {2} & {8} \\ {\bar{0}} & {\bar{0}} & {\bar{0}} & {\bar{1}} & {\bar{1}}\end{bmatrix}$&${10}$
\\
\hline
${268}$&${T_{1}^{2}T_{2}T_{3}+T_{4}^{10}+T_{5}^{2}}$
&$\arraycolsep=1mm\begin{bmatrix}{2} & {2} & {4} & {1} & {5} \\ {\bar{1}} & {\bar{0}} & {\bar{0}} & {\bar{0}} & {\bar{1}}\end{bmatrix}$&${4}$
&${347}$&${T_{1}^{2}T_{2}^{2}T_{3}+T_{4}^{3}+T_{5}^{2}}$
&$\arraycolsep=1mm\begin{bmatrix}{2} & {1} & {12} & {6} & {9} \\ {\bar{1}} & {\bar{0}} & {\bar{0}} & {\bar{0}} & {\bar{1}}\end{bmatrix}$&${12}$
\\
\hline
${288}$&${T_{1}^{2}T_{2}^{2}T_{3}+T_{4}^{4}+T_{5}^{2}}$
&$\arraycolsep=1mm\begin{bmatrix}{2} & {1} & {6} & {3} & {6} \\ {\bar{0}} & {\bar{0}} & {\bar{0}} & {\bar{1}} & {\bar{1}}\end{bmatrix}$&${6}$
&${348}$&${T_{1}T_{2}T_{3}+T_{4}^{4}+T_{5}^{2}}$
&$\arraycolsep=1mm\begin{bmatrix}{1} & {3} & {12} & {4} & {8} \\ {\bar{0}} & {\bar{0}} & {\bar{0}} & {\bar{1}} & {\bar{1}}\end{bmatrix}$&${12}$
\\
\hline
${289}$&${T_{1}^{4}T_{2}^{2}T_{3}+T_{4}^{2}+T_{5}^{2}}$
&$\arraycolsep=1mm\begin{bmatrix}{1} & {3} & {2} & {6} & {6} \\ {\bar{0}} & {\bar{1}} & {\bar{0}} & {\bar{1}} & {\bar{0}}\end{bmatrix}$&${6}$
&${380}$&${T_{1}^{2}T_{2}^{2}T_{3}+T_{4}^{3}+T_{5}^{2}}$
&$\arraycolsep=1mm\begin{bmatrix}{1} & {1} & {8} & {4} & {6} \\ {\bar{1}} & {\bar{0}} & {\bar{0}} & {\bar{0}} & {\bar{1}}\end{bmatrix}$&${16}$
\\
\hline
${290}$&${T_{1}T_{2}T_{3}+T_{4}^{4}+T_{5}^{2}}$
&$\arraycolsep=1mm\begin{bmatrix}{2} & {3} & {3} & {2} & {4} \\ {\bar{0}} & {\bar{0}} & {\bar{0}} & {\bar{1}} & {\bar{1}}\end{bmatrix}$&${6}$
&${399}$&${T_{1}T_{2}T_{3}+T_{4}^{4}+T_{5}^{2}}$
&$\arraycolsep=1mm\begin{bmatrix}{1} & {1} & {6} & {2} & {4} \\ {\bar{0}} & {\bar{0}} & {\bar{0}} & {\bar{1}} & {\bar{1}}\end{bmatrix}$&${18}$\\
\hline
${311}$&${T_{1}^{3}T_{2}T_{3}+T_{4}^{4}+T_{5}^{2}}$
&$\arraycolsep=1mm\begin{bmatrix}{1} & {4} & {1} & {2} & {4} \\ {\bar{0}} & {\bar{0}} & {\bar{0}} & {\bar{1}} & {\bar{1}}\end{bmatrix}$&${8}$
\\
\cline{1-4}
\end{tabu}
\end{table}
\end{class-list}

\newpage

\begin{class-list}
Non-toric, $\QQ$-factorial, Gorenstein, log terminal Fano
threefolds of Picard number one with an effective two-torus action:
Specifying data for divisor class group $\ZZ\oplus \ZZ/2\ZZ$ and
format $(2,1,1,1)$.
\begin{table}[h!]\centering\setlength{\tabcolsep}{3.5pt}
\renewcommand{\arraystretch}{1.25}\small\tabulinesep=1mm
\begin{tabu}{c|c|c|c||c|c|c|c}
\hline
ID&relations&gd-matrix&$-\mathcal{K}^3$&ID&relations&gd-matrix&$-\mathcal{K}^3$
\\
\hline
${246}$&${T_{1}^{3}T_{2}^{2}+T_{3}^{4}+T_{4}^{2}}$
&$\arraycolsep=1mm\begin{bmatrix}{2} & {1} & {2} & {4} & {1} \\ {\bar{0}} & {\bar{1}} & {\bar{0}} & {\bar{1}} & {\bar{0}}\end{bmatrix}$&${2}$
&${285}$&${T_{1}^{5}T_{2}^{2}+T_{3}^{2}+T_{4}^{2}}$
&$\arraycolsep=1mm\begin{bmatrix}{2} & {1} & {6} & {6} & {3} \\ {\bar{0}} & {\bar{1}} & {\bar{1}} & {\bar{0}} & {\bar{0}}\end{bmatrix}$&${6}$
\\
\hline
${247}$&${T_{1}^{3}T_{2}^{2}+T_{3}^{4}+T_{4}^{2}}$
&$\arraycolsep=1mm\begin{bmatrix}{2} & {1} & {2} & {4} & {1} \\ {\bar{0}} & {\bar{0}} & {\bar{1}} & {\bar{1}} & {\bar{0}}\end{bmatrix}$&${2}$
&${286}$&${T_{1}^{3}T_{2}^{2}+T_{3}^{2}+T_{4}^{2}}$
&$\arraycolsep=1mm\begin{bmatrix}{2} & {3} & {6} & {6} & {1} \\ {\bar{0}} & {\bar{1}} & {\bar{1}} & {\bar{0}} & {\bar{0}}\end{bmatrix}$&${6}$
\\
\hline
${256}$&${T_{1}^{5}T_{2}+T_{3}^{3}+T_{4}^{2}}$
&$\arraycolsep=1mm\begin{bmatrix}{1} & {1} & {2} & {3} & {1} \\ {\bar{0}} & {\bar{0}} & {\bar{0}} & {\bar{1}} & {\bar{1}}\end{bmatrix}$&${4}$
&${287}$&${T_{1}T_{2}+T_{3}^{6}+T_{4}^{4}}$
&$\arraycolsep=1mm\begin{bmatrix}{6} & {6} & {2} & {3} & {1} \\ {\bar{0}} & {\bar{0}} & {\bar{1}} & {\bar{1}} & {\bar{0}}\end{bmatrix}$&${6}$
\\
\hline
${257}$&${T_{1}^{4}T_{2}^{2}+T_{3}^{3}+T_{4}^{2}}$
&$\arraycolsep=1mm\begin{bmatrix}{1} & {4} & {4} & {6} & {1} \\ {\bar{1}} & {\bar{0}} & {\bar{0}} & {\bar{1}} & {\bar{0}}\end{bmatrix}$&${4}$
&${297}$&${T_{1}^{4}T_{2}+T_{3}^{3}+T_{4}^{2}}$
&$\arraycolsep=1mm\begin{bmatrix}{1} & {2} & {2} & {3} & {2} \\ {\bar{0}} & {\bar{0}} & {\bar{0}} & {\bar{1}} & {\bar{1}}\end{bmatrix}$&${8}$
\\
\hline
${258}$&${T_{1}^{4}T_{2}^{2}+T_{3}^{3}+T_{4}^{2}}$
&$\arraycolsep=1mm\begin{bmatrix}{1} & {1} & {2} & {3} & {1} \\ {\bar{1}} & {\bar{1}} & {\bar{0}} & {\bar{0}} & {\bar{0}}\end{bmatrix}$&${4}$
&${298}$&${T_{1}^{3}T_{2}+T_{3}^{4}+T_{4}^{2}}$
&$\arraycolsep=1mm\begin{bmatrix}{1} & {1} & {1} & {2} & {1} \\ {\bar{0}} & {\bar{0}} & {\bar{1}} & {\bar{1}} & {\bar{0}}\end{bmatrix}$&${8}$
\\
\hline
${259}$&${T_{1}^{4}T_{2}^{2}+T_{3}^{3}+T_{4}^{2}}$
&$\arraycolsep=1mm\begin{bmatrix}{1} & {1} & {2} & {3} & {1} \\ {\bar{1}} & {\bar{0}} & {\bar{0}} & {\bar{1}} & {\bar{0}}\end{bmatrix}$&${4}$
&${299}$&${T_{1}^{3}T_{2}+T_{3}^{4}+T_{4}^{2}}$
&$\arraycolsep=1mm\begin{bmatrix}{1} & {1} & {1} & {2} & {1} \\ {\bar{1}} & {\bar{1}} & {\bar{1}} & {\bar{1}} & {\bar{0}}\end{bmatrix}$&${8}$
\\
\hline
${260}$&${T_{1}^{2}T_{2}^{2}+T_{3}^{5}+T_{4}^{2}}$
&$\arraycolsep=1mm\begin{bmatrix}{4} & {1} & {2} & {5} & {2} \\ {\bar{0}} & {\bar{0}} & {\bar{0}} & {\bar{1}} & {\bar{1}}\end{bmatrix}$&${4}$
&${300}$&${T_{1}^{2}T_{2}^{2}+T_{3}^{3}+T_{4}^{2}}$
&$\arraycolsep=1mm\begin{bmatrix}{2} & {1} & {2} & {3} & {2} \\ {\bar{0}} & {\bar{0}} & {\bar{0}} & {\bar{1}} & {\bar{1}}\end{bmatrix}$&${8}$
\\
\hline
${261}$&${T_{1}^{2}T_{2}+T_{3}^{12}+T_{4}^{2}}$
&$\arraycolsep=1mm\begin{bmatrix}{4} & {4} & {1} & {6} & {1} \\ {\bar{0}} & {\bar{0}} & {\bar{1}} & {\bar{1}} & {\bar{0}}\end{bmatrix}$&${4}$
&${301}$&${T_{1}^{2}T_{2}^{2}+T_{3}^{3}+T_{4}^{2}}$
&$\arraycolsep=1mm\begin{bmatrix}{2} & {1} & {2} & {3} & {2} \\ {\bar{1}} & {\bar{0}} & {\bar{0}} & {\bar{0}} & {\bar{1}}\end{bmatrix}$&${8}$
\\
\hline
${262}$&${T_{1}^{2}T_{2}+T_{3}^{10}+T_{4}^{2}}$
&$\arraycolsep=1mm\begin{bmatrix}{4} & {2} & {1} & {5} & {2} \\ {\bar{0}} & {\bar{0}} & {\bar{0}} & {\bar{1}} & {\bar{1}}\end{bmatrix}$&${4}$
&${302}$&${T_{1}^{2}T_{2}+T_{3}^{8}+T_{4}^{2}}$
&$\arraycolsep=1mm\begin{bmatrix}{2} & {4} & {1} & {4} & {1} \\ {\bar{1}} & {\bar{0}} & {\bar{0}} & {\bar{1}} & {\bar{0}}\end{bmatrix}$&${8}$
\\
\hline
${277}$&${T_{1}^{3}T_{2}+T_{3}^{4}+T_{4}^{2}}$
&$\arraycolsep=1mm\begin{bmatrix}{2} & {6} & {3} & {6} & {1} \\ {\bar{0}} & {\bar{0}} & {\bar{1}} & {\bar{1}} & {\bar{0}}\end{bmatrix}$&${6}$
&${303}$&${T_{1}^{2}T_{2}+T_{3}^{8}+T_{4}^{2}}$
&$\arraycolsep=1mm\begin{bmatrix}{2} & {4} & {1} & {4} & {1} \\ {\bar{0}} & {\bar{0}} & {\bar{1}} & {\bar{1}} & {\bar{0}}\end{bmatrix}$&${8}$
\\
\hline
${278}$&${T_{1}^{2}T_{2}+T_{3}^{12}+T_{4}^{2}}$
&$\arraycolsep=1mm\begin{bmatrix}{3} & {6} & {1} & {6} & {2} \\ {\bar{0}} & {\bar{0}} & {\bar{1}} & {\bar{1}} & {\bar{0}}\end{bmatrix}$&${6}$
&${304}$&${T_{1}^{2}T_{2}+T_{3}^{6}+T_{4}^{2}}$
&$\arraycolsep=1mm\begin{bmatrix}{2} & {2} & {1} & {3} & {2} \\ {\bar{1}} & {\bar{0}} & {\bar{0}} & {\bar{1}} & {\bar{0}}\end{bmatrix}$&${8}$
\\
\hline
${279}$&${T_{1}^{2}T_{2}+T_{3}^{6}+T_{4}^{2}}$
&$\arraycolsep=1mm\begin{bmatrix}{3} & {6} & {2} & {6} & {1} \\ {\bar{1}} & {\bar{0}} & {\bar{0}} & {\bar{1}} & {\bar{0}}\end{bmatrix}$&${6}$
&${305}$&${T_{1}^{2}T_{2}+T_{3}^{6}+T_{4}^{2}}$
&$\arraycolsep=1mm\begin{bmatrix}{2} & {2} & {1} & {3} & {2} \\ {\bar{0}} & {\bar{0}} & {\bar{0}} & {\bar{1}} & {\bar{1}}\end{bmatrix}$&${8}$
\\
\hline
${280}$&${T_{1}^{2}T_{2}+T_{3}^{6}+T_{4}^{2}}$
&$\arraycolsep=1mm\begin{bmatrix}{3} & {6} & {2} & {6} & {1} \\ {\bar{1}} & {\bar{0}} & {\bar{1}} & {\bar{0}} & {\bar{0}}\end{bmatrix}$&${6}$
&${306}$&${T_{1}^{2}T_{2}+T_{3}^{3}+T_{4}^{2}}$
&$\arraycolsep=1mm\begin{bmatrix}{2} & {2} & {2} & {3} & {1} \\ {\bar{1}} & {\bar{0}} & {\bar{0}} & {\bar{1}} & {\bar{0}}\end{bmatrix}$&${8}$
\\
\hline
${281}$&${T_{1}^{2}T_{2}+T_{3}^{4}+T_{4}^{2}}$
&$\arraycolsep=1mm\begin{bmatrix}{3} & {2} & {2} & {4} & {3} \\ {\bar{0}} & {\bar{0}} & {\bar{0}} & {\bar{1}} & {\bar{1}}\end{bmatrix}$&${6}$
&${307}$&${T_{1}^{7}T_{2}+T_{3}^{2}+T_{4}^{2}}$
&$\arraycolsep=1mm\begin{bmatrix}{1} & {1} & {4} & {4} & {2} \\ {\bar{0}} & {\bar{0}} & {\bar{1}} & {\bar{0}} & {\bar{1}}\end{bmatrix}$&${8}$
\\
\hline
${282}$&${T_{1}^{2}T_{2}+T_{3}^{4}+T_{4}^{2}}$
&$\arraycolsep=1mm\begin{bmatrix}{3} & {2} & {2} & {4} & {3} \\ {\bar{1}} & {\bar{0}} & {\bar{1}} & {\bar{1}} & {\bar{1}}\end{bmatrix}$&${6}$
&${308}$&${T_{1}^{6}T_{2}+T_{3}^{2}+T_{4}^{2}}$
&$\arraycolsep=1mm\begin{bmatrix}{1} & {2} & {4} & {4} & {1} \\ {\bar{1}} & {\bar{0}} & {\bar{1}} & {\bar{0}} & {\bar{0}}\end{bmatrix}$&${8}$
\\
\hline
${283}$&${T_{1}^{10}T_{2}+T_{3}^{2}+T_{4}^{2}}$
&$\arraycolsep=1mm\begin{bmatrix}{1} & {2} & {6} & {6} & {3} \\ {\bar{0}} & {\bar{0}} & {\bar{1}} & {\bar{0}} & {\bar{1}}\end{bmatrix}$&${6}$
&${309}$&${T_{1}^{5}T_{2}^{3}+T_{3}^{2}+T_{4}^{2}}$
&$\arraycolsep=1mm\begin{bmatrix}{1} & {1} & {4} & {4} & {2} \\ {\bar{0}} & {\bar{0}} & {\bar{1}} & {\bar{0}} & {\bar{1}}\end{bmatrix}$&${8}$
\\
\hline
${284}$&${T_{1}^{6}T_{2}^{3}+T_{3}^{2}+T_{4}^{2}}$
&$\arraycolsep=1mm\begin{bmatrix}{1} & {2} & {6} & {6} & {3} \\ {\bar{0}} & {\bar{0}} & {\bar{1}} & {\bar{0}} & {\bar{1}}\end{bmatrix}$&${6}$
&${310}$&${T_{1}^{3}T_{2}^{2}+T_{3}^{2}+T_{4}^{2}}$
&$\arraycolsep=1mm\begin{bmatrix}{2} & {1} & {4} & {4} & {1} \\ {\bar{0}} & {\bar{1}} & {\bar{1}} & {\bar{0}} & {\bar{0}}\end{bmatrix}$&${8}$
\\
\hline
\end{tabu}\vspace{-4mm}
\end{table}

\begin{table}[h!]\centering\setlength{\tabcolsep}{3.5pt}
\renewcommand{\arraystretch}{1.3}\small\tabulinesep=1mm
\begin{tabu}{c|c|c|c||c|c|c|c}
\hline
ID&relations&gd-matrix&$-\mathcal{K}^3$&ID&relations&gd-matrix&$-\mathcal{K}^3$
\\
\hline
${336}$&${T_{1}T_{2}+T_{3}^{6}+T_{4}^{4}}$
&$\arraycolsep=1mm\begin{bmatrix}{2} & {10} & {2} & {3} & {5} \\ {\bar{0}} & {\bar{0}} & {\bar{1}} & {\bar{0}} & {\bar{1}}\end{bmatrix}$&${10}$
&${374}$&${T_{1}T_{2}+T_{3}^{12}+T_{4}^{2}}$
&$\arraycolsep=1mm\begin{bmatrix}{4} & {8} & {1} & {6} & {1} \\ {\bar{0}} & {\bar{0}} & {\bar{1}} & {\bar{1}} & {\bar{0}}\end{bmatrix}$&${16}$
\\
\hline
${340}$&${T_{1}^{3}T_{2}+T_{3}^{3}+T_{4}^{2}}$
&$\arraycolsep=1mm\begin{bmatrix}{1} & {3} & {2} & {3} & {3} \\ {\bar{0}} & {\bar{0}} & {\bar{0}} & {\bar{1}} & {\bar{1}}\end{bmatrix}$&${12}$
&${375}$&${T_{1}T_{2}+T_{3}^{10}+T_{4}^{2}}$
&$\arraycolsep=1mm\begin{bmatrix}{2} & {8} & {1} & {5} & {2} \\ {\bar{0}} & {\bar{0}} & {\bar{0}} & {\bar{1}} & {\bar{1}}\end{bmatrix}$&${16}$
\\
\hline
${341}$&${T_{1}^{4}T_{2}+T_{3}^{2}+T_{4}^{2}}$
&$\arraycolsep=1mm\begin{bmatrix}{1} & {2} & {3} & {3} & {3} \\ {\bar{0}} & {\bar{0}} & {\bar{0}} & {\bar{1}} & {\bar{1}}\end{bmatrix}$&${12}$
&${376}$&${T_{1}T_{2}+T_{3}^{6}+T_{4}^{2}}$
&$\arraycolsep=1mm\begin{bmatrix}{2} & {4} & {1} & {3} & {4} \\ {\bar{1}} & {\bar{1}} & {\bar{0}} & {\bar{1}} & {\bar{1}}\end{bmatrix}$&${16}$
\\
\hline
${342}$&${T_{1}T_{2}+T_{3}^{18}+T_{4}^{2}}$
&$\arraycolsep=1mm\begin{bmatrix}{6} & {12} & {1} & {9} & {2} \\ {\bar{0}} & {\bar{0}} & {\bar{0}} & {\bar{1}} & {\bar{1}}\end{bmatrix}$&${12}$
&${377}$&${T_{1}T_{2}+T_{3}^{6}+T_{4}^{2}}$
&$\arraycolsep=1mm\begin{bmatrix}{2} & {4} & {1} & {3} & {4} \\ {\bar{0}} & {\bar{0}} & {\bar{0}} & {\bar{1}} & {\bar{1}}\end{bmatrix}$&${16}$
\\
\hline
${343}$&${T_{1}T_{2}+T_{3}^{16}+T_{4}^{2}}$
&$\arraycolsep=1mm\begin{bmatrix}{4} & {12} & {1} & {8} & {3} \\ {\bar{0}} & {\bar{0}} & {\bar{0}} & {\bar{1}} & {\bar{1}}\end{bmatrix}$&${12}$
&${378}$&${T_{1}T_{2}+T_{3}^{4}+T_{4}^{2}}$
&$\arraycolsep=1mm\begin{bmatrix}{2} & {2} & {1} & {2} & {1} \\ {\bar{1}} & {\bar{1}} & {\bar{1}} & {\bar{1}} & {\bar{0}}\end{bmatrix}$&${16}$
\\
\hline
${344}$&${T_{1}T_{2}+T_{3}^{10}+T_{4}^{2}}$
&$\arraycolsep=1mm\begin{bmatrix}{4} & {6} & {1} & {5} & {6} \\ {\bar{0}} & {\bar{0}} & {\bar{0}} & {\bar{1}} & {\bar{1}}\end{bmatrix}$&${12}$
&${379}$&${T_{1}T_{2}+T_{3}^{4}+T_{4}^{2}}$
&$\arraycolsep=1mm\begin{bmatrix}{2} & {2} & {1} & {2} & {1} \\ {\bar{0}} & {\bar{0}} & {\bar{1}} & {\bar{1}} & {\bar{0}}\end{bmatrix}$&${16}$
\\
\hline
${345}$&${T_{1}T_{2}+T_{3}^{3}+T_{4}^{2}}$
&$\arraycolsep=1mm\begin{bmatrix}{3} & {3} & {2} & {3} & {1} \\ {\bar{1}} & {\bar{1}} & {\bar{0}} & {\bar{0}} & {\bar{0}}\end{bmatrix}$&${12}$
&${393}$&${T_{1}^{2}T_{2}+T_{3}^{8}+T_{4}^{2}}$
&$\arraycolsep=1mm\begin{bmatrix}{1} & {6} & {1} & {4} & {2} \\ {\bar{0}} & {\bar{0}} & {\bar{1}} & {\bar{1}} & {\bar{0}}\end{bmatrix}$&${18}$
\\
\hline
${346}$&${T_{1}T_{2}+T_{3}^{2}+T_{4}^{2}}$
&$\arraycolsep=1mm\begin{bmatrix}{2} & {4} & {3} & {3} & {6} \\ {\bar{0}} & {\bar{0}} & {\bar{0}} & {\bar{1}} & {\bar{1}}\end{bmatrix}$&${12}$
&${394}$&${T_{1}^{2}T_{2}+T_{3}^{4}+T_{4}^{2}}$
&$\arraycolsep=1mm\begin{bmatrix}{1} & {6} & {2} & {4} & {1} \\ {\bar{1}} & {\bar{0}} & {\bar{0}} & {\bar{1}} & {\bar{0}}\end{bmatrix}$&${18}$
\\
\hline
${365}$&${T_{1}^{4}T_{2}+T_{3}^{3}+T_{4}^{2}}$
&$\arraycolsep=1mm\begin{bmatrix}{1} & {8} & {4} & {6} & {1} \\ {\bar{1}} & {\bar{0}} & {\bar{0}} & {\bar{1}} & {\bar{0}}\end{bmatrix}$&${16}$
&${395}$&${T_{1}^{2}T_{2}+T_{3}^{4}+T_{4}^{2}}$
&$\arraycolsep=1mm\begin{bmatrix}{1} & {6} & {2} & {4} & {1} \\ {\bar{0}} & {\bar{0}} & {\bar{1}} & {\bar{1}} & {\bar{0}}\end{bmatrix}$&${18}$
\\
\hline
${366}$&${T_{1}^{2}T_{2}+T_{3}^{6}+T_{4}^{2}}$
&$\arraycolsep=1mm\begin{bmatrix}{1} & {4} & {1} & {3} & {1} \\ {\bar{1}} & {\bar{0}} & {\bar{0}} & {\bar{1}} & {\bar{0}}\end{bmatrix}$&${16}$
&${396}$&${T_{1}^{2}T_{2}+T_{3}^{2}+T_{4}^{2}}$
&$\arraycolsep=1mm\begin{bmatrix}{1} & {2} & {2} & {2} & {3} \\ {\bar{1}} & {\bar{0}} & {\bar{1}} & {\bar{0}} & {\bar{0}}\end{bmatrix}$&${18}$
\\
\hline
${367}$&${T_{1}^{2}T_{2}+T_{3}^{6}+T_{4}^{2}}$
&$\arraycolsep=1mm\begin{bmatrix}{1} & {4} & {1} & {3} & {1} \\ {\bar{1}} & {\bar{0}} & {\bar{1}} & {\bar{0}} & {\bar{0}}\end{bmatrix}$&${16}$
&${397}$&${T_{1}T_{2}+T_{3}^{8}+T_{4}^{2}}$
&$\arraycolsep=1mm\begin{bmatrix}{2} & {6} & {1} & {4} & {1} \\ {\bar{0}} & {\bar{0}} & {\bar{1}} & {\bar{1}} & {\bar{0}}\end{bmatrix}$&${18}$
\\
\hline
${368}$&${T_{1}^{2}T_{2}+T_{3}^{5}+T_{4}^{2}}$
&$\arraycolsep=1mm\begin{bmatrix}{1} & {8} & {2} & {5} & {2} \\ {\bar{0}} & {\bar{0}} & {\bar{0}} & {\bar{1}} & {\bar{1}}\end{bmatrix}$&${16}$
&${398}$&${T_{1}T_{2}+T_{3}^{4}+T_{4}^{2}}$
&$\arraycolsep=1mm\begin{bmatrix}{2} & {2} & {1} & {2} & {3} \\ {\bar{0}} & {\bar{0}} & {\bar{0}} & {\bar{1}} & {\bar{1}}\end{bmatrix}$&${18}$
\\
\hline
${369}$&${T_{1}^{2}T_{2}+T_{3}^{4}+T_{4}^{2}}$
&$\arraycolsep=1mm\begin{bmatrix}{1} & {2} & {1} & {2} & {2} \\ {\bar{0}} & {\bar{0}} & {\bar{1}} & {\bar{1}} & {\bar{0}}\end{bmatrix}$&${16}$
&${404}$&${T_{1}T_{2}+T_{3}^{3}+T_{4}^{2}}$
&$\arraycolsep=1mm\begin{bmatrix}{1} & {5} & {2} & {3} & {5} \\ {\bar{1}} & {\bar{1}} & {\bar{0}} & {\bar{0}} & {\bar{0}}\end{bmatrix}$&${20}$
\\
\hline
${370}$&${T_{1}^{2}T_{2}+T_{3}^{4}+T_{4}^{2}}$
&$\arraycolsep=1mm\begin{bmatrix}{1} & {2} & {1} & {2} & {2} \\ {\bar{0}} & {\bar{0}} & {\bar{0}} & {\bar{1}} & {\bar{1}}\end{bmatrix}$&${16}$
&${406}$&${T_{1}T_{2}+T_{3}^{4}+T_{4}^{2}}$
&$\arraycolsep=1mm\begin{bmatrix}{1} & {3} & {1} & {2} & {3} \\ {\bar{1}} & {\bar{1}} & {\bar{0}} & {\bar{1}} & {\bar{1}}\end{bmatrix}$&${24}$
\\
\hline
${371}$&${T_{1}^{2}T_{2}+T_{3}^{3}+T_{4}^{2}}$
&$\arraycolsep=1mm\begin{bmatrix}{1} & {4} & {2} & {3} & {4} \\ {\bar{1}} & {\bar{0}} & {\bar{0}} & {\bar{0}} & {\bar{1}}\end{bmatrix}$&${16}$
&${407}$&${T_{1}T_{2}+T_{3}^{4}+T_{4}^{2}}$
&$\arraycolsep=1mm\begin{bmatrix}{1} & {3} & {1} & {2} & {3} \\ {\bar{0}} & {\bar{0}} & {\bar{0}} & {\bar{1}} & {\bar{1}}\end{bmatrix}$&${24}$
\\
\hline
${372}$&${T_{1}^{3}T_{2}+T_{3}^{2}+T_{4}^{2}}$
&$\arraycolsep=1mm\begin{bmatrix}{1} & {1} & {2} & {2} & {2} \\ {\bar{0}} & {\bar{0}} & {\bar{0}} & {\bar{1}} & {\bar{1}}\end{bmatrix}$&${16}$
&${412}$&${T_{1}T_{2}+T_{3}^{2}+T_{4}^{2}}$
&$\arraycolsep=1mm\begin{bmatrix}{1} & {1} & {1} & {1} & {2} \\ {\bar{0}} & {\bar{0}} & {\bar{0}} & {\bar{1}} & {\bar{1}}\end{bmatrix}$&${32}$
\\
\hline
${373}$&${T_{1}^{2}T_{2}+T_{3}^{2}+T_{4}^{2}}$
&$\arraycolsep=1mm\begin{bmatrix}{1} & {2} & {2} & {2} & {1} \\ {\bar{1}} & {\bar{0}} & {\bar{0}} & {\bar{1}} & {\bar{0}}\end{bmatrix}$&${16}$
\\
\cline{1-4}
\end{tabu}\vspace{1mm}
\end{table}
\end{class-list}

\begin{class-list}
Non-toric, $\QQ$-factorial, Gorenstein, log terminal Fano
threefolds of Picard number one with an effective two-torus action:
Specifying data for divisor class group $\ZZ\oplus \ZZ/2\ZZ$ and
format $(2,1,1,1,1)$.
\begin{table}[h!]\centering\setlength{\tabcolsep}{5.5pt}
\renewcommand{\arraystretch}{1.25}\small\tabulinesep=1mm
\begin{tabu}{c|c|c|c}
\hline
ID&relations&gd-matrix&$-\mathcal{K}^3$
\\
\hline
${255}$&$\setlength{\arraycolsep}{1pt}\begin{array}{r}{T_{1}T_{2}+T_{3}^{3}+T_{4}^{2}},\\ {\lambda T_{3}^{3}+T_{4}^{2}+T_{5}^{2}}\end{array}$
&$\arraycolsep=1mm\begin{bmatrix}{4} & {2} & {2} & {3} & {3} & {2} \\ {\bar{0}} & {\bar{0}} & {\bar{0}} & {\bar{1}} & {\bar{0}} & {\bar{1}}\end{bmatrix}$&${4}$
\\
\hline
\end{tabu}\vspace{-2mm}
\end{table}
\end{class-list}

\begin{class-list}
Non-toric, $\QQ$-factorial, Gorenstein, log terminal Fano
threefolds of Picard number one with an effective two-torus action:
Specifying data for divisor class group $\ZZ\oplus \ZZ/2\ZZ$ and format $(1,1,1,2)$.
\begin{table}[h!]\centering\setlength{\tabcolsep}{3.5pt}
\renewcommand{\arraystretch}{1.15}\small\tabulinesep=1mm
\begin{tabu}{c|c|c|c||c|c|c|c}
\hline
ID&relations&gd-matrix&$-\mathcal{K}^3$&ID&relations&gd-matrix&$-\mathcal{K}^3$
\\
\hline
${338}$&${T_{1}^{3}+T_{2}^{2}+T_{3}^{2}}$
&$\arraycolsep=1mm\begin{bmatrix}{2} & {3} & {3} & {6} & {4} \\ {\bar{0}} & {\bar{0}} & {\bar{1}} & {\bar{1}} & {\bar{0}}\end{bmatrix}$&${12}$
&${339}$&${T_{1}^{3}+T_{2}^{2}+T_{3}^{2}}$
&$\arraycolsep=1mm\begin{bmatrix}{2} & {3} & {3} & {3} & {1} \\ {\bar{0}} & {\bar{1}} & {\bar{0}} & {\bar{1}} & {\bar{0}}\end{bmatrix}$&${12}$
\\
\hline
\end{tabu}\vspace{-3mm}
\end{table}
\end{class-list}

\begin{class-list}\label{class-list(0,2,2)}
Non-toric, $\QQ$-factorial, Gorenstein, log terminal Fano
threefolds of Picard number one with an effective two-torus
action: Specifying data for divisor class group
$\ZZ\oplus \ZZ/2\ZZ\oplus \ZZ/2\ZZ$ and format $(2,2,1,0)$.
\begin{table}[h!]\centering\setlength{\tabcolsep}{3.5pt}
\renewcommand{\arraystretch}{1.15}\small\tabulinesep=1mm
\begin{tabu}{c|c|c|c||c|c|c|c}
\hline
ID&relations&gd-matrix&$-\mathcal{K}^3$&ID&relations&gd-matrix&$-\mathcal{K}^3$
\\
\hline
${480}$&${T_{1}^{4}T_{2}^{2}+T_{3}^{2}T_{4}^{2}+T_{5}^{2}}$
&$\arraycolsep=1mm\begin{bmatrix}{1} &{1} &{2} &{1} &{3}\\
{\bar{0}} &{\bar{1}} &{\bar{0}} &{\bar{1}} &{\bar{0}}\\
{\bar{0}} &{\bar{0}} &{\bar{0}} &{\bar{1}} &{\bar{1}}\end{bmatrix}$&${2}$
&${495}$&${T_{1}^{2}T_{2}+T_{3}^{2}T_{4}^{2}+T_{5}^{2}}$
&$\arraycolsep=1mm\begin{bmatrix}{2} & {2} & {2} & {1} & {3} \\ {\bar{0}} & {\bar{0}} & {\bar{1}} & {\bar{0}} & {\bar{1}} \\ {\bar{1}} & {\bar{0}} & {\bar{0}} & {\bar{0}} & {\bar{1}}\end{bmatrix}$&${4}$
\\
\hline
${491}$&${T_{1}^{6}T_{2}^{2}+T_{3}^{2}T_{4}+T_{5}^{2}}$
&$\arraycolsep=1mm\begin{bmatrix}{1} &{1} &{2} &{4} &{4}\\
{\bar{0}} &{\bar{1}} &{\bar{0}} &{\bar{0}} &{\bar{1}}\\
{\bar{0}} &{\bar{1}} &{\bar{1}} &{\bar{0}} &{\bar{0}}\end{bmatrix}$&${4}$
&${501}$&${T_{1}^{3}T_{2}+T_{3}T_{4}+T_{5}^{3}}$
&$\arraycolsep=1mm\begin{bmatrix}{1} & {3} & {3} & {3} & {2} \\ {\bar{0}} & {\bar{1}} & {\bar{0}} & {\bar{1}} & {\bar{1}} \\ {\bar{0}} & {\bar{1}} & {\bar{1}} & {\bar{0}} & {\bar{1}}\end{bmatrix}$
&${6}$
\\
\hline
${492}$&${T_{1}^{4}T_{2}^{4}+T_{3}^{2}T_{4}+T_{5}^{2}}$
&$\arraycolsep=1mm\begin{bmatrix}{1} &{1} &{2} &{4} &{4}\\
{\bar{1}} &{\bar{0}} &{\bar{1}} &{\bar{0}} &{\bar{0}}\\
{\bar{1}} &{\bar{0}} &{\bar{0}} &{\bar{0}} &{\bar{1}}\end{bmatrix}$&${4}$
&${518}$&${T_{1}^{4}T_{2}^{2}+T_{3}^{2}T_{4}+T_{5}^{2}}$
&$\arraycolsep=1mm\begin{bmatrix}{1} & {1} & {1} & {4} & {3} \\ {\bar{0}} & {\bar{1}} & {\bar{0}} & {\bar{0}} & {\bar{1}} \\ {\bar{0}} & {\bar{1}} & {\bar{1}} & {\bar{0}} & {\bar{0}}\end{bmatrix}$
&${8}$
\\
\hline
${493}$&${T_{1}^{2}T_{2}^{2}+T_{3}^{2}T_{4}^{2}+T_{5}^{2}}$
&$\arraycolsep=1mm\begin{bmatrix}{1} &{1} &{1} &{1} &{2}\\
{\bar{1}} &{\bar{0}} &{\bar{1}} &{\bar{0}} &{\bar{0}}\\
{\bar{0}} &{\bar{0}} &{\bar{0}} &{\bar{1}} &{\bar{1}}\end{bmatrix}$&${4}$
&${520}$&${T_{1}T_{2}+T_{3}T_{4}+T_{5}^{3}}$
&$\arraycolsep=1mm\begin{bmatrix}{1} & {2} & {1} & {2} & {1} \\ {\bar{1}} & {\bar{1}} & {\bar{1}} & {\bar{1}} & {\bar{0}} \\ {\bar{1}} & {\bar{1}} & {\bar{0}} & {\bar{0}} & {\bar{0}}\end{bmatrix}$
&${12}$
\\
\hline
${494}$&${T_{1}^{2}T_{2}^{2}+T_{3}^{2}T_{4}+T_{5}^{4}}$
&$\arraycolsep=1mm\begin{bmatrix}{1} &{1} &{1} &{2} &{1}\\
{\bar{1}} &{\bar{1}} &{\bar{0}} &{\bar{0}} &{\bar{0}}\\
{\bar{1}} &{\bar{0}} &{\bar{1}} &{\bar{0}} &{\bar{0}}\end{bmatrix}$&${4}$
\\
\cline{1-4}
\end{tabu}\vspace{-3mm}
\end{table}
\end{class-list}

\begin{class-list}
Non-toric, $\QQ$-factorial, Gorenstein, log terminal Fano
threefolds of Picard number one with an effective two-torus action:
Specifying data for divisor class group
$\ZZ\oplus \ZZ/2\ZZ\oplus \ZZ/2\ZZ$ and format $(2,2,1,1,0)$.
\begin{table}[h!]\centering\setlength{\tabcolsep}{5.5pt}
\renewcommand{\arraystretch}{1.15}\small\tabulinesep=1mm
\begin{tabu}{c|c|c|c}
\hline
ID&relations&gd-matrix&$-\mathcal{K}^3$
\\
\hline
${479}$&$\begin{array}{r}{T_{1}^{2}T_{2}^{2}+T_{3}T_{4}+T_{5}^{2}},\\ {\lambda T_{3}T_{4}+T_{5}^{2}+T_{6}^{2}}\end{array}$
&
$\arraycolsep=1mm\begin{bmatrix}{1} & {1} & {2} & {2} & {2} & {2} \\ {\bar{0}} & {\bar{1}} & {\bar{0}} & {\bar{0}} & {\bar{1}} & {\bar{0}} \\ {\bar{0}} & {\bar{0}} & {\bar{0}} & {\bar{0}} & {\bar{1}} & {\bar{1}}\end{bmatrix}$&${2}$
\\
\hline
\end{tabu}\vspace{-4mm}
\end{table}
\end{class-list}

\begin{class-list}
Non-toric, $\QQ$-factorial, Gorenstein, log terminal Fano
threefolds of Picard number one with an effective two-torus action:
Specifying data for divisor class group
$\ZZ\oplus \ZZ/2\ZZ\oplus \ZZ/2\ZZ$ and format $(2,2,1,1,1,0)$.
\begin{table}[h!]\centering\setlength{\tabcolsep}{5.5pt}
\renewcommand{\arraystretch}{1.25}\small\tabulinesep=1mm
\begin{tabu}{c|c|c|c}
\hline
ID&relations&gd-matrix&$-\mathcal{K}^3$
\\
\hline
${478}$&
$\begin{array}{r}{T_{1}T_{2}+T_{3}T_{4}+T_{5}^{2}},\\ {\lambda_{1} T_{3}T_{4}+T_{5}^{2}+T_{6}^{2}},\\ {\lambda_{2} T_{5}^{2}+T_{6}^{2}+T_{7}^{2}}\end{array}$
&$\begin{bmatrix}{1} & {1} & {1} & {1} & {1} & {1} & {1} \\ {\bar{0}} & {\bar{0}} & {\bar{0}} & {\bar{0}} & {\bar{1}} & {\bar{1}} & {\bar{0}} \\ {\bar{1}} & {\bar{1}} & {\bar{1}} & {\bar{1}} & {\bar{1}} & {\bar{0}} & {\bar{0}}\end{bmatrix}$&${2}$
\\
\hline
\end{tabu}
\end{table}
\end{class-list}

\begin{class-list}
Non-toric, $\QQ$-factorial, Gorenstein, log terminal Fano
threefolds of Picard number one with an effective two-torus action:
Specifying data for divisor class group $\ZZ\oplus \ZZ/2\ZZ\oplus \ZZ/2\ZZ$
and format $(3,1,1,0)$.
\begin{table}[h!]\centering\setlength{\tabcolsep}{3.5pt}
\renewcommand{\arraystretch}{1.15}\small\tabulinesep=1mm
\begin{tabu}{c|c|c|c||c|c|c|c}
\hline
ID&relations&gd-matrix&$-\mathcal{K}^3$&ID&relations&gd-matrix&$-\mathcal{K}^3$
\\
\hline
${477}$&${T_{1}^{2}T_{2}^{2}T_{3}+T_{4}^{6}+T_{5}^{2}}$
&$\arraycolsep=1mm\begin{bmatrix}{1} & {1} & {2} & {1} & {3} \\ {\bar{0}} & {\bar{0}} & {\bar{0}} & {\bar{1}} & {\bar{1}} \\ {\bar{1}} & {\bar{0}} & {\bar{0}} & {\bar{0}} & {\bar{1}}\end{bmatrix}$&${2}$
&${490}$&${T_{1}^{2}T_{2}T_{3}+T_{4}^{4}+T_{5}^{2}}$
&$\arraycolsep=1mm\begin{bmatrix}{1} & {1} & {1} & {1} & {2} \\ {\bar{0}} & {\bar{1}} & {\bar{1}} & {\bar{0}} & {\bar{0}} \\ {\bar{0}} & {\bar{1}} & {\bar{1}} & {\bar{1}} & {\bar{1}}\end{bmatrix}$
&${4}$
\\
\hline
${489}$&${T_{1}^{2}T_{2}^{2}T_{3}+T_{4}^{4}+T_{5}^{2}}$
&$\arraycolsep=1mm\begin{bmatrix}{1} & {1} & {4} & {2} & {4} \\ {\bar{1}} & {\bar{0}} & {\bar{0}} & {\bar{1}} & {\bar{0}} \\ {\bar{0}} & {\bar{0}} & {\bar{0}} & {\bar{1}} & {\bar{1}}\end{bmatrix}$
&${4}$
\\
\cline{1-4}
\end{tabu}
\end{table}
\end{class-list}

\begin{class-list}
Non-toric, $\QQ$-factorial, Gorenstein, log terminal Fano
threefolds of Picard number one with an effective two-torus action:
Specifying data for divisor class group
$\ZZ\oplus \ZZ/2\ZZ\oplus \ZZ/2\ZZ$ and format $(3,1,1,1,0)$.
\begin{table}[h!]\centering\setlength{\tabcolsep}{5.5pt}
\renewcommand{\arraystretch}{1.25}\small\tabulinesep=1mm
\begin{tabu}{c|c|c|c}
\hline
ID&relations&gd-matrix&$-\mathcal{K}^3$
\\
\hline
${476}$&$\begin{array}{r}{T_{1}T_{2}T_{3}+T_{4}^{2}+T_{5}^{2}},\\ {\lambda T_{4}^{2}+T_{5}^{2}+T_{6}^{2}}\end{array}$
&$\arraycolsep=1mm\begin{bmatrix}{1} & {1} & {2} & {2} & {2} & {2} \\ {\bar{0}} & {\bar{0}} & {\bar{0}} & {\bar{1}} & {\bar{1}} & {\bar{0}} \\ {\bar{0}} & {\bar{0}} & {\bar{0}} & {\bar{0}} & {\bar{1}} & {\bar{1}}\end{bmatrix}$&${2}$
\\
\hline
\end{tabu}
\end{table}
\end{class-list}

\begin{class-list}
Non-toric, $\QQ$-factorial, Gorenstein, log terminal Fano
threefolds of Picard number one with an effective two-torus action:
Specifying data for divisor class group
$\ZZ\oplus \ZZ/2\ZZ\oplus \ZZ/2\ZZ$ and format $(2,1,1,1)$.
\begin{table}[h!]\centering\setlength{\tabcolsep}{3.5pt}
\renewcommand{\arraystretch}{1.15}\small\tabulinesep=1mm
\begin{tabu}{c|c|c|c||c|c|c|c}
\hline
ID&relations&gd-matrix&$-\mathcal{K}^3$&ID&relations&gd-matrix&$-\mathcal{K}^3$
\\
\hline
${474}$&${T_{1}^{4}T_{2}^{2}+T_{3}^{3}+T_{4}^{2}}$
&$\arraycolsep=1mm\begin{bmatrix}{1} & {1} & {2} & {3} & {1} \\ {\bar{0}} & {\bar{1}} & {\bar{0}} & {\bar{1}} & {\bar{0}} \\ {\bar{1}} & {\bar{0}} & {\bar{0}} & {\bar{1}} & {\bar{0}}\end{bmatrix}$
&${2}$
&${481}$&${T_{1}^{3}T_{2}+T_{3}^{4}+T_{4}^{2}}$
&$\arraycolsep=1mm\begin{bmatrix}{1} & {1} & {1} & {2} & {1} \\ {\bar{1}} & {\bar{1}} & {\bar{0}} & {\bar{0}} & {\bar{0}} \\ {\bar{0}} & {\bar{0}} & {\bar{1}} & {\bar{1}} & {\bar{0}}\end{bmatrix}$
&${4}$
\\
\hline
${475}$&${T_{1}^{2}T_{2}^{2}+T_{3}^{6}+T_{4}^{2}}$
&$\arraycolsep=1mm\begin{bmatrix}{2} & {1} & {1} & {3} & {1} \\ {\bar{0}} & {\bar{1}} & {\bar{1}} & {\bar{0}} & {\bar{0}} \\ {\bar{0}} & {\bar{0}} & {\bar{1}} & {\bar{1}} & {\bar{0}}\end{bmatrix}$
&${2}$
&${482}$&${T_{1}^{2}T_{2}^{2}+T_{3}^{4}+T_{4}^{2}}$
&$\arraycolsep=1mm\begin{bmatrix}{1} & {1} & {1} & {2} & {1} \\ {\bar{0}} & {\bar{0}} & {\bar{1}} & {\bar{1}} & {\bar{0}} \\ {\bar{0}} & {\bar{1}} & {\bar{1}} & {\bar{0}} & {\bar{0}}\end{bmatrix}$
&${4}$
\\
\hline
\end{tabu}\vspace{2mm}
\end{table}

\begin{table}[t!]\centering\setlength{\tabcolsep}{3.5pt}
\renewcommand{\arraystretch}{1.05}\small\tabulinesep=1mm
\begin{tabu}{c|c|c|c||c|c|c|c}
\hline
ID&relations&gd-matrix&$-\mathcal{K}^3$&ID&relations&gd-matrix&$-\mathcal{K}^3$
\\
\hline
${483}$&${T_{1}^{2}T_{2}^{2}+T_{3}^{4}+T_{4}^{2}}$
&$\arraycolsep=1mm\begin{bmatrix}{1} & {1} & {1} & {2} & {1} \\ {\bar{1}} & {\bar{1}} & {\bar{1}} & {\bar{1}} & {\bar{0}} \\ {\bar{0}} & {\bar{1}} & {\bar{1}} & {\bar{0}} & {\bar{0}}\end{bmatrix}$
&${4}$
&${511}$&${T_{1}^{3}T_{2}+T_{3}^{2}+T_{4}^{2}}$
&$\arraycolsep=1mm\begin{bmatrix}{1} & {1} & {2} & {2} & {2} \\ {\bar{0}} & {\bar{0}} & {\bar{1}} & {\bar{1}} & {\bar{0}} \\ {\bar{0}} & {\bar{0}} & {\bar{1}} & {\bar{0}} & {\bar{1}}\end{bmatrix}$
&${8}$
\\
\hline
${484}$&${T_{1}^{2}T_{2}^{2}+T_{3}^{3}+T_{4}^{2}}$
&$\arraycolsep=1mm\begin{bmatrix}{2} & {1} & {2} & {3} & {2} \\ {\bar{1}} & {\bar{0}} & {\bar{0}} & {\bar{1}} & {\bar{0}} \\ {\bar{1}} & {\bar{0}} & {\bar{0}} & {\bar{0}} & {\bar{1}}\end{bmatrix}$
&${4}$
&${512}$&${T_{1}^{2}T_{2}^{2}+T_{3}^{2}+T_{4}^{2}}$
&$\arraycolsep=1mm\begin{bmatrix}{1} & {1} & {2} & {2} & {2} \\ {\bar{1}} & {\bar{0}} & {\bar{1}} & {\bar{0}} & {\bar{0}} \\ {\bar{0}} & {\bar{0}} & {\bar{1}} & {\bar{1}} & {\bar{0}}\end{bmatrix}$
&${8}$
\\
\hline
${485}$&${T_{1}^{2}T_{2}+T_{3}^{8}+T_{4}^{2}}$
&$\arraycolsep=1mm\begin{bmatrix}{2} & {4} & {1} & {4} & {1} \\ {\bar{1}} & {\bar{0}} & {\bar{1}} & {\bar{0}} & {\bar{0}} \\ {\bar{1}} & {\bar{0}} & {\bar{0}} & {\bar{1}} & {\bar{0}}\end{bmatrix}$
&${4}$
&${513}$&${T_{1}^{2}T_{2}^{2}+T_{3}^{2}+T_{4}^{2}}$
&$\arraycolsep=1mm\begin{bmatrix}{1} & {1} & {2} & {2} & {2} \\ {\bar{0}} & {\bar{0}} & {\bar{0}} & {\bar{1}} & {\bar{1}} \\ {\bar{0}} & {\bar{0}} & {\bar{1}} & {\bar{1}} & {\bar{0}}\end{bmatrix}$
&${8}$
\\
\hline
${486}$&${T_{1}^{2}T_{2}+T_{3}^{6}+T_{4}^{2}}$&
$\arraycolsep=1mm\begin{bmatrix}{2} & {2} & {1} & {3} & {2} \\ {\bar{1}} & {\bar{0}} & {\bar{0}} & {\bar{0}} & {\bar{1}} \\ {\bar{1}} & {\bar{0}} & {\bar{0}} & {\bar{1}} & {\bar{0}}\end{bmatrix}$
&${4}$
&${514}$&${T_{1}^{2}T_{2}^{2}+T_{3}^{2}+T_{4}^{2}}$
&$\arraycolsep=1mm\begin{bmatrix}{1} & {1} & {2} & {2} & {2} \\ {\bar{0}} & {\bar{0}} & {\bar{0}} & {\bar{1}} & {\bar{1}} \\ {\bar{1}} & {\bar{0}} & {\bar{0}} & {\bar{0}} & {\bar{1}}\end{bmatrix}$
&${8}$
\\
\hline
${487}$&${T_{1}^{6}T_{2}^{2}+T_{3}^{2}+T_{4}^{2}}$
&$\arraycolsep=1mm\begin{bmatrix}{1} & {1} & {4} & {4} & {2} \\ {\bar{0}} & {\bar{1}} & {\bar{1}} & {\bar{0}} & {\bar{0}} \\ {\bar{0}} & {\bar{1}} & {\bar{0}} & {\bar{0}} & {\bar{1}}\end{bmatrix}$
&${4}$
&${515}$&${T_{1}^{2}T_{2}+T_{3}^{2}+T_{4}^{2}}$
&$\arraycolsep=1mm\begin{bmatrix}{1} & {2} & {2} & {2} & {1} \\ {\bar{0}} & {\bar{0}} & {\bar{1}} & {\bar{1}} & {\bar{0}} \\ {\bar{1}} & {\bar{0}} & {\bar{1}} & {\bar{0}} & {\bar{0}}\end{bmatrix}$
&${8}$
\\
\hline
${488}$&${T_{1}^{4}T_{2}^{2}+T_{3}^{2}+T_{4}^{2}}$
&$\arraycolsep=1mm\begin{bmatrix}{1} & {2} & {4} & {4} & {1} \\ {\bar{0}} & {\bar{1}} & {\bar{1}} & {\bar{0}} & {\bar{0}} \\ {\bar{1}} & {\bar{1}} & {\bar{0}} & {\bar{0}} & {\bar{0}}\end{bmatrix}$
&${4}$
&${516}$&${T_{1}T_{2}+T_{3}^{6}+T_{4}^{2}}$
&$\arraycolsep=1mm\begin{bmatrix}{2} & {4} & {1} & {3} & {4} \\ {\bar{1}} & {\bar{1}} & {\bar{0}} & {\bar{0}} & {\bar{0}} \\ {\bar{1}} & {\bar{1}} & {\bar{0}} & {\bar{1}} & {\bar{1}}\end{bmatrix}$
&${8}$
\\
\hline
${499}$&${T_{1}^{4}T_{2}+T_{3}^{2}+T_{4}^{2}}$
&$\arraycolsep=1mm\begin{bmatrix}{1} & {2} & {3} & {3} & {3} \\ {\bar{0}} & {\bar{0}} & {\bar{1}} & {\bar{1}} & {\bar{0}} \\ {\bar{0}} & {\bar{0}} & {\bar{1}} & {\bar{0}} & {\bar{1}}\end{bmatrix}$
&${6}$
&${517}$&${T_{1}T_{2}+T_{3}^{4}+T_{4}^{2}}$
&$\arraycolsep=1mm\begin{bmatrix}{2} & {2} & {1} & {2} & {1} \\ {\bar{1}} & {\bar{1}} & {\bar{0}} & {\bar{0}} & {\bar{0}} \\ {\bar{1}} & {\bar{1}} & {\bar{1}} & {\bar{1}} & {\bar{0}}\end{bmatrix}$
&${8}$
\\
\hline
${500}$&${T_{1}^{2}T_{2}^{2}+T_{3}^{2}+T_{4}^{2}}$
&$\arraycolsep=1mm\begin{bmatrix}{2} & {1} & {3} & {3} & {3} \\ {\bar{0}} & {\bar{0}} & {\bar{0}} & {\bar{1}} & {\bar{1}} \\ {\bar{0}} & {\bar{0}} & {\bar{1}} & {\bar{1}} & {\bar{0}}\end{bmatrix}$
&${6}$
&${519}$&${T_{1}T_{2}+T_{3}^{4}+T_{4}^{2}}$
&$\arraycolsep=1mm\begin{bmatrix}{1} & {3} & {1} & {2} & {3} \\ {\bar{1}} & {\bar{1}} & {\bar{0}} & {\bar{0}} & {\bar{0}} \\ {\bar{1}} & {\bar{1}} & {\bar{0}} & {\bar{1}} & {\bar{1}}\end{bmatrix}$
&${12}$
\\
\hline
${509}$&${T_{1}^{2}T_{2}+T_{3}^{6}+T_{4}^{2}}$
&$\arraycolsep=1mm\begin{bmatrix}{1} & {4} & {1} & {3} & {1} \\ {\bar{0}} & {\bar{0}} & {\bar{1}} & {\bar{1}} & {\bar{0}} \\ {\bar{1}} & {\bar{0}} & {\bar{0}} & {\bar{1}} & {\bar{0}}\end{bmatrix}$
&${8}$
&${523}$&${T_{1}T_{2}+T_{3}^{2}+T_{4}^{2}}$
&$\arraycolsep=1mm\begin{bmatrix}{1} & {1} & {1} & {1} & {2} \\ {\bar{0}} & {\bar{0}} & {\bar{1}} & {\bar{1}} & {\bar{0}} \\ {\bar{0}} & {\bar{0}} & {\bar{1}} & {\bar{0}} & {\bar{1}}\end{bmatrix}$
&${16}$
\\
\hline
${510}$&${T_{1}^{2}T_{2}+T_{3}^{4}+T_{4}^{2}}$
&$\arraycolsep=1mm\begin{bmatrix}{1} & {2} & {1} & {2} & {2} \\ {\bar{0}} & {\bar{0}} & {\bar{1}} & {\bar{0}} & {\bar{1}} \\ {\bar{0}} & {\bar{0}} & {\bar{1}} & {\bar{1}} & {\bar{0}}\end{bmatrix}$
&${8}$
\\
\cline{1-4}
\end{tabu}\vspace{-2mm}
\end{table}
\end{class-list}

\begin{class-list}
Non-toric, $\QQ$-factorial, Gorenstein, log terminal Fano
threefolds of Picard number one with an effective two-torus action:
Specifying data for divisor class group
$\ZZ\oplus \ZZ/2\ZZ\oplus \ZZ/2\ZZ$ and format $(2,1,1,1,1)$.
$$\text{\centering\setlength{\tabcolsep}{3.5pt}
\renewcommand{\arraystretch}{1.25}\small\tabulinesep=1mm\begin{tabu}{c|c|c|c||c|c|c|c}
\hline
ID&relations&gd-matrix&$-\mathcal{K}^3$&ID&relations&gd-matrix&$-\mathcal{K}^3$
\\
\hline
${473}$&
$\begin{array}{r}{T_{1}T_{2}+T_{3}^{4}+T_{4}^{2}},\\ {\lambda T_{3}^{4}+T_{4}^{2}+T_{5}^{2}}\end{array}$
&$\arraycolsep=1mm\begin{bmatrix}{2} & {2} & {1} & {2} & {2} & {1} \\ {\bar{0}} & {\bar{0}} & {\bar{1}} & {\bar{0}} & {\bar{1}} & {\bar{0}} \\ {\bar{0}} & {\bar{0}} & {\bar{0}} & {\bar{1}} & {\bar{1}} & {\bar{0}}\end{bmatrix}$&${2}$
&
${508}$&
$\begin{array}{r}{T_{1}T_{2}+T_{3}^{2}+T_{4}^{2}},\\ {\lambda T_{3}^{2}+T_{4}^{2}+T_{5}^{2}}\end{array}$
&$\arraycolsep=1mm\begin{bmatrix}{1} & {1} & {1} & {1} & {1} & {1} \\ {\bar{1}} & {\bar{1}} & {\bar{1}} & {\bar{1}} & {\bar{0}} & {\bar{0}} \\ {\bar{0}} & {\bar{0}} & {\bar{0}} & {\bar{1}} & {\bar{1}} & {\bar{0}}\end{bmatrix}$&${8}$
\\
\hline
${507}$&
$\begin{array}{r}{T_{1}T_{2}+T_{3}^{2}+T_{4}^{2}},\\ {\lambda T_{3}^{2}+T_{4}^{2}+T_{5}^{2}}\end{array}$
&$\arraycolsep=1mm\begin{bmatrix}{1} & {1} & {1} & {1} & {1} & {1} \\ {\bar{0}} & {\bar{0}} & {\bar{1}} & {\bar{1}} & {\bar{0}} & {\bar{0}} \\ {\bar{0}} & {\bar{0}} & {\bar{0}} & {\bar{1}} & {\bar{1}} & {\bar{0}}\end{bmatrix}$&${8}$
\\
\cline{1-4}
\end{tabu}\vspace{-2mm}
}$$
\end{class-list}

\begin{class-list}
Non-toric, $\QQ$-factorial, Gorenstein, log terminal Fano
threefolds of Picard number one with an effective two-torus action:
Specifying data for divisor class group $
\ZZ\oplus \ZZ/2\ZZ\oplus \ZZ/2\ZZ$ and format $(1,1,1,2)$.
\begin{table}[h!]\centering\setlength{\tabcolsep}{3.5pt}
\renewcommand{\arraystretch}{1.15}\small\tabulinesep=1mm
\begin{tabu}{c|c|c|c||c|c|c|c}
\hline
ID&relations&gd-matrix&$-\mathcal{K}^3$&ID&relations&gd-matrix&$-\mathcal{K}^3$
\\
\hline
${496}$&${T_{1}^{6}+T_{2}^{2}+T_{3}^{2}}$
&$\arraycolsep=1mm\begin{bmatrix}{1} & {3} & {3} & {3} & {2} \\ {\bar{0}} & {\bar{0}} & {\bar{1}} & {\bar{1}} & {\bar{0}} \\ {\bar{0}} & {\bar{1}} & {\bar{1}} & {\bar{0}} & {\bar{0}}\end{bmatrix}$
&${6}$
&${505}$&${T_{1}^{2}+T_{2}^{2}+T_{3}^{2}}$
&$\arraycolsep=1mm\begin{bmatrix}{2} & {2} & {2} & {1} & {1} \\ {\bar{1}} & {\bar{1}} & {\bar{0}} & {\bar{0}} & {\bar{0}} \\ {\bar{0}} & {\bar{1}} & {\bar{1}} & {\bar{0}} & {\bar{0}}\end{bmatrix}$
&${8}$
\\
\hline
${497}$&${T_{1}^{3}+T_{2}^{2}+T_{3}^{2}}$
&$\arraycolsep=1mm\begin{bmatrix}{2} & {3} & {3} & {3} & {1} \\ {\bar{0}} & {\bar{1}} & {\bar{1}} & {\bar{0}} & {\bar{0}} \\ {\bar{0}} & {\bar{1}} & {\bar{0}} & {\bar{1}} & {\bar{0}}\end{bmatrix}$
&${6}$
&${506}$&${T_{1}^{2}+T_{2}^{2}+T_{3}^{2}}$
&$\arraycolsep=1mm\begin{bmatrix}{2} & {2} & {2} & {1} & {1} \\ {\bar{0}} & {\bar{0}} & {\bar{1}} & {\bar{1}} & {\bar{0}} \\ {\bar{0}} & {\bar{1}} & {\bar{1}} & {\bar{0}} & {\bar{0}}\end{bmatrix}$
&${8}$
\\
\hline
${498}$&${T_{1}^{2}+T_{2}^{2}+T_{3}^{2}}$
&$\arraycolsep=1mm\begin{bmatrix}{3} & {3} & {3} & {2} & {1} \\ {\bar{1}} & {\bar{1}} & {\bar{0}} & {\bar{0}} & {\bar{0}} \\ {\bar{0}} & {\bar{1}} & {\bar{1}} & {\bar{0}} & {\bar{0}}\end{bmatrix}$
&${6}$
&${521}$&${T_{1}^{2}+T_{2}^{2}+T_{3}^{2}}$
&$\arraycolsep=1mm\begin{bmatrix}{1} & {1} & {1} & {1} & {2} \\ {\bar{1}} & {\bar{1}} & {\bar{0}} & {\bar{0}} & {\bar{0}} \\ {\bar{1}} & {\bar{0}} & {\bar{0}} & {\bar{1}} & {\bar{0}}\end{bmatrix}$
&${16}$
\\
\hline
${502}$&${T_{1}^{4}+T_{2}^{2}+T_{3}^{2}}$
&$\arraycolsep=1mm\begin{bmatrix}{1} & {2} & {2} & {1} & {2} \\ {\bar{0}} & {\bar{0}} & {\bar{1}} & {\bar{1}} & {\bar{0}} \\ {\bar{0}} & {\bar{1}} & {\bar{1}} & {\bar{0}} & {\bar{0}}\end{bmatrix}$
&${8}$
&${522}$&${T_{1}^{2}+T_{2}^{2}+T_{3}^{2}}$
&$\arraycolsep=1mm\begin{bmatrix}{1} & {1} & {1} & {2} & {1} \\ {\bar{0}} & {\bar{0}} & {\bar{1}} & {\bar{1}} & {\bar{0}} \\ {\bar{0}} & {\bar{1}} & {\bar{1}} & {\bar{0}} & {\bar{0}}\end{bmatrix}$
&${16}$
\\
\hline
${503}$&
${T_{1}^{4}+T_{2}^{2}+T_{3}^{2}}$
&$\arraycolsep=1mm\begin{bmatrix}{1} & {2} & {2} & {2} & {1} \\ {\bar{1}} & {\bar{0}} & {\bar{1}} & {\bar{0}} & {\bar{0}} \\ {\bar{1}} & {\bar{1}} & {\bar{1}} & {\bar{1}} & {\bar{0}}\end{bmatrix}$
&${8}$
&${524}$&${T_{1}^{2}+T_{2}^{2}+T_{3}^{2}}$
&$\arraycolsep=1mm\begin{bmatrix}{1} & {1} & {1} & {3} & {2} \\ {\bar{1}} & {\bar{1}} & {\bar{0}} & {\bar{0}} & {\bar{0}} \\ {\bar{1}} & {\bar{0}} & {\bar{0}} & {\bar{1}} & {\bar{0}}\end{bmatrix}$
&${18}$
\\
\hline
${504}$&${T_{1}^{4}+T_{2}^{2}+T_{3}^{2}}$
&$\arraycolsep=1mm\begin{bmatrix}{1} & {2} & {2} & {2} & {1} \\ {\bar{0}} & {\bar{0}} & {\bar{1}} & {\bar{1}} & {\bar{0}} \\ {\bar{0}} & {\bar{1}} & {\bar{1}} & {\bar{0}} & {\bar{0}}\end{bmatrix}$
&${8}$
\\
\cline{1-4}
\end{tabu}
\end{table}
\end{class-list}

\begin{class-list}\label{class-list(0,2,2,2)}
Non-toric, $\QQ$-factorial, Gorenstein, log terminal Fano
threefolds of Picard number one with an effective two-torus action:
Specifying data for divisor class group
$\ZZ\oplus \ZZ/2\ZZ\oplus \ZZ/2\ZZ\oplus \ZZ/2\ZZ$ and
format $(2,2,1,0)$.
\begin{table}[h!]\centering\setlength{\tabcolsep}{5.5pt}
\renewcommand{\arraystretch}{1.15}\small\tabulinesep=1mm
\begin{tabu}{c|c|c|c}
\hline
ID&relations&gd-matrix&$-\mathcal{K}^3$
\\
\hline
${533}$&${T_{1}^{2}T_{2}^{2}+T_{3}^{2}T_{4}^{2}+T_{5}^{2}}$
&$\arraycolsep=1mm\begin{bmatrix}{1} & {1} & {1} & {1} & {2} \\ {\bar{0}} & {\bar{0}} & {\bar{1}} & {\bar{1}} & {\bar{0}} \\ {\bar{1}} & {\bar{0}} & {\bar{0}} & {\bar{1}} & {\bar{0}} \\ {\bar{0}} & {\bar{0}} & {\bar{1}} & {\bar{0}} & {\bar{1}}\end{bmatrix}$
&${2}$
\\
\hline
\end{tabu}
\end{table}
\end{class-list}

\begin{class-list}
Non-toric, $\QQ$-factorial, Gorenstein, log terminal Fano
threefolds of Picard number one with an effective two-torus action:
Specifying data for divisor class group
$\ZZ\oplus \ZZ/2\ZZ\oplus \ZZ/2\ZZ\oplus \ZZ/2\ZZ$ and format $(2,1,1,1)$.
\begin{table}[h!]\centering\setlength{\tabcolsep}{3.5pt}
\renewcommand{\arraystretch}{1.15}\small\tabulinesep=1mm
\begin{tabu}{c|c|c|c||c|c|c|c}
\hline
ID&relations&gd-matrix&$-\mathcal{K}^3$&ID&relations&gd-matrix&$-\mathcal{K}^3$
\\
\hline
${532}$&${T_{1}^{2}T_{2}^{2}+T_{3}^{4}+T_{4}^{2}}$
&$\arraycolsep=1mm\begin{bmatrix}{1} & {1} & {1} & {2} & {1} \\ {\bar{1}} & {\bar{1}} & {\bar{0}} & {\bar{0}} & {\bar{0}} \\ {\bar{0}} & {\bar{0}} & {\bar{1}} & {\bar{1}} & {\bar{0}} \\ {\bar{1}} & {\bar{0}} & {\bar{1}} & {\bar{0}} & {\bar{0}}\end{bmatrix}$
&${2}$
&${537}$&${T_{1}^{2}T_{2}^{2}+T_{3}^{2}+T_{4}^{2}}$
&$\arraycolsep=1mm\begin{bmatrix}{1} & {1} & {2} & {2} & {2} \\ {\bar{1}} & {\bar{0}} & {\bar{1}} & {\bar{1}} & {\bar{1}} \\ {\bar{1}} & {\bar{0}} & {\bar{1}} & {\bar{0}} & {\bar{0}} \\ {\bar{0}} & {\bar{0}} & {\bar{1}} & {\bar{1}} & {\bar{0}}\end{bmatrix}$
&${4}$
\\
\hline
\end{tabu}
\end{table}
\end{class-list}

\begin{class-list}
Non-toric, $\QQ$-factorial, Gorenstein, log terminal Fano
threefolds of Picard number one with an effective two-torus action:
Specifying data for divisor class group
$\ZZ\oplus \ZZ/2\ZZ\oplus \ZZ/2\ZZ\oplus \ZZ/2\ZZ$ and format $(2,1,1,1,1)$.
\begin{table}[h!]\centering\setlength{\tabcolsep}{5.5pt}
\renewcommand{\arraystretch}{1.1}\small\tabulinesep=1mm
\begin{tabu}{c|c|c|c}
\hline
ID&relations&gd-matrix&$-\mathcal{K}^3$
\\
\hline
${536}$&
$\begin{array}{r}{T_{1}T_{2}+T_{3}^{2}+T_{4}^{2}},\\ {\lambda T_{3}^{2}+T_{4}^{2}+T_{5}^{2}}\end{array}$
&$\arraycolsep=1mm\begin{bmatrix}{1} & {1} & {1} & {1} & {1} & {1} \\ {\bar{1}} & {\bar{1}} & {\bar{0}} & {\bar{0}} & {\bar{0}} & {\bar{0}} \\ {\bar{0}} & {\bar{0}} & {\bar{1}} & {\bar{1}} & {\bar{0}} & {\bar{0}} \\ {\bar{1}} & {\bar{1}} & {\bar{0}} & {\bar{1}} & {\bar{1}} & {\bar{0}}\end{bmatrix}$
&${4}$
\\
\hline
\end{tabu}\vspace{-2mm}
\end{table}
\end{class-list}

\begin{class-list}
Non-toric, $\QQ$-factorial, Gorenstein, log terminal Fano
threefolds of Picard number one with an effective two-torus action:
Specifying data for divisor class group
$\ZZ\oplus \ZZ/2\ZZ\oplus \ZZ/2\ZZ\oplus \ZZ/2\ZZ$ and format $(1,1,1,2)$.
\begin{table}[h!]\centering\setlength{\tabcolsep}{3.5pt}
\renewcommand{\arraystretch}{1.1}\small\tabulinesep=1mm
\begin{tabu}{c|c|c|c}
\hline
ID&relations&gd-matrix&$-\mathcal{K}^3$
\\
\hline
${534}$&${T_{1}^{4}+T_{2}^{2}+T_{3}^{2}}$
&$\arraycolsep=1mm\begin{bmatrix}{1} & {2} & {2} & {2} & {1} \\ {\bar{1}} & {\bar{0}} & {\bar{0}} & {\bar{1}} & {\bar{0}} \\ {\bar{1}} & {\bar{0}} & {\bar{1}} & {\bar{0}} & {\bar{0}} \\ {\bar{1}} & {\bar{1}} & {\bar{1}} & {\bar{1}} & {\bar{0}}\end{bmatrix}$
&${4}$
\\
\hline
${535}$&${T_{1}^{2}+T_{2}^{2}+T_{3}^{2}}$
&$\arraycolsep=1mm\begin{bmatrix}{2} & {2} & {2} & {1} & {1} \\ {\bar{1}} & {\bar{1}} & {\bar{1}} & {\bar{1}} & {\bar{0}} \\ {\bar{1}} & {\bar{1}} & {\bar{0}} & {\bar{0}} & {\bar{0}} \\ {\bar{1}} & {\bar{0}} & {\bar{0}} & {\bar{1}} & {\bar{0}}\end{bmatrix}$
&${4}$
\\
\hline
${538}$&${T_{1}^{2}+T_{2}^{2}+T_{3}^{2}}$
&$\arraycolsep=1mm\begin{bmatrix}{1} & {1} & {1} & {2} & {1} \\ {\bar{1}} & {\bar{1}} & {\bar{1}} & {\bar{1}} & {\bar{0}} \\ {\bar{1}} & {\bar{1}} & {\bar{0}} & {\bar{0}} & {\bar{0}} \\ {\bar{1}} & {\bar{0}} & {\bar{0}} & {\bar{1}} & {\bar{0}}\end{bmatrix}$
&${8}$
\\
\hline
\end{tabu}\vspace{-2mm}
\end{table}
\end{class-list}

\begin{class-list}\label{class-list(0,2,4)}
Non-toric, $\QQ$-factorial, Gorenstein, log terminal Fano
threefolds of Picard number one with an effective two-torus action:
Specifying data for divisor class group
$\ZZ\oplus \ZZ/2\ZZ\oplus \ZZ/4\ZZ$ and format $(2,2,1,0)$.
\begin{table}[h!]\centering\setlength{\tabcolsep}{5.5pt}
\renewcommand{\arraystretch}{1.1}\small\tabulinesep=1mm
\begin{tabu}{c|c|c|c}
\hline
ID&relations&gd-matrix&$-\mathcal{K}^3$
\\
\hline
${529}$&${T_{1}^{2}T_{2}^{2}+T_{3}T_{4}+T_{5}^{2}}$
&$\arraycolsep=1mm\begin{bmatrix}{1} & {1} & {2} & {2} & {2} \\ {\bar{0}} & {\bar{0}} & {\bar{1}} & {\bar{1}} & {\bar{0}} \\ {\bar{0}} & {\bar{3}} & {\bar{2}} & {\bar{0}} & {\bar{1}}\end{bmatrix}$
&${4}$
\\
\hline
\end{tabu}\vspace{-2mm}
\end{table}
\end{class-list}

\begin{class-list}
Non-toric, $\QQ$-factorial, Gorenstein, log terminal Fano
threefolds of Picard number one with an effective two-torus action:
Specifying data for divisor class group
$\ZZ\oplus \ZZ/2\ZZ\oplus \ZZ/4\ZZ$ and format $(2,2,1,1,0)$.
\begin{table}[h!]\centering\setlength{\tabcolsep}{5.5pt}
\renewcommand{\arraystretch}{1.1}\small\tabulinesep=1mm
\begin{tabu}{c|c|c|c}
\hline
ID&relations&gd-matrix&$-\mathcal{K}^3$
\\
\hline
${528}$&$\begin{array}{r}{T_{1}T_{2}+T_{3}T_{4}+T_{5}^{2}},\\ {\lambda T_{3}T_{4}+T_{5}^{2}+T_{6}^{2}}\end{array}$
&$\arraycolsep=1mm\begin{bmatrix}{1} & {1} & {1} & {1} & {1} & {1} \\ {\bar{1}} & {\bar{1}} & {\bar{0}} & {\bar{0}} & {\bar{0}} & {\bar{0}} \\ {\bar{3}} & {\bar{1}} & {\bar{3}} & {\bar{1}} & {\bar{2}} & {\bar{0}}\end{bmatrix}$&${4}$
\\
\hline
\end{tabu}\vspace{-2mm}
\end{table}
\end{class-list}

\begin{class-list}
Non-toric, $\QQ$-factorial, Gorenstein, log terminal Fano
threefolds of Picard number one with an effective two-torus action:
Specifying data for divisor class group
$\ZZ\oplus \ZZ/2\ZZ\oplus \ZZ/4\ZZ$ and format $(3,1,1,0)$.
\begin{table}[h!]\centering\setlength{\tabcolsep}{5.5pt}
\renewcommand{\arraystretch}{1.1}\small\tabulinesep=1mm
\begin{tabu}{c|c|c|c}
\hline
ID&relations&gd-matrix&$-\mathcal{K}^3$
\\
\hline
${527}$&${T_{1}T_{2}T_{3}+T_{4}^{2}+T_{5}^{2}}$
&$\arraycolsep=1mm\begin{bmatrix}{1} & {1} & {2} & {2} & {2} \\ {\bar{1}} & {\bar{0}} & {\bar{1}} & {\bar{1}} & {\bar{1}} \\ {\bar{2}} & {\bar{0}} & {\bar{0}} & {\bar{3}} & {\bar{1}}\end{bmatrix}$
&${4}$
\\
\hline
\end{tabu}
\end{table}
\end{class-list}

\begin{class-list}
Non-toric, $\QQ$-factorial, Gorenstein, log terminal Fano
threefolds of Picard number one with an effective two-torus action:
Specifying data for divisor class group
$\ZZ\oplus \ZZ/2\ZZ\oplus \ZZ/4\ZZ$ and format $(2,1,1,1)$.
\begin{table}[h!]\centering\setlength{\tabcolsep}{3.5pt}
\renewcommand{\arraystretch}{1.15}\small\tabulinesep=1mm
\begin{tabu}{c|c|c|c||c|c|c|c}
\hline
ID&relations&gd-matrix&$-\mathcal{K}^3$&ID&relations&gd-matrix&$-\mathcal{K}^3$
\\
\hline
${526}$&${T_{1}T_{2}+T_{3}^{4}+T_{4}^{4}}$
&$\arraycolsep=1mm\begin{bmatrix}{2} & {2} & {1} & {1} & {2} \\ {\bar{1}} & {\bar{1}} & {\bar{0}} & {\bar{0}} & {\bar{0}} \\ {\bar{3}} & {\bar{1}} & {\bar{3}} & {\bar{0}} & {\bar{3}}\end{bmatrix}$
&${4}$
&${530}$&${T_{1}T_{2}+T_{3}^{2}+T_{4}^{2}}$
&$\arraycolsep=1mm\begin{bmatrix}{1} & {1} & {1} & {1} & {2} \\ {\bar{0}} & {\bar{0}} & {\bar{0}} & {\bar{1}} & {\bar{1}} \\ {\bar{2}} & {\bar{0}} & {\bar{1}} & {\bar{3}} & {\bar{0}}\end{bmatrix}$
&${8}$
\\
\hline
\end{tabu}
\end{table}
\end{class-list}

\begin{class-list}\label{class-list(0,2,6)}
Non-toric, $\QQ$-factorial, Gorenstein, log terminal Fano
threefolds of Picard number one with an effective two-torus action:
Specifying data for divisor class group
$\ZZ\oplus \ZZ/2\ZZ\oplus \ZZ/6\ZZ$ and format $(3,1,1,0)$.
\begin{table}[h!]\centering\setlength{\tabcolsep}{5.5pt}
\renewcommand{\arraystretch}{1.1}\small\tabulinesep=1mm
\begin{tabu}{c|c|c|c}
\hline
ID&relations&gd-matrix&$-\mathcal{K}^3$
\\
\hline
${531}$&${T_{1}T_{2}T_{3}+T_{4}^{3}+T_{5}^{3}}$
&$\arraycolsep=1mm\begin{bmatrix}{1} & {1} & {1} & {1} & {1} \\ {\bar{1}} & {\bar{0}} & {\bar{1}} & {\bar{0}} & {\bar{0}} \\ {\bar{5}} & {\bar{5}} & {\bar{2}} & {\bar{4}} & {\bar{0}}\end{bmatrix}$
&${2}$
\\
\hline
\end{tabu}
\end{table}
\end{class-list}

\begin{class-list}\label{class-list(0,3)}
Non-toric, $\QQ$-factorial, Gorenstein, log terminal Fano
threefolds of Picard number one with an effective two-torus action:
Specifying data for divisor class group
$\ZZ\oplus \ZZ/3\ZZ$ and format $(2,2,1,0)$.
\begin{table}[h!]\centering\setlength{\tabcolsep}{3.5pt}
\renewcommand{\arraystretch}{1.15}\small\tabulinesep=1mm
\begin{tabu}{c|c|c|c||c|c|c|c}
\hline
ID&relations&gd-matrix&$-\mathcal{K}^3$&ID&relations&gd-matrix&$-\mathcal{K}^3$
\\
\hline
${424}$&${T_{1}^{3}T_{2}+T_{3}^{3}T_{4}^{3}+T_{5}^{2}}$
&$\arraycolsep=1mm\begin{bmatrix}{1} & {3} & {1} & {1} & {3} \\ {\bar{0}} & {\bar{0}} & {\bar{2}} & {\bar{1}} & {\bar{0}}\end{bmatrix}$&${6}$
&${434}$&${T_{1}T_{2}+T_{3}T_{4}+T_{5}^{2}}$
&$\arraycolsep=1mm\begin{bmatrix}{1} & {3} & {2} & {2} & {2} \\ {\bar{2}} & {\bar{0}} & {\bar{2}} & {\bar{0}} & {\bar{1}}\end{bmatrix}$&${12}$
\\
\hline
${425}$&${T_{1}^{9}T_{2}^{3}+T_{3}T_{4}+T_{5}^{3}}$
&$\arraycolsep=1mm\begin{bmatrix}{1} & {1} & {6} & {6} & {4} \\ {\bar{0}} & {\bar{1}} & {\bar{0}} & {\bar{0}} & {\bar{2}}\end{bmatrix}$&${6}$
&${437}$&${T_{1}^{7}T_{2}+T_{3}T_{4}+T_{5}^{2}}$
&$\arraycolsep=1mm\begin{bmatrix}{1} & {3} & {1} & {9} & {5} \\ {\bar{0}} & {\bar{1}} & {\bar{1}} & {\bar{0}} & {\bar{2}}\end{bmatrix}$&${18}$
\\
\hline
${426}$&${T_{1}^{3}T_{2}^{3}+T_{3}T_{4}+T_{5}^{3}}$
&$\arraycolsep=1mm\begin{bmatrix}{1} & {2} & {3} & {6} & {3} \\ {\bar{0}} & {\bar{1}} & {\bar{0}} & {\bar{0}} & {\bar{2}}\end{bmatrix}$&${6}$
&${438}$&${T_{1}T_{2}+T_{3}T_{4}+T_{5}^{2}}$
&$\arraycolsep=1mm\begin{bmatrix}{1} & {1} & {1} & {1} & {1} \\ {\bar{0}} & {\bar{0}} & {\bar{1}} & {\bar{2}} & {\bar{0}}\end{bmatrix}$&${18}$
\\
\hline
${432}$&${T_{1}^{10}T_{2}+T_{3}T_{4}+T_{5}^{2}}$
&$\arraycolsep=1mm\begin{bmatrix}{1} & {4} & {2} & {12} & {7} \\ {\bar{0}} & {\bar{1}} & {\bar{1}} & {\bar{0}} & {\bar{2}}\end{bmatrix}$&${12}$
&${439}$&${T_{1}T_{2}+T_{3}T_{4}+T_{5}^{2}}$
&$\arraycolsep=1mm\begin{bmatrix}{1} & {1} & {1} & {1} & {1} \\ {\bar{2}} & {\bar{1}} & {\bar{2}} & {\bar{1}} & {\bar{0}}\end{bmatrix}$&${18}$
\\
\hline
${433}$&${T_{1}^{7}T_{2}+T_{3}T_{4}+T_{5}^{2}}$
&$\arraycolsep=1mm\begin{bmatrix}{1} & {1} & {2} & {6} & {4} \\ {\bar{0}} & {\bar{1}} & {\bar{1}} & {\bar{0}} & {\bar{2}}\end{bmatrix}$&${12}$
\\
\cline{1-4}
\end{tabu}
\end{table}
\end{class-list}

\begin{class-list}
Non-toric, $\QQ$-factorial, Gorenstein, log terminal Fano
threefolds of Picard number one with an effective two-torus action:
Specifying data for divisor class group $\ZZ\oplus \ZZ/3\ZZ$
and format $(2,2,1,1,0)$.
\begin{table}[h!]\centering\setlength{\tabcolsep}{5.5pt}
\renewcommand{\arraystretch}{1.35}\small\tabulinesep=1mm
\begin{tabu}{c|c|c|c}
\hline
ID&relations&gd-matrix&$-\mathcal{K}^3$
\\
\hline
${423}$&$\begin{array}{r}{T_{1}T_{2}+T_{3}T_{4}+T_{5}^{3}},\\ {\lambda T_{3}T_{4}+T_{5}^{3}+T_{6}^{3}}\end{array}$
&$\arraycolsep=1mm\begin{bmatrix}{2} & {1} & {1} & {2} & {1} & {1} \\ {\bar{2}} & {\bar{1}} & {\bar{1}} & {\bar{2}} & {\bar{2}} & {\bar{0}}\end{bmatrix}$&${6}$
\\
\hline
\end{tabu}\vspace{2mm}
\end{table}
\end{class-list}

\begin{class-list}
Non-toric, $\QQ$-factorial, Gorenstein, log terminal Fano
threefolds of Picard number one with an effective two-torus action:
Specifying data for divisor class group $\ZZ\oplus \ZZ/3\ZZ$
and format $(3,1,1,0)$.
\begin{table}[h!]\centering\setlength{\tabcolsep}{3.5pt}
\renewcommand{\arraystretch}{1.35}\small\tabulinesep=1mm
\begin{tabu}{c|c|c|c||c|c|c|c}
\hline
ID&relations&gd-matrix&$-\mathcal{K}^3$&ID&relations&gd-matrix&$-\mathcal{K}^3$
\\
\hline
${414}$&${T_{1}^{2}T_{2}^{2}T_{3}+T_{4}^{3}+T_{5}^{3}}$
&$\arraycolsep=1mm\begin{bmatrix}{1} & {1} & {2} & {2} & {2} \\ {\bar{0}} & {\bar{0}} & {\bar{0}} & {\bar{1}} & {\bar{2}}\end{bmatrix}$&${2}$
&${422}$&${T_{1}T_{2}T_{3}+T_{4}^{3}+T_{5}^{3}}$
&$\arraycolsep=1mm\begin{bmatrix}{1} & {2} & {6} & {3} & {3} \\ {\bar{0}} & {\bar{0}} & {\bar{0}} & {\bar{1}} & {\bar{2}}\end{bmatrix}$&${6}$
\\
\hline
${415}$&${T_{1}^{2}T_{2}T_{3}+T_{4}^{3}+T_{5}^{3}}$
&$\arraycolsep=1mm\begin{bmatrix}{2} & {1} & {1} & {2} & {2} \\ {\bar{1}} & {\bar{2}} & {\bar{2}} & {\bar{2}} & {\bar{0}}\end{bmatrix}$&${2}$
&${429}$&${T_{1}T_{2}T_{3}+T_{4}^{3}+T_{5}^{3}}$
&$\arraycolsep=1mm\begin{bmatrix}{1} & {1} & {4} & {2} & {2} \\ {\bar{0}} & {\bar{0}} & {\bar{0}} & {\bar{1}} & {\bar{2}}\end{bmatrix}$&${8}$
\\
\hline
${416}$&${T_{1}T_{2}T_{3}+T_{4}^{9}+T_{5}^{3}}$
&$\arraycolsep=1mm\begin{bmatrix}{1} & {4} & {4} & {1} & {3} \\ {\bar{1}} & {\bar{1}} & {\bar{1}} & {\bar{0}} & {\bar{1}}\end{bmatrix}$&${4}$
&${430}$&${T_{1}T_{2}T_{3}+T_{4}^{3}+T_{5}^{3}}$
&$\arraycolsep=1mm\begin{bmatrix}{1} & {1} & {1} & {1} & {1} \\ {\bar{1}} & {\bar{1}} & {\bar{1}} & {\bar{2}} & {\bar{0}}\end{bmatrix}$&${8}$
\\
\hline
${421}$&${T_{1}^{3}T_{2}^{3}T_{3}+T_{4}^{3}+T_{5}^{2}}$
&$\arraycolsep=1mm\begin{bmatrix}{1} & {1} & {6} & {4} & {6} \\ {\bar{1}} & {\bar{0}} & {\bar{0}} & {\bar{2}} & {\bar{0}}\end{bmatrix}$&${6}$
\\
\cline{1-4}
\end{tabu}\vspace{2mm}
\end{table}
\end{class-list}

\begin{class-list}
Non-toric, $\QQ$-factorial, Gorenstein, log terminal Fano
threefolds of Picard number one with an effective two-torus action:
Specifying data for divisor class group $\ZZ\oplus \ZZ/3\ZZ$ and
format $(2,1,1,1)$.
\begin{table}[h!]\centering\setlength{\tabcolsep}{3.5pt}
\renewcommand{\arraystretch}{1.35}\small\tabulinesep=1mm
\begin{tabu}{c|c|c|c||c|c|c|c}
\hline
ID&relations&gd-matrix&$-\mathcal{K}^3$&ID&relations&gd-matrix&$-\mathcal{K}^3$
\\
\hline
${413}$&${T_{1}^{2}T_{2}^{2}+T_{3}^{3}+T_{4}^{3}}$
&$\arraycolsep=1mm\begin{bmatrix}{2} & {1} & {2} & {2} & {1} \\ {\bar{0}} & {\bar{0}} & {\bar{1}} & {\bar{2}} & {\bar{0}}\end{bmatrix}$&${2}$
&${428}$&${T_{1}^{2}T_{2}+T_{3}^{3}+T_{4}^{3}}$
&$\arraycolsep=1mm\begin{bmatrix}{1} & {1} & {1} & {1} & {1} \\ {\bar{0}} & {\bar{0}} & {\bar{1}} & {\bar{2}} & {\bar{0}}\end{bmatrix}$&${8}$
\\
\hline
${418}$&${T_{1}T_{2}+T_{3}^{12}+T_{4}^{3}}$
&$\arraycolsep=1mm\begin{bmatrix}{6} & {6} & {1} & {4} & {1} \\ {\bar{0}} & {\bar{0}} & {\bar{1}} & {\bar{2}} & {\bar{0}}\end{bmatrix}$&${6}$
&${431}$&${T_{1}T_{2}+T_{3}^{12}+T_{4}^{3}}$
&$\arraycolsep=1mm\begin{bmatrix}{2} & {10} & {1} & {4} & {5} \\ {\bar{1}} & {\bar{2}} & {\bar{0}} & {\bar{1}} & {\bar{1}}\end{bmatrix}$&${10}$
\\
\hline
${419}$&${T_{1}T_{2}+T_{3}^{9}+T_{4}^{3}}$
&$\arraycolsep=1mm\begin{bmatrix}{3} & {6} & {1} & {3} & {2} \\ {\bar{0}} & {\bar{0}} & {\bar{0}} & {\bar{1}} & {\bar{2}}\end{bmatrix}$&${6}$
&${435}$&${T_{1}T_{2}+T_{3}^{9}+T_{4}^{3}}$
&$\arraycolsep=1mm\begin{bmatrix}{1} & {8} & {1} & {3} & {4} \\ {\bar{1}} & {\bar{2}} & {\bar{0}} & {\bar{1}} & {\bar{1}}\end{bmatrix}$&${16}$
\\
\hline
${420}$&${T_{1}T_{2}+T_{3}^{3}+T_{4}^{3}}$
&$\arraycolsep=1mm\begin{bmatrix}{3} & {3} & {2} & {2} & {2} \\ {\bar{0}} & {\bar{0}} & {\bar{2}} & {\bar{0}} & {\bar{1}}\end{bmatrix}$&${6}$
&${436}$&${T_{1}T_{2}+T_{3}^{3}+T_{4}^{3}}$
&$\arraycolsep=1mm\begin{bmatrix}{1} & {2} & {1} & {1} & {2} \\ {\bar{1}} & {\bar{2}} & {\bar{2}} & {\bar{0}} & {\bar{2}}\end{bmatrix}$&${16}$
\\
\hline
${427}$&${T_{1}^{2}T_{2}+T_{3}^{3}+T_{4}^{3}}$
&$\arraycolsep=1mm\begin{bmatrix}{1} & {4} & {2} & {2} & {1} \\ {\bar{0}} & {\bar{0}} & {\bar{1}} & {\bar{2}} & {\bar{0}}\end{bmatrix}$&${8}$
\\
\cline{1-4}
\end{tabu}\vspace{2mm}
\end{table}
\end{class-list}

\begin{class-list}
Non-toric, $\QQ$-factorial, Gorenstein, log terminal Fano
threefolds of Picard number one with an effective two-torus action:
Specifying data for divisor class group $\ZZ\oplus \ZZ/3\ZZ$ and
format $(1,1,1,2)$.
\begin{table}[h!]\centering\setlength{\tabcolsep}{5.5pt}
\renewcommand{\arraystretch}{1.35}\small\tabulinesep=1mm
\begin{tabu}{c|c|c|c}
\hline
ID&relations&gd-matrix&$-\mathcal{K}^3$
\\
\hline
${417}$&${T_{1}^{3}+T_{2}^{3}+T_{3}^{2}}$
&$\arraycolsep=1mm\begin{bmatrix}{2} & {2} & {3} & {3} & {2} \\ {\bar{1}} & {\bar{0}} & {\bar{0}} & {\bar{0}} & {\bar{2}}\end{bmatrix}$&${6}$
\\
\hline
\end{tabu}\vspace{4mm}
\end{table}
\end{class-list}

\begin{class-list}\label{class-list(0,3,3)}
Non-toric, $\QQ$-factorial, Gorenstein, log terminal Fano
threefolds of Picard number one with an effective two-torus action:
Specifying data for divisor class group
$\ZZ\oplus \ZZ/3\ZZ\oplus \ZZ/3\ZZ$ and format $(2,2,1,0)$.
\begin{table}[h!]\centering\setlength{\tabcolsep}{5.5pt}
\renewcommand{\arraystretch}{1.35}\small\tabulinesep=1mm
\begin{tabu}{c|c|c|c}
\hline
ID&relations&gd-matrix&$-\mathcal{K}^3$
\\
\hline
${525}$&${T_{1}T_{2}+T_{3}T_{4}+T_{5}^{2}}$
&$\arraycolsep=1mm\begin{bmatrix}{1} & {1} & {1} & {1} & {1} \\ {\bar{2}} & {\bar{1}} & {\bar{2}} & {\bar{1}} & {\bar{0}} \\ {\bar{2}} & {\bar{1}} & {\bar{0}} & {\bar{0}} & {\bar{0}}\end{bmatrix}$
&${6}$
\\
\hline
\end{tabu}\vspace{4mm}
\end{table}
\end{class-list}

\begin{class-list}\label{class-list(0,4)}
Non-toric, $\QQ$-factorial, Gorenstein, log terminal Fano
threefolds of Picard number one with an effective two-torus action:
Specifying data for divisor class group $\ZZ\oplus \ZZ/4\ZZ$ and
format $(2,2,1,0)$.
\begin{table}[h!]\centering\setlength{\tabcolsep}{3.5pt}
\renewcommand{\arraystretch}{1.35}\small\tabulinesep=1mm
\begin{tabu}{c|c|c|c||c|c|c|c}
\hline
ID&relations&gd-matrix&$-\mathcal{K}^3$&ID&relations&gd-matrix&$-\mathcal{K}^3$
\\
\hline
${452}$&${T_{1}^{5}T_{2}+T_{3}^{2}T_{4}+T_{5}^{2}}$
&$\arraycolsep=1mm\begin{bmatrix}{1} & {1} & {1} & {4} & {3} \\ {\bar{0}} & {\bar{2}} & {\bar{3}} & {\bar{0}} & {\bar{1}}\end{bmatrix}$&${8}$
&${455}$&${T_{1}^{10}T_{2}^{2}+T_{3}T_{4}+T_{5}^{2}}$
&$\arraycolsep=1mm\begin{bmatrix}{1} & {1} & {4} & {8} & {6} \\ {\bar{0}} & {\bar{1}} & {\bar{2}} & {\bar{0}} & {\bar{3}}\end{bmatrix}$&${8}$
\\
\hline
${453}$&${T_{1}^{3}T_{2}^{3}+T_{3}^{2}T_{4}+T_{5}^{2}}$
&$\arraycolsep=1mm\begin{bmatrix}{1} & {1} & {1} & {4} & {3} \\ {\bar{2}} & {\bar{0}} & {\bar{3}} & {\bar{0}} & {\bar{1}}\end{bmatrix}$&${8}$
&${456}$&${T_{1}^{6}T_{2}^{2}+T_{3}T_{4}+T_{5}^{2}}$
&$\arraycolsep=1mm\begin{bmatrix}{1} & {2} & {2} & {8} & {5} \\ {\bar{0}} & {\bar{1}} & {\bar{2}} & {\bar{0}} & {\bar{3}}\end{bmatrix}$&${8}$
\\
\hline
${454}$&${T_{1}^{2}T_{2}+T_{3}^{2}T_{4}+T_{5}^{2}}$
&$\arraycolsep=1mm\begin{bmatrix}{1} & {2} & {1} & {2} & {2} \\ {\bar{0}} & {\bar{2}} & {\bar{3}} & {\bar{0}} & {\bar{1}}\end{bmatrix}$&${8}$
&${457}$&${T_{1}^{2}T_{2}^{2}+T_{3}T_{4}+T_{5}^{2}}$
&$\arraycolsep=1mm\begin{bmatrix}{1} & {1} & {2} & {2} & {2} \\ {\bar{0}} & {\bar{3}} & {\bar{0}} & {\bar{2}} & {\bar{1}}\end{bmatrix}$&${8}$
\\
\hline
\end{tabu}\vspace{4mm}
\end{table}
\end{class-list}

\begin{class-list}
Non-toric, $\QQ$-factorial, Gorenstein, log terminal Fano
threefolds of Picard number one with an effective two-torus action:
Specifying data for divisor class group $\ZZ\oplus \ZZ/4\ZZ$ and
format $(2,2,1,1,0)$.
\begin{table}[h!]\centering\setlength{\tabcolsep}{5.5pt}
\renewcommand{\arraystretch}{1.35}\small\tabulinesep=1mm
\begin{tabu}{c|c|c|c}
\hline
ID&relations&gd-matrix&$-\mathcal{K}^3$
\\
\hline
${451}$&$\begin{array}{r}{T_{1}T_{2}+T_{3}T_{4}+T_{5}^{2}},\\ {\lambda T_{3}T_{4}+T_{5}^{2}+T_{6}^{2}}\end{array}$
&$\arraycolsep=1mm\begin{bmatrix}{1} & {1} & {1} & {1} & {1} & {1} \\ {\bar{1}} & {\bar{3}} & {\bar{3}} & {\bar{1}} & {\bar{2}} & {\bar{0}}\end{bmatrix}$&${8}$
\\
\hline
\end{tabu}\vspace{4mm}
\end{table}
\end{class-list}

\newpage

\begin{class-list}
Non-toric, $\QQ$-factorial, Gorenstein, log terminal Fano
threefolds of Picard number one with an effective two-torus action:
Specifying data for divisor class group $\ZZ\oplus \ZZ/4\ZZ$ and
format $(3,1,1,0)$.
\begin{table}[h!]\centering\setlength{\tabcolsep}{3.5pt}
\renewcommand{\arraystretch}{1.15}\small\tabulinesep=1mm
\begin{tabu}{c|c|c|c||c|c|c|c}
\hline
ID&relations&gd-matrix&$-\mathcal{K}^3$&ID&relations&gd-matrix&$-\mathcal{K}^3$
\\
\hline
${443}$&${T_{1}T_{2}T_{3}+T_{4}^{10}+T_{5}^{2}}$
&$\arraycolsep=1mm\begin{bmatrix}{1} & {3} & {6} & {1} & {5} \\ {\bar{3}} & {\bar{3}} & {\bar{2}} & {\bar{0}} & {\bar{2}}\end{bmatrix}$&${6}$
&${449}$&${T_{1}T_{2}T_{3}+T_{4}^{2}+T_{5}^{2}}$
&$\arraycolsep=1mm\begin{bmatrix}{1} & {1} & {2} & {2} & {2} \\ {\bar{3}} & {\bar{3}} & {\bar{0}} & {\bar{3}} & {\bar{1}}\end{bmatrix}$&${8}$
\\
\hline
${444}$&${T_{1}T_{2}T_{3}+T_{4}^{2}+T_{5}^{2}}$
&$\arraycolsep=1mm\begin{bmatrix}{1} & {2} & {3} & {3} & {3} \\ {\bar{3}} & {\bar{2}} & {\bar{3}} & {\bar{2}} & {\bar{0}}\end{bmatrix}$&${6}$
&${450}$&${T_{1}T_{2}T_{3}+T_{4}^{2}+T_{5}^{2}}$
&$\arraycolsep=1mm\begin{bmatrix}{1} & {1} & {2} & {2} & {2} \\ {\bar{3}} & {\bar{1}} & {\bar{2}} & {\bar{3}} & {\bar{1}}\end{bmatrix}$&${8}$
\\
\hline
${448}$&${T_{1}^{2}T_{2}^{2}T_{3}+T_{4}^{3}+T_{5}^{2}}$
&$\arraycolsep=1mm\begin{bmatrix}{1} & {1} & {8} & {4} & {6} \\ {\bar{1}} & {\bar{0}} & {\bar{0}} & {\bar{2}} & {\bar{3}}\end{bmatrix}$&${8}$
\\
\cline{1-4}
\end{tabu}
\end{table}
\end{class-list}

\begin{class-list}
Non-toric, $\QQ$-factorial, Gorenstein, log terminal Fano
threefolds of Picard number one with an effective two-torus action:
Specifying data for divisor class group $\ZZ\oplus \ZZ/4\ZZ$ and
format $(2,1,1,1)$.
\begin{table}[h!]\centering\setlength{\tabcolsep}{3.5pt}
\renewcommand{\arraystretch}{1.15}\small\tabulinesep=1mm
\begin{tabu}{c|c|c|c||c|c|c|c}
\hline
ID&relations&gd-matrix&$-\mathcal{K}^3$&ID&relations&gd-matrix&$-\mathcal{K}^3$
\\
\hline
${440}$&${T_{1}T_{2}+T_{3}^{8}+T_{4}^{4}}$
&$\arraycolsep=1mm\begin{bmatrix}{4} & {4} & {1} & {2} & {1} \\ {\bar{0}} & {\bar{0}} & {\bar{1}} & {\bar{3}} & {\bar{0}}\end{bmatrix}$&${4}$
&${446}$&${T_{1}^{2}T_{2}+T_{3}^{2}+T_{4}^{2}}$
&$\arraycolsep=1mm\begin{bmatrix}{1} & {2} & {2} & {2} & {1} \\ {\bar{0}} & {\bar{2}} & {\bar{1}} & {\bar{3}} & {\bar{0}}\end{bmatrix}$&${8}$
\\
\hline
${441}$&${T_{1}^{3}T_{2}+T_{3}^{2}+T_{4}^{2}}$
&$\arraycolsep=1mm\begin{bmatrix}{1} & {3} & {3} & {3} & {2} \\ {\bar{0}} & {\bar{2}} & {\bar{1}} & {\bar{3}} & {\bar{0}}\end{bmatrix}$&${6}$
&${447}$&${T_{1}T_{2}+T_{3}^{4}+T_{4}^{4}}$
&$\arraycolsep=1mm\begin{bmatrix}{2} & {2} & {1} & {1} & {2} \\ {\bar{1}} & {\bar{3}} & {\bar{3}} & {\bar{0}} & {\bar{3}}\end{bmatrix}$&${8}$
\\
\hline
${442}$&${T_{1}T_{2}+T_{3}^{8}+T_{4}^{4}}$
&$\arraycolsep=1mm\begin{bmatrix}{2} & {6} & {1} & {2} & {3} \\ {\bar{2}} & {\bar{2}} & {\bar{0}} & {\bar{1}} & {\bar{1}}\end{bmatrix}$&${6}$
&${458}$&${T_{1}T_{2}+T_{3}^{2}+T_{4}^{2}}$
&$\arraycolsep=1mm\begin{bmatrix}{1} & {1} & {1} & {1} & {2} \\ {\bar{2}} & {\bar{0}} & {\bar{3}} & {\bar{1}} & {\bar{0}}\end{bmatrix}$&${16}$
\\
\hline
${445}$&${T_{1}^{3}T_{2}+T_{3}^{2}+T_{4}^{2}}$
&$\arraycolsep=1mm\begin{bmatrix}{1} & {1} & {2} & {2} & {2} \\ {\bar{0}} & {\bar{2}} & {\bar{1}} & {\bar{3}} & {\bar{0}}\end{bmatrix}$&${8}$
\\
\cline{1-4}
\end{tabu}
\end{table}
\end{class-list}

\begin{class-list}\label{class-list(0,5)}
Non-toric, $\QQ$-factorial, Gorenstein, log terminal Fano
threefolds of Picard number one with an effective two-torus action:
Specifying data for divisor class group $\ZZ\oplus \ZZ/5\ZZ$ and
format $(2,2,1,0)$.
\begin{table}[h!]\centering\setlength{\tabcolsep}{5.5pt}
\renewcommand{\arraystretch}{1.15}\small\tabulinesep=1mm
\begin{tabu}{c|c|c|c}
\hline
ID&relations&gd-matrix&$-\mathcal{K}^3$
\\
\hline
${462}$&${T_{1}^{3}T_{2}^{3}+T_{3}T_{4}+T_{5}^{2}}$
&$\arraycolsep=1mm\begin{bmatrix}{1} & {1} & {1} & {5} & {3} \\ {\bar{3}} & {\bar{0}} & {\bar{4}} & {\bar{0}} & {\bar{2}}\end{bmatrix}$&${10}$
\\
\hline
\end{tabu}
\end{table}
\end{class-list}

\begin{class-list}
Non-toric, $\QQ$-factorial, Gorenstein, log terminal Fano
threefolds of Picard number one with an effective two-torus action:
Specifying data for divisor class group $\ZZ\oplus \ZZ/5\ZZ$ and
format $(3,1,1,0)$.
\begin{table}[h!]\centering\setlength{\tabcolsep}{3.5pt}
\renewcommand{\arraystretch}{1.15}\small\tabulinesep=1mm
\begin{tabu}{c|c|c|c||c|c|c|c}
\hline
ID&relations&gd-matrix&$-\mathcal{K}^3$&ID&relations&gd-matrix&$-\mathcal{K}^3$
\\
\hline
${459}$&${T_{1}T_{2}T_{3}+T_{4}^{5}+T_{5}^{5}}$
&$\arraycolsep=1mm\begin{bmatrix}{1} & {2} & {2} & {1} & {1} \\ {\bar{2}} & {\bar{4}} & {\bar{4}} & {\bar{4}} & {\bar{0}}\end{bmatrix}$&${2}$
&${461}$&${T_{1}T_{2}T_{3}+T_{4}^{3}+T_{5}^{2}}$
&$\arraycolsep=1mm\begin{bmatrix}{1} & {1} & {10} & {4} & {6} \\ {\bar{1}} & {\bar{0}} & {\bar{0}} & {\bar{2}} & {\bar{3}}\end{bmatrix}$&${10}$
\\
\hline
\end{tabu}
\end{table}
\end{class-list}

\begin{class-list}
Non-toric, $\QQ$-factorial, Gorenstein, log terminal Fano
threefolds of Picard number one with an effective two-torus action:
Specifying data for divisor class group $\ZZ\oplus \ZZ/5\ZZ$ and
format $(2,1,1,1)$.
\begin{table}[h!]\centering\setlength{\tabcolsep}{5.5pt}
\renewcommand{\arraystretch}{1.1}\small\tabulinesep=1mm
\begin{tabu}{c|c|c|c}
\hline
ID&relations&gd-matrix&$-\mathcal{K}^3$
\\
\hline
${460}$&${T_{1}T_{2}+T_{3}^{5}+T_{4}^{5}}$
&$\arraycolsep=1mm\begin{bmatrix}{1} & {4} & {1} & {1} & {2} \\ {\bar{2}} & {\bar{3}} & {\bar{4}} & {\bar{0}} & {\bar{4}}\end{bmatrix}$&${8}$
\\
\hline
\end{tabu}\vspace{-2mm}
\end{table}
\end{class-list}

\begin{class-list}\label{class-list(0,6)}
Non-toric, $\QQ$-factorial, Gorenstein, log terminal Fano
threefolds of Picard number one with an effective two-torus action:
Specifying data for divisor class group $\ZZ\oplus \ZZ/6\ZZ$ and
format $(2,2,1,0)$.
\begin{table}[h!]\centering\setlength{\tabcolsep}{5.5pt}
\renewcommand{\arraystretch}{1.1}\small\tabulinesep=1mm
\begin{tabu}{c|c|c|c}
\hline
ID&relations&gd-matrix&$-\mathcal{K}^3$
\\
\hline
${469}$&${T_{1}^{4}T_{2}^{4}+T_{3}T_{4}+T_{5}^{2}}$
&$\arraycolsep=1mm\begin{bmatrix}{1} & {1} & {2} & {6} & {4} \\ {\bar{1}} & {\bar{0}} & {\bar{4}} & {\bar{0}} & {\bar{5}}\end{bmatrix}$&${6}$
\\
\hline
\end{tabu}\vspace{-2mm}
\end{table}
\end{class-list}

\begin{class-list}
Non-toric, $\QQ$-factorial, Gorenstein, log terminal Fano
threefolds of Picard number one with an effective two-torus action:
Specifying data for divisor class group $\ZZ\oplus \ZZ/6\ZZ$ and
format $(3,1,1,0)$.
\begin{table}[h!]\centering\setlength{\tabcolsep}{3.5pt}
\renewcommand{\arraystretch}{1.1}\small\tabulinesep=1mm
\begin{tabu}{c|c|c|c||c|c|c|c}
\hline
ID&relations&gd-matrix&$-\mathcal{K}^3$&ID&relations&gd-matrix&$-\mathcal{K}^3$
\\
\hline
${465}$&${T_{1}T_{2}T_{3}+T_{4}^{3}+T_{5}^{3}}$
&$\arraycolsep=1mm\begin{bmatrix}{1} & {1} & {4} & {2} & {2} \\ {\bar{3}} & {\bar{0}} & {\bar{0}} & {\bar{5}} & {\bar{1}}\end{bmatrix}$&${4}$
&${468}$&${T_{1}T_{2}T_{3}+T_{4}^{4}+T_{5}^{2}}$
&$\arraycolsep=1mm\begin{bmatrix}{1} & {1} & {6} & {2} & {4} \\ {\bar{2}} & {\bar{0}} & {\bar{0}} & {\bar{5}} & {\bar{1}}\end{bmatrix}$&${6}$
\\
\hline
${466}$&${T_{1}T_{2}T_{3}+T_{4}^{3}+T_{5}^{3}}$
&$\arraycolsep=1mm\begin{bmatrix}{1} & {1} & {1} & {1} & {1} \\ {\bar{2}} & {\bar{5}} & {\bar{5}} & {\bar{4}} & {\bar{0}}\end{bmatrix}$&${4}$
\\
\cline{1-4}
\end{tabu}\vspace{-2mm}
\end{table}
\end{class-list}

\begin{class-list}
Non-toric, $\QQ$-factorial, Gorenstein, log terminal Fano
threefolds of Picard number one with an effective two-torus action:
Specifying data for divisor class group $\ZZ\oplus \ZZ/6\ZZ$ and
format $(2,1,1,1)$.
\begin{table}[h!]\centering\setlength{\tabcolsep}{3.5pt}
\renewcommand{\arraystretch}{1.1}\small\tabulinesep=1mm
\begin{tabu}{c|c|c|c||c|c|c|c}
\hline
ID&relations&gd-matrix&$-\mathcal{K}^3$&ID&relations&gd-matrix&$-\mathcal{K}^3$
\\
\hline
${463}$&${T_{1}^{2}T_{2}+T_{3}^{3}+T_{4}^{3}}$
&$\arraycolsep=1mm\begin{bmatrix}{1} & {1} & {1} & {1} & {1} \\ {\bar{0}} & {\bar{3}} & {\bar{5}} & {\bar{1}} & {\bar{0}}\end{bmatrix}$&${4}$
&${467}$&${T_{1}T_{2}+T_{3}^{4}+T_{4}^{2}}$
&$\arraycolsep=1mm\begin{bmatrix}{2} & {2} & {1} & {2} & {3} \\ {\bar{2}} & {\bar{4}} & {\bar{0}} & {\bar{3}} & {\bar{3}}\end{bmatrix}$&${6}$
\\
\hline
${464}$&${T_{1}T_{2}+T_{3}^{6}+T_{4}^{6}}$
&$\arraycolsep=1mm\begin{bmatrix}{2} & {4} & {1} & {1} & {2} \\ {\bar{2}} & {\bar{4}} & {\bar{5}} & {\bar{0}} & {\bar{5}}\end{bmatrix}$&${4}$
&${470}$&${T_{1}T_{2}+T_{3}^{3}+T_{4}^{3}}$
&$\arraycolsep=1mm\begin{bmatrix}{1} & {2} & {1} & {1} & {2} \\ {\bar{5}} & {\bar{1}} & {\bar{4}} & {\bar{0}} & {\bar{4}}\end{bmatrix}$&${8}$
\\
\hline
\end{tabu}\vspace{-2mm}
\end{table}
\end{class-list}

\begin{class-list}\label{class-list(0,8)}
Non-toric, $\QQ$-factorial, Gorenstein, log terminal Fano
threefolds of Picard number one with an effective two-torus action:
Specifying data for divisor class group $\ZZ\oplus \ZZ/8\ZZ$ and
format $(3,1,1,0)$.
\begin{table}[h!]\centering\setlength{\tabcolsep}{3.5pt}
\renewcommand{\arraystretch}{1.1}\small\tabulinesep=1mm
\begin{tabu}{c|c|c|c||c|c|c|c}
\hline
ID&relations&gd-matrix&$-\mathcal{K}^3$&ID&relations&gd-matrix&$-\mathcal{K}^3$
\\
\hline
${471}$&${T_{1}T_{2}T_{3}+T_{4}^{4}+T_{5}^{4}}$
&$\arraycolsep=1mm\begin{bmatrix}{1} & {1} & {2} & {1} & {1} \\ {\bar{3}} & {\bar{7}} & {\bar{6}} & {\bar{6}} & {\bar{0}}\end{bmatrix}$&${2}$
&${472}$&${T_{1}T_{2}T_{3}+T_{4}^{6}+T_{5}^{2}}$
&$\arraycolsep=1mm\begin{bmatrix}{1} & {1} & {4} & {1} & {3} \\ {\bar{2}} & {\bar{0}} & {\bar{0}} & {\bar{3}} & {\bar{5}}\end{bmatrix}$&${4}$
\\
\hline
\end{tabu}\vspace{-4mm}
\end{table}
\end{class-list}

\section{Hilbert--Poincar\'{e} series}
\label{sec:hilbert-series}

Here we present the Hilbert--Poincar\'{e} series of our
Fano varieties.
Recall that the Hilbert--Poincar\'{e} series of a
finitely generated $\ZZ_{\ge 0}$-graded $\KK$-algebra
$A = \oplus_{k} A_k$ is the formal power series
\[
\HP_A(t) := \sum_{k \ge 0} \dim_\KK(A_k) t^k.
\]
Assume that $f_1,\dots,f_r \in A$ are homogeneous
of degrees $w_1,\dots,w_r$ respectively and generate
$A$ as an algebra.
Then there is a polynomial $q_A \in \ZZ[t]$
such that
\[
\HP_A(t) = \frac{q_A(t)}{\prod_{i=1}^r (1-t^{w_i})}.
\]

Given a Fano variety $X$,
we associate with it the Hilbert--Poincar\'{e}
series $\HP_X(t)$ of its anticanonical ring $A_X$ and
we define the corresponding polynomial $q_X(t)$
with respect to a~minimal system of homogeneous
generators of the anticanonical ring~$A_X$.

\begin{prop}
The following table lists for each possible pair $(g,c)$
of genus and codimension the classification IDs
from Section~$\ref{sec:class-lists}$ of the varieties
$X$ attaining $(g,c)$ and the cancelled presentation
of the associated Hilbert--Poincar\'{e}
series together with its first eight terms.
\begin{table}[h!]\centering\setlength{\tabcolsep}{4.5pt}
\renewcommand{\arraystretch}{1.4}\small\tabulinesep=1mm
\begin{tabu}{c|c|c}\hline
 $(g,c)$ & $\HP_X(t)$ & {\rm IDs}
 \\
 \hline
$(2,1)$\label{T(1,2)}
&$\begin{array}{c}
 \dfrac{1+t^3}{1-t^4} \\[10pt]
 {1+4t+10t^2+21t^3+39t^4+66t^5+104t^6+155t^7+\cdots}
\end{array}$
& \begin{minipage}{4cm} \rm
 2, 9, 10, 11, 246, 251, 252, 253, 254, 459,
 471, 473, 474, 475, 477, 478, 479, 480, 531, 532, 533
 \end{minipage}
\\
\hline
$(2,2)$\label{T(2,2)}
&$\begin{array}{c}
\dfrac{1+t^3}{1-t^4} \\[10pt]
{1+4t+10t^2+21t^3+39t^4+66t^5+104t^6+155t^7+\cdots}
\end{array}$
& \begin{minipage}{4cm} \rm
 1, 3, 4, 5, 6, 7, 8, 247, 248, 249,
 250, 413, 414, 415, 476
 \end{minipage}
\\
\hline
$(3,1)$\label{T(1,3)}
&$\begin{array}{c}
\dfrac{1+t+t^2+t^3}{1-t^4} \\[10pt]
{1+5t+15t^2+35t^3+69t^4+121t^5+195t^6+295t^7+\cdots}
\end{array}$
&\begin{minipage}{4cm} \rm
 12, 416, 440, 472, 529
 \end{minipage}
\\
\hline
$(3,2)$\label{T(2,3)}
&$\begin{array}{c}
 \dfrac{1+t+t^2+t^3}{1-t^4} \\[15pt]
{1+5t+15t^2+35t^3+69t^4+121t^5+195t^6+295t^7+\cdots}
\end{array}$
& \begin{minipage}{4cm} \rm
 13, 14, 15, 16, 17, 18, 255, 256, 257, 258, 259,
 260, 261, 262, 263, 264, 265, 266, 267, 268, 269,
 270, 271, 272, 273, 274, 275, 276, 463, 464, 465,
 466, 481, 482, 483, 484, 485, 486, 487, 488, 489,
 490, 491, 492, 493, 494, 495, 526, 527, 528, 534,
 535, 536, 537
 \end{minipage}
\\
\hline
$(4,2)$\label{T(2,4)}
&$\begin{array}{c}
	 \dfrac{1+2t+2t^2+t^3}{1-t^4} \\[10pt]
	 {1+6t+20t^2+49t^3+99t^4+176t^5+286t^6+435t^7+\cdots}
\end{array}$
& \begin{minipage}{4cm} \rm
 19, 21, 22, 28, 30, 32, 33, 34, 35, 36,
	 37, 280, 282, 287, 290, 417, 418, 419, 420, 421,
	 423, 424, 425, 426, 442, 443, 467, 468, 469, 501,
	 525
 \end{minipage}
\\
\hline
\end{tabu}
\end{table}

\begin{table}[h!]\centering\setlength{\tabcolsep}{4.5pt}
\renewcommand{\arraystretch}{1.3}\small\tabulinesep=1mm
\begin{tabu}{c|c|c}\hline
 $(g,c)$ & $\HP_X(t)$ & {\rm IDs}
 \\
 \hline
	$(4,4)$\label{T(4,4)}
	&
	$\begin{array}{c}
	 \dfrac{1+2t+2t^2+t^3}{1-t^4} \\[3pt]
	 {1+6t+20t^2+49t^3+99t^4+176t^5+286t^6+435t^7+\cdots}
	\end{array}$
	&
	 \begin{minipage}{4cm}
 \rm
 	20, 23, 24, 25, 26, 27, 29, 31, 38, 39,
 	40, 41, 42, 277, 278, 279, 281, 283, 284, 285,
	 286, 288, 289, 291, 292, 293, 294, 295, 296, 422,
	 441, 444, 496, 497, 498, 499, 500
 \end{minipage}
	\\ \hline
	$(5,3)$\label{T(3,5)}
	&
	$\begin{array}{c}
	 \dfrac{1+3t+3t^2+t^3}{1-t^4} \\[3pt]
	 {1+7t+25t^2+63t^3+129t^4+231t^5+377t^6+575t^7+\cdots}
	\end{array}$
	&
	 \begin{minipage}{4cm}
 \rm
 301, 317, 321, 322, 325, 330, 332, 333, 334, 335,
	 448, 451, 452, 453, 454, 455, 456, 457, 503, 516,
	 517, 530
 \end{minipage}
	\\ \hline
	$(5,6)$\label{T(6,5)}
	&
	$\begin{array}{c}
	 \dfrac{1+3t+3t^2+t^3}{1-t^4} \\[15pt]
	 {1+7t+25t^2+63t^3+129t^4+231t^5+377t^6+575t^7+\cdots}
	\end{array}$
	&
	 \begin{minipage}{4cm}
 \rm
 	43, 44, 45, 46, 47, 48, 49, 50, 51, 52,
	 53, 54, 55, 56, 57, 58, 297, 298, 299, 300,
	 302, 303, 304, 305, 306, 307, 308, 309, 310, 311,
	 312, 313, 314, 315, 316, 318, 319, 320, 323, 324,
	 326, 327, 328, 329, 331, 427, 428, 429, 430, 445,
	 446, 447, 449, 450, 460, 470, 502, 504, 505, 506,
	 507, 508, 509, 510, 511, 512, 513, 514, 515, 518,
	 538
 \end{minipage}
	\\ \hline
	$(6,4)$\label{T(4,6)}
	&
	$\begin{array}{c}
	 \dfrac{1+4t+4t^2+t^3}{1-t^4} \\[3pt]
	 {1+8t+30t^2+77t^3+159t^4+286t^5+468t^6+715t^7+\cdots}
	\end{array}$
	&
	 \begin{minipage}{4cm}
 \rm
 	59, 60, 61, 62, 63, 64, 65, 66, 67, 68,
	 69, 70, 336, 337, 431, 461, 462
 \end{minipage}
	\\ \hline
	$(7,5)$\label{T(5,7)}
	&
	$\begin{array}{c}
	 \dfrac{1+5t+5t^2+t^3}{1-t^4} \\[15pt]
	 {1+9t+35t^2+91t^3+189t^4+341t^5+559t^6+855t^7+\cdots}
	\end{array}$
	&
	 \begin{minipage}{4cm}
 \rm
 	71, 72, 73, 74, 75, 76, 77, 78, 79, 80,
	 81, 82, 83, 84, 85, 86, 87, 88, 89, 90,
	 91, 92, 93, 94, 95, 96, 97, 98, 99, 100,
	 338, 339, 340, 341, 342, 343, 344, 345, 346, 347,
	 348, 349, 350, 351, 352, 353, 354, 355, 356, 357,
	 358, 359, 360, 361, 362, 363, 364, 432, 433, 434,
	 519, 520
 \end{minipage}
\\ \hline
	$(9,7)$\label{T(7,9)}
	&
	$\begin{array}{c}
	 \dfrac{1+7t+7t^2+t^3}{1-t^4} \\[15pt]
	 {1+11t+45t^2+119t^3+249t^4+451t^5+741t^6+1135t^7+\cdots}
	\end{array}$
	&
	 \begin{minipage}{4cm}
 \rm
 	101, 102, 103, 104, 105, 106, 107, 108, 109, 110,
	 111, 112, 113, 114, 115, 116, 117, 365, 366, 367,
	 368, 369, 370, 371, 372, 373, 374, 375, 376, 377,
	 378, 379, 380, 381, 382, 383, 384, 385, 386, 387,
	 388, 389, 390, 391, 392, 435, 436, 458, 521, 522,
	 523
 \end{minipage}
 \\
\hline
\end{tabu}
\end{table}

\begin{table}[h!]\centering\setlength{\tabcolsep}{3.5pt}
\renewcommand{\arraystretch}{1.3}\small\tabulinesep=1mm
\begin{tabu}{c|c|c}\hline
 $(g,c)$ & $\HP_X(t)$ & {\rm IDs}
 \\
 \hline
	$(10,8)$\label{T(8,10)}
	&
	$\begin{array}{c}
	 \dfrac{1+8t+8t^2+t^3}{1-t^4} \\[10pt]
	 {1+12t+50t^2+133t^3+279t^4+506t^5+832t^6+1275t^7+\cdots}
	\end{array}$
	&
	 \begin{minipage}{4cm}
 \rm
 	118, 119, 120, 121, 122, 123, 124, 125, 126, 127,
	 128, 129, 130, 131, 132, 133, 134, 135, 136, 137,
	 138, 139, 140, 141, 142, 143, 393, 394, 395, 396,
	 397, 398, 399, 400, 401, 402, 403, 437, 438, 439,
	 524
 \end{minipage}
	\\ \hline
	$(11,9)$\label{T(9,11)}
	&
	$\begin{array}{c}
	 \dfrac{1+9t+9t^2+t^3}{1-t^4} \\[3pt]
	 {1+13t+55t^2+147t^3+309t^4+561t^5+923t^6+1415t^7+\cdots}
	\end{array}$
	&
	 \begin{minipage}{4cm}
 \rm
 	144, 145, 146, 147, 148, 149, 150, 151, 404, 405
 \end{minipage}
	\\ \hline
	$(12,10)$\label{T(10,12)}
	&
	$\begin{array}{c}
	 \dfrac{1+10t+10t^2+t^3}{1-t^4} \\[3pt]
	 {1+14t+60t^2+161t^3+339t^4+616t^5+1014t^6+1555t^7+\cdots}
	\end{array}$
	&
	\begin{minipage}{4cm}
 \rm
 	152
 \end{minipage}
	\\ \hline
	$(13,11)$\label{T(11,13)}
	&
	$\begin{array}{c}
	 \dfrac{1+11t+11t^2+t^3}{1-t^4} \\[3pt]
	 {1+15t+65t^2+175t^3+369t^4+671t^5+1105t^6+1695t^7+\cdots}
	\end{array}$
	&
	 \begin{minipage}{4cm}
 \rm
 	153, 154, 155, 156, 157, 158, 159, 160, 161, 162,
	 163, 164, 165, 166, 167, 168, 169, 170, 171, 172,
	 173, 174, 175, 176, 177, 178, 179, 180, 181, 182,
	 183, 184, 185, 186, 187, 406, 407, 408, 409, 410,
	 411
 \end{minipage}
\\ \hline
	$(16,14)$\label{T(14,16)}
	&
	$\begin{array}{c}
	 \dfrac{1\!+\!14t\!+\!14t^2\!+\!t^3}{1-t^4} \\[3pt]
	 {1\!+\!18t\!+\!80t^2\!+\!217t^3\!+\!459t^4\!+\!836t^5\!+\!1378t^6\!+\!2115t^7\!+\!\cdots}
	\end{array}$
	&
	 \begin{minipage}{4cm}
 \rm
 	188, 189, 190, 191, 192, 193, 194, 195, 196
 \end{minipage}
	\\ \hline
	$(17,15)$\label{T(15,17)}
	&
	$\begin{array}{c}
	 \dfrac{1\!+\!15t\!+\!15t^2\!+\!t^3}{1-t^4} \\[3pt]
	 {1\!+\!19t\!+\!85t^2\!+\!231t^3\!+\!489t^4\!+\!891t^5\!+\!1469t^6\!+\!2255t^7\!+\!\cdots}
	\end{array}$
	&
	 \begin{minipage}{4cm}
 \rm
 	197, 198, 199, 200, 201, 202, 203, 204, 205, 206,
	 207, 208, 209, 210, 211, 212, 213, 214, 215, 412
 \end{minipage}
	\\ \hline
	$(19,17)$\label{T(17,19)}
	&
	$\begin{array}{c}
	 \dfrac{1\!+\!17t\!+\!17t^2\!+\!t^3}{1-t^4} \\[3pt]
	 {1\!+\!21t\!+\!95t^2\!+\!259t^3\!+\!549t^4\!+\!1001t^5\!+\!1651t^6\!+\!2535t^7\!+\!\cdots}
	\end{array}$
	&
	 \begin{minipage}{4cm}
 \rm
 	216, 217, 218, 219, 220, 221, 222, 223, 224, 225,
	 226, 227
 \end{minipage}
	\\ \hline
	$(21,19)$\label{T(19,21)}
	&
	$\begin{array}{c}
	 \dfrac{1\!+\!19t\!+\!19t^2\!+\!t^3}{1-t^4} \\[3pt]
	 {1\!+\!23t\!+\!105t^2\!+\!287t^3\!+\!609t^4\!+\!1111t^5\!+\!1833t^6\!+\!2815t^7\!+\!\cdots}
	\end{array}$
	&
	 \begin{minipage}{4cm}
 \rm
 	228, 229, 230, 231
 \end{minipage}
	\\ \hline
	$(22,20)$\label{T(20,22)}
	&
	$\begin{array}{c}
	 \dfrac{1\!+\!20t\!+\!20t^2\!+\!t^3}{1-t^4} \\[3pt]
	 {1\!+\!24t\!+\!110t^2\!+\!301t^3\!+\!639t^4\!+\!1166t^5\!+\!1924t^6\!+\!2955t^7\!+\!\cdots}
	\end{array}$
	&
	 \begin{minipage}{4cm}
 \rm
 	232, 233, 234
 \end{minipage}
	\\ \hline
	$(25,23)$\label{T(23,25)}
	&
	$\begin{array}{c}
	 \dfrac{1\!+\!23t\!+\!23t^2\!+\!t^3}{1-t^4} \\[3pt]
	 {1\!+\!27t\!+\!125t^2\!+\!343t^3\!+\!729t^4\!+\!1331t^5\!+\!2197t^6\!+\!3375t^7\!+\!\cdots}
	\end{array}$
	&
	 \begin{minipage}{4cm}
 \rm
 	235, 236, 237, 238
 \end{minipage}
	\\ \hline
\end{tabu}\vspace{4mm}
\end{table}

\begin{table}[h!]\centering\setlength{\tabcolsep}{3.5pt}
\renewcommand{\arraystretch}{1.45}\small\tabulinesep=1mm
\begin{tabu}{c|c|c}\hline
 $(g,c)$ & $\HP_X(t)$ & {\rm IDs}
 \\
 \hline
$(26,24)$\label{T(24,26)}
	&
	$\begin{array}{c}
	 \dfrac{1\!+\!24t\!+\!24t^2\!+\!t^3}{1-t^4} \\[3pt]
	 {1\!+\!28t\!+\!130t^2\!+\!357t^3\!+\!759t^4\!+\!1386t^5\!+\!2288t^6\!+\!3515t^7\!+\!\cdots}
	\end{array}$
	&
	 \begin{minipage}{4cm}
 \rm
 	239, 240, 241, 242
 \end{minipage}
	\\
\hline	
$(28,26)$\label{T(26,28)}
	&
	$\begin{array}{c}
	 \dfrac{1\!+\!26t\!+\!26t^2\!+\!t^3}{1-t^4} \\[3pt]
	 {1\!+\!30t\!+\!140t^2\!+\!385t^3\!+\!819t^4\!+\!1496t^5\!+\!2470t^6\!+\!3795t^7\!+\!\cdots}
	\end{array}$
	&
	 \begin{minipage}{4cm}
 \rm
 	243, 244, 245
 \end{minipage}
	\\
\hline
\end{tabu}
\end{table}
\end{prop}

\begin{proof}
Observe that the anticanonical ring is the Veronese subalgebra
of the Cox ring associated to the subgroup generated by the
anticanonical class.
Thus, we can use the Cox ring data from the classification lists
in Section~\ref{sec:class-lists} to compute a minimal system of
generators and the associated relations.
This provides us in particular with genus and codimension.
Moreover, it allows us to compute the Hilbert--Poincar\'{e}
series; we used the computer algebra system Singular.
\end{proof}

\begin{coro}
The Hilbert--Poincar\'e series of a non-toric, $\QQ$-factorial,
log-terminal, Gorenstein, Fano threefold $X$ of Picard number
one with an effective action of a two-dimensional torus only depends on the genus $g$
of $X$ and can be explicitly written down as
\[
\HP_X(t) = \frac{1 + (g-2)(t+t^2) + t^3}{(1-t)^4}.
\]
\end{coro}

\subsection*{Acknowledgements}
We would like to thank Gavin Brown and
Al Kasprzyk for interesting and helpful
discussions.
Moreover, we are grateful to the referees
for many very helpful remarks and suggestions.

\pdfbookmark[1]{References}{ref}
\LastPageEnding

\end{document}